\documentclass[letterpaper,11pt]{report}
\usepackage{rotating}
\usepackage[spanish]{babel}
\usepackage{amssymb,amsmath}
\usepackage{epsfig}
\usepackage{fancyhdr}

\usepackage[pdftex,pdftitle={anexoc},colorlinks=true,breaklinks=true,linkcolor=blue]{hyperref}

\newcommand{\RR}{I\!\!R}

\newcommand{\supp}{\textrm{sop }}

\newtheorem{Th}{Teorema}
\newtheorem{Le}{Lema}

\newcommand{\dem}{{\bf Dem: }}
\newcommand{\qedo}{\hfill$\square$}

\begin{document}
\pagenumbering{roman}
\section*{}
\thispagestyle{empty}

\vskip-2.9cm
\begin{figure}[h]
\begin{center}
\includegraphics[height=1in]{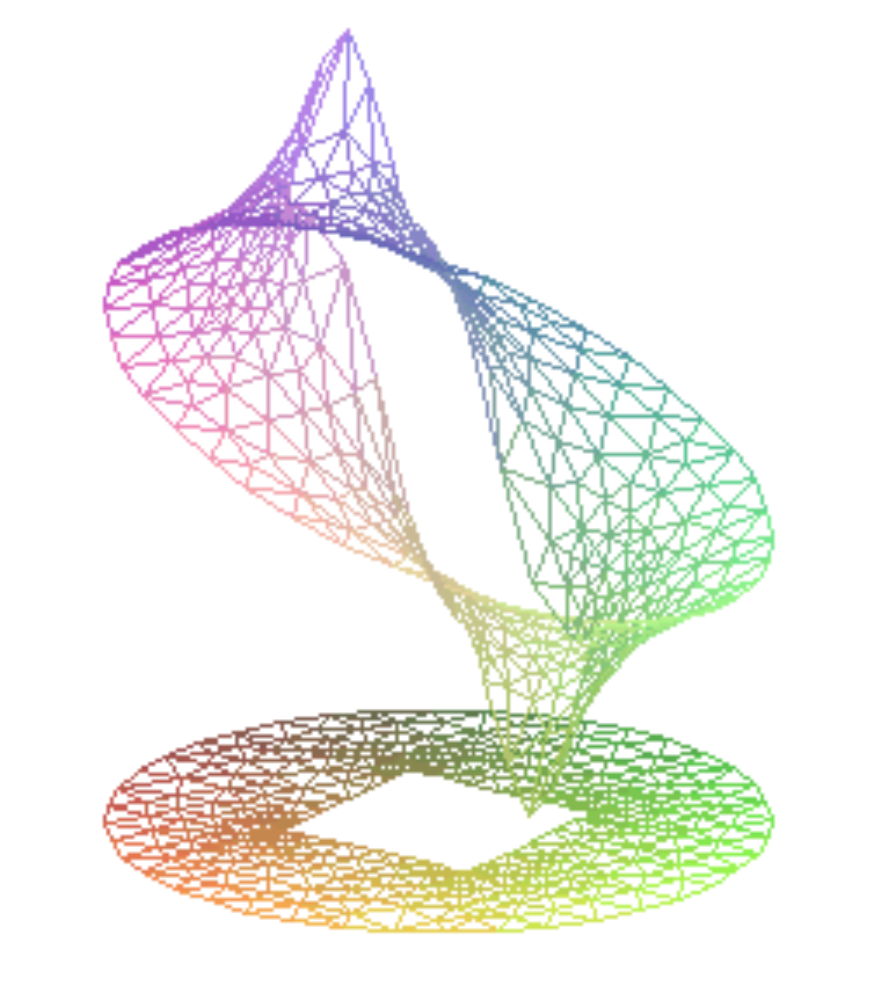}
\end{center}
\end{figure}
\vskip-2.5cm
\begin{center}
{\bf{\large FACULTAD DE CIENCIAS F\'ISICAS Y MATEM\'ATICAS\\ \hspace{.3cm} DEPARTAMENTO DE INGENIER\'IA MATEM\'ATICA}}
\vskip-1cm
         \vskip110pt

         {\Huge M\'etodos de Multiresoluci\'on\\ y su Aplicaci\'on a un Modelo de Ingenier{\'i}a }

        \vskip60pt

         {\sl Tesis para optar al t{\'i}tulo de Ingeniero Matem\'atico}

         \vskip30pt
         {\bf Ricardo Esteban Ruiz Baier}
         \vskip60pt
         \bf{Marzo 2005}
\end{center}
\vskip70pt
\hrule
\vskip25pt

\noindent \makebox(2,2)[l]{\includegraphics[height=0.5in]{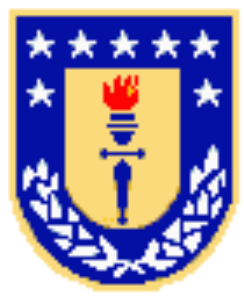}}
\vskip-0.5cm
{\bf{ \qquad \quad UNIVERSIDAD DE CONCEPCI\'ON, CONCEPCI\'ON, CHILE}}

\normalsize
\newpage


\pagestyle{fancy}
\renewcommand{\chaptermark}[1]{\markboth{Cap{\'i}tulo \thechapter:\
 #1}{}}
\renewcommand{\sectionmark}[1]{\markright{\thesection\ #1}}
\fancyhead[LE,RO]{\thepage}
\fancyhead[LE]{\sl{\rightmark}}
\fancyhead[RE]{\sl{\leftmark}}
\tableofcontents
\clearpage{\pagestyle{empty}\cleardoublepage}
\chapter*{Agradecimientos}
A Dios en primer lugar por su constante cuidado y demostraci\'on de fidelidad. 

Deseo agradecer al Departamento de Ingenier{\'i}a Matem\'atica de la Universidad de Concepci\'on por el apoyo brindado. En este marco, deseo agradecer especialmente a mi profesor gu{\'i}a Dr. Mauricio Sep\'ulveda por su apoyo infinito en este \'ultimo per{\'i}odo de estudios, y a los asesores de tesis Dr. Freddy Paiva, Dr. Raimund B\"urger y Dra. Alice Kozakevicius por sus valios{\'i}simos consejos y comentarios. Agradezco tambi\'en el apoyo de FONDECYT mediante su proyecto 1030718 por el financiamiento de este trabajo.

Tambi\'en quiero mencionar (en orden aleatorio) a varios profesores que fueron siempre un apoyo para m{\'i}: Dr. Alberto Foppiano, Dr. Gabriel Gatica, Dr. Manuel Campos, Dr. Gabriel Barrenechea, Dr. Fabi\'an Flores y Dr. Rodolfo Rodr{\'i}guez. 

Agradezco de igual forma a mis amigos: David, Hern\'an, Bollis, Hans, Hiram, D.E.J.A., Alianza, Xime, Ro \& Fa; y compa\~neros: R.R.R., ing-mat. promociones 1997-1998 y otras. No s\'e si habr\'an contribuido al t\'ermino o al retraso de esta tesis, pero es mi deber mencionarlos. 

Estoy en especial agradecido de mi polola quien me motiv\'o constantemente a esforzarme y dedicarme a este trabajo. Espero poder volver a agradecerle en mi pr\'oxima tesis.
 
Finalmente quiero agradecer profundamente a mi familia, que me ha apoyado y alentado para que concluya \'este y todos mis proyectos de vida. Cari\~nos para pap\'a, mam\'a, Fe\~na, Alexis, Joel, Welo, Abuela, Vero, Cutys.  

\clearpage{\pagestyle{empty}\cleardoublepage}
\pagenumbering{arabic}
\chapter{Introducci\'on}
El principal objetivo de este trabajo es presentar una adaptaci\'on de los m\'etodos de vol\'umenes finitos utilizados en la resoluci\'on de problemas provenientes de los procesos de sedimentaci\'on de suspensiones floculadas (o sedimentaci\'on con compresi\'on). Esta adaptaci\'on est\'a basada en la utilizaci\'on de t\'ecnicas de multiresoluci\'on, originalmente ideadas para rebajar el costo computacional en la resoluci\'on num\'erica de leyes de conservaci\'on hiperb\'olicas, en conjunto con esquemas de alta resoluci\'on. 

Se introducir\'an los m\'etodos utilizados para la resoluci\'on num\'erica de leyes de conservaci\'on y ecuaciones parab\'olicas y la importancia del algoritmo de multiresoluci\'on en la aplicaci\'on de estos m\'etodos.

\subsection*{Leyes de conservaci\'on hiperb\'olicas} 
Los sistemas de leyes de conservaci\'on son modelos matem\'aticos para situaciones f{\'i}sicas en que la cantidad total de la variable no var{\'i}a con respecto al tiempo. En este tipo de situaciones, la cantidad de una variable f{\'i}sica contenida en una regi\'on acotada del espacio s\'olo puede variar debido al flujo de la variable a trav\'es de la frontera de dicha regi\'on. Esto puede traducirse en una formulaci\'on integral que, bajo ciertas hip\'otesis de regularidad, se convierte en un sistema de ecuaciones en derivadas parciales.  Si se toma el caso unidimensional (en espacio), las ecuaciones correspondientes son de la forma
\begin{equation}\label{edp1}
\partial_tu(x,t)+\partial_xf(u(x,t))=0,
\end{equation}
donde $u:\RR\times\RR\to\RR^m$ es el vector de variables conservadas o variables de estado, y $f:\RR^m\to\RR^m$ es el vector de flujos. En problemas de din\'amica de fluidos, estas variables son densidad, momento y energ{\'i}a.

La ecuaci\'on (\ref{edp1}) est\'a provista de condiciones iniciales y posiblemente condiciones de frontera en el dominio espacial acotado. 

Un ejemplo cl\'asico para ilustrar el comportamiento de las soluciones en leyes de conservaci\'on, es el problema de Riemann en un tubo de shock: din\'amica de los gases. Se tiene un tubo lleno con gas, inicialmente dividido en dos secciones por una membrana. El gas tiene densidad y presi\'on, en reposo, m\'as alta en una mitad del tubo que en la otra. En el tiempo $t=0$ se rompe la membrana y el gas fluye. Si se supone que el flujo es uniforme a lo largo del tubo, la variaci\'on se produce s\'olo en una direcci\'on y pueden aplicarse las ecuaciones de Euler en una dimensi\'on.

La estructura de la soluci\'on del problema de Riemann implica tres ondas distintas que separan regiones en las que las variables son constantes. La \emph{onda de choque} se propaga hacia la regi\'on de m\'as baja presi\'on; a trav\'es de esta onda, la densidad y la presi\'on asumen valores m\'as altos y todas las variables son discontinuas. Luego aparece una \emph{discontinuidad de contacto}, a trav\'es de la cual la densidad es discontinua, pero las dem\'as variables son constantes. La tercera es la \emph{onda de rarefacci\'on} (recibe este nombre debido a que la densidad del gas decrece cuando esta onda pasa a trav\'es de \'el) que se mueve en direcci\'on contraria a las otras dos y tiene una estructura diferente: todas las variables son continuas y presentan una suave transici\'on \cite{GR}. 


\subsection*{Ecuaciones parab\'olicas}
Se quiere estudiar un problema de valores iniciales para una ecuaci\'on parab\'olica. Para ello, para $(x,t)\in \Omega\times[0,\infty[$, consid\'erese la ecuaci\'on 
\begin{equation*}
\begin{split}
\partial_t u(x,t)+\partial_xF(u(x,t),\partial_x u(x,t))&=S(u),\\
u(x,0)=u_0(x)&
\end{split}
\end{equation*}
donde ahora el flujo $F$ incluye a la derivada de $u$ y este se define por un operador diferencial con difusividad constante $\nu>0$, es decir, 
\begin{equation*}
F(u(x,t),\partial_x u(x,t)):=f(u)-\nu\partial_x u(x,t).
\end{equation*}
Se tienen versiones lineales y no lineales. Para la ecuaci\'on de convecci\'on-difusi\'on unidimensional, se tiene 
\begin{eqnarray*}
f(u)&=&cu,\\
S(u)&=&0,
\end{eqnarray*}
con $c>0$. Este tipo de ecuaciones es de gran utilidad, por ejemplo, para calcular el transporte de sedimentos as{\'i} como el transporte de constituyentes en estudios de calidad de agua \cite{Carrillo}.
 
En el caso de la ecuaci\'on viscosa de Burgers unidimensional, se tiene
\begin{eqnarray*}
f(u)&=&\frac{u^2}{2},\\
S(u)&=&0,
\end{eqnarray*}
Esta ecuaci\'on es un modelo sencillo para la propagaci\'on de fluidos, tomando en cuenta que existe viscosidad constante en el fluido.

Para la ecuaci\'on de reacci\'on-difusi\'on ($\alpha>0,\beta>0$),
\begin{eqnarray*}
f(u)&=&0,\\
S(u)&=&\frac{\beta^2}{2}(1-u)\exp\frac{\beta(1-u)}{\alpha(1-u)-1}.
\end{eqnarray*}
Esta ecuaci\'on conduce al modelo unidimensional de la propagaci\'on de llama premezclada  \cite{Sch}, donde las difusividades de masa y calor son iguales. La funci\'on $u$ representa la temperatura adimensional, que var{\'i}a entre 0 y 1, y la masa parcial de gas sin quemar es representada por $1-u$.

\subsection*{Ecuaciones parab\'olicas fuertemente degeneradas}
Consid\'erese una ecuaci\'on parab\'olica de la forma
\begin{equation}\label{introentro}
\partial_tu+\partial_xf(u)=\partial_{xx}^2A(u),
\end{equation}
con $(x,t)\in ]0,1[\times[0,T[$ y 
\begin{equation*}
A(u):=\int_0^ua(s)ds,\quad a(u)\geqslant 0.
\end{equation*}
En general se permite que $a(u)$ sea cero en incluso un intervalo $[0,u_c]$, en el cual la ecuaci\'on es de naturaleza hiperb\'olica, y $a(u)$ es discontinua en $u=u_c$. Dada la forma degenerada de $a(u)$ y la naturaleza generalmente no lineal de $f(u)$, las soluciones de la ecuaci\'on son generalmente discontinuas y es necesario considerar soluciones entr\'opicas. 

Una ecuaci\'on de convecci\'on-difusi\'on fuertemente degenerada, con una funci\'on de flujo no necesariamente convexa que depende del tiempo, asociada a ciertas condiciones iniciales y de frontera como (\ref{introentro}) se considera como el modelo cl\'asico para los procesos de sedimentaci\'on-consolidaci\'on. La sedimentaci\'on es, a grandes rasgos, un proceso en que part{\'i}culas o agregados son separados de un l{\'i}quido bajo la acci\'on de la fuerza de gravedad. Este es probablemente el m\'etodo industrial a gran escala m\'as importante  utilizado en qu{\'i}mica y miner{\'i}a \cite{PFP}. En soluciones relativamente diluidas, las part{\'i}culas no se comportan en forma discreta sino que tienden a agregarse unas a otras durante el proceso de sedimentaci\'on. Conforme se produce la floculaci\'on, la masa de part{\'i}culas va aumentando, y se deposita a mayor velocidad. La medida en que se desarrolle el fen\'omeno de floculaci\'on depende de la posibilidad de contacto entre las diferentes part{\'i}culas, que a su vez es funci\'on de la carga de superficie, de la profundidad del tanque, del gradiente de velocidad del sistema, de la concentraci\'on de part{\'i}culas y de los tama\~nos de las mismas. El efecto de estas variables sobre el proceso s\'olo puede determinarse mediante ensayos de sedimentaci\'on. Esto hace que sea de gran utilidad en la modelaci\'on de estos fen\'omenos, la teor{\'i}a de problemas inversos (ver \cite{BBCS,Mauricio} entre otros).

Desde hace ya varios a\~nos se ha estado trabajado con mucho \'enfasis en mejorar los fundamentos de los modelos existentes para este tipo de procesos. Grandes avances se deben al trabajo de B\"urger \emph{et al.} \cite{BBK,BCS,BEK,BEKL,BK,BWC} entre otros. Para una descripci\'on detallada de estos procesos y su modelaci\'on, se recomienda consultar \cite{BBK,Thick}.  

Por las caracter{\'i}sticas de este tipo de ecuaciones, no es posible aplicar ni la teor{\'i}a de ecuaciones estrictamente parab\'olicas, ni la teor{\'i}a establecida de soluciones de entrop{\'i}a de leyes de conservaci\'on \cite{BWC}.

\subsection*{M\'etodo de multiresoluci\'on: Motivaci\'on}
Generalmente, el vector de flujos en una ecuaci\'on hiperb\'olica o parab\'olica, est\'a formado por funciones cuya dependencia de las variables de estado es no lineal y esto hace que no sea posible deducir soluciones exactas para estas ecuaciones. De aqu{\'i} nace la necesidad de dise\~nar m\'etodos num\'ericos que aproximen convenientemente estas soluciones. Este es un problema general que afecta a la mayor parte de las ecuaciones en derivadas parciales no lineales, sin embargo, existen razones para estudiar esta clase particular de sistemas:
\begin{itemize}
\item Muchos problemas pr\'acticos en ingenier{\'i}a y ciencia involucran cantidades que se conservan y conducen a problemas del tipo ley de conservaci\'on.
\item Existen dificultades especiales y espec{\'i}ficas a esta clase de sistemas (por ejemplo la formaci\'on de ondas de choque) que no se observan en otros problemas no lineales y que deben tenerse en cuenta en el dise\~no de m\'etodos num\'ericos que aproximen sus soluciones.
\item Aunque se conocen pocas soluciones exactas, la estructura matem\'atica de las ecuaciones y sus soluciones es cada d{\'i}a m\'as estudiada. Este conocimiento se puede aprovechar para desarrollar m\'etodos adecuados a las caracter{\'i}sticas de estos sistemas y sus soluciones.
\end{itemize}
El hecho de que las soluciones de este tipo de ecuaciones admitan discontinuidades plantea varios problemas, tanto desde el punto de vista matem\'atico como num\'erico. Es evidente que una soluci\'on discontinua no puede satisfacer la ecuaci\'on en derivadas parciales en el sentido cl\'asico. La teor{\'i}a de distribuciones provee de una herramienta matem\'atica muy \'util, pues permite caracterizar las discontinuidades admisibles y definir el concepto de soluci\'on d\'ebil de un problema diferencial.

Sin embargo, la clase de funciones continuas a trozos es demasiado amplia para garantizar unicidad de soluci\'on. Generalmente existen soluciones d\'ebiles con los mismos datos iniciales. Puesto que estas ecuaciones son modelos para situaciones f{\'i}sicas reales (o al menos esa es la motivaci\'on), es obvio que s\'olo una de estas soluciones puede ser aceptable desde el punto de vista f{\'i}sico. El hecho de que existan otras soluciones esp\'ureas es consecuencia de que nuestras ecuaciones son tan s\'olo un modelo que ignora algunos efectos f{\'i}sicos, particularmente en el caso de leyes de conservaci\'on, los efectos difusivos y viscosos. Aunque estos efectos (y otros) pueden ignorarse en la mayor parte del fluido, cerca de las discontinuidades juegan un rol esencial. 

Estas consideraciones conducen a la imposici\'on de determinados criterios basados en consideraciones f{\'i}sicas que permiten aislar la soluci\'on f{\'i}sicamente relevante entre todas las posibles soluciones d\'ebiles. Este tipo de criterios se conocen como \emph{condiciones de entrop{\'i}a} de nuevo por analog{\'i}a con la din\'amica de gases (en este caso, la segunda ley de la Termodin\'amica: La entrop{\'i}a nunca decrece). En particular cuando las mol\'eculas del gas pasan a trav\'es de una onda de choque, su entrop{\'i}a deber\'a aumentar, y esto proporciona el principio f{\'i}sico adecuado para determinar de manera un{\'i}voca la soluci\'on con sentido f{\'i}sico. 

La aproximaci\'on num\'erica de este tipo de soluciones incorpora un nuevo conjunto de problemas. Las discretizaciones de la ecuaci\'on en derivadas parciales mediante diferencias finitas ocasionar\'an problemas si las soluciones que se quieren aproximar son discontinuas. Estos problemas son de dos tipos. En general, los m\'etodos num\'ericos de primer orden incorporan difusi\'on num\'erica; esto facilita la convergencia a la soluci\'on entr\'opica, pero limita la utilidad real de estos m\'etodos.
Los metodos cl\'asicos de orden superior reducen la \emph{viscosidad num\'erica} pero incorporan t\'erminos dispersivos y dan lugar a oscilaciones num\'ericas que pueden desencadenar inestabilidades no lineales o hacer que las aproximaciones num\'ericas no converjan a la soluci\'on f{\'i}sicamente relevante.

Los esquemas num\'ericos dise\~nados para la aproximaci\'on de las soluciones de este tipo de ecuaciones deben poder escribirse en forma conservativa. Esto garantiza que si las aproximaciones num\'ericas convergen, lo hacen a una soluci\'on d\'ebil de la ecuaci\'on (Teorema de Lax-Wendroff).

Si un m\'etodo conservativo satisface adem\'as alg\'un an\'alogo discreto de las condiciones de entrop{\'i}a, el l{\'i}mite de las aproximaciones num\'ericas ser\'a precisamente la soluci\'on relevante desde el punto de vista f{\'i}sico.

Una excelente clase de m\'etodos conservativos para la aproximaci\'on num\'erica de las ecuaciones hiperb\'olicas y parab\'olicas, son los m\'etodos de alto orden de precisi\'on. Estos proporcionan perfiles bien delimitados y sin oscilaciones cerca de las discontinuidades. Un aspecto importante a tener en cuenta de los m\'etodos de alto orden de precisi\'on, es su elevado costo computacional, el cual es a\'un mayor bajo las siguientes condiciones:
\begin{itemize}
\item Sistemas de ecuaciones.
\item M\'as de una dimensi\'on.
\item Un gran n\'umero de puntos en la malla.
\item Extensos per{\'i}odos de simulaci\'on.
\end{itemize}
\subsection*{M\'etodo de multiresoluci\'on: Descripci\'on}
El m\'etodo de multiresoluci\'on es una t\'ecnica destinada (al menos, originalmente) a rebajar el costo computacional asociado a los m\'etodos de alta resoluci\'on. En situaciones est\'andar, el comportamiento de la soluci\'on $w(x,t)$ como funci\'on de $x$ es altamente no uniforme, con fuertes variaciones en regiones  puntuales y un comportamiento suave en la mayor parte del intervalo computacional. La t\'ecnica  de multiresoluci\'on (al menos, en la forma en que ser\'a utilizada en este trabajo) fue dise\~nada originalmente por Harten \cite{Harten1} para ecuaciones hiperb\'olicas y utilizada por Bihari \cite{Bihari} y Roussel \emph{et al.} \cite{Sch} para ecuaciones parab\'olicas. Se desea estudiar la aplicaci\'on del m\'etodo de multiresoluci\'on a los m\'etodos existentes para modelar fen\'omenos de sedimentaci\'on de suspensiones floculadas \cite{Thick}. 

Dado un m\'etodo en forma conservativa y una malla uniforme apropiada para la soluci\'on num\'erica del problema de valores iniciales para una ley de conservaci\'on hiperb\'olica escalar o una ecuaci\'on parab\'olica, el m\'etodo de multiresoluci\'on aproxima la soluci\'on a una tolerancia prescrita de una forma m\'as eficiente, entendiendo por eficiencia una reducci\'on en  el n\'umero de veces que se calcula el flujo num\'erico con el m\'etodo de alta resoluci\'on. Para ello se consideran los valores puntuales o medias en celda de la soluci\'on num\'erica mediante un proceso jer\'arquico de mallas anidadas di\'adicas, en el cual la malla dada es la m\'as fina, y se introduce una representaci\'on que contiene la misma informaci\'on.

La representaci\'on de multiresoluci\'on de la soluci\'on num\'erica est\'a formada por sus valores puntuales en la malla m\'as gruesa y el conjunto de errores por interpolar los valores puntuales de cada nivel de resoluci\'on a partir de los del nivel pr\'oximo m\'as grueso. La compresi\'on de datos es realizada haciendo cero las componentes de la representaci\'on que est\'an por debajo de una tolerancia prescrita, e incluso eliminando de la malla a los puntos cuyos errores son menores a esta tolerancia prescrita; por consiguiente en lugar de calcular la evoluci\'on en tiempo de la soluci\'on num\'erica en la malla dada, se calcula la evoluci\'on de su representaci\'on de multiresoluci\'on comprimida. Como la transformaci\'on entre una funci\'on y su representaci\'on de ondelette es r\'apida, la proposici\'on de efectuar la gran parte de los c\'alculos en la representaci\'on de multiresoluci\'on es factible y atractiva.

La informaci\'on contenida en el an\'alisis de multiresoluci\'on de la soluci\'on num\'erica es utilizada para identificar la localizaci\'on de las discontinuidades en la soluci\'on num\'erica, y dise\~nar m\'etodos que mejoren el c\'alculo del flujo num\'erico. Esta informaci\'on es de gran utilidad al momento de calcular los flujos, pues el procedimiento correspondiente toma en cuenta la regularidad de la funci\'on. Adem\'as, la eficiencia computacional del m\'etodo de multiresoluci\'on est\'a directamente relacionada con la raz\'on de compresi\'on de los datos iniciales, es decir, la soluci\'on num\'erica en la malla m\'as fina \cite{Harten1}. La eficiencia del algoritmo se mide mediante la tasa de compresi\'on y el tiempo de CPU. 

\subsection*{Programa}
Este trabajo se organiza del siguiente modo: En el cap{\'i}tulo 2 se revisar\'an los conceptos b\'asicos necesarios para el an\'alisis de multiresoluci\'on propuesto por Harten \cite{Harten1}. En el cap{\'i}tulo 3, se utiliza este an\'alisis para desarrollar un m\'etodo de alta resoluci\'on en mallas generadas por multiresoluci\'on, dise\~nado por Kozakevicius y Santos \cite{Alice}, el que ser\'a aplicado a leyes de conservaci\'on hiperb\'olicas escalares. Se muestran resultados de los test num\'ericos realizados. En el cap{\'i}tulo 4 se analizan las ecuaciones parab\'olicas escalares y un m\'etodo num\'erico que utiliza la multiresoluci\'on y la alta resoluci\'on (esquemas ENO de segundo orden y esquemas Runge-Kutta de segundo orden) como herramientas principales. Se utiliza una nueva estructura de datos desarrollada por Cohen \emph{et al.} \cite{Cohen}. Se muestran los resultados de los experimentos num\'ericos realizados, coincidentes con los resultados obtenidos por Roussel \emph{et al.} \cite{Sch}. En el cap{\'i}tulo 5 se presentan los supuestos b\'asicos para el problema de la sedimentaci\'on, analizando varios casos test. Se simula un proceso de sedimantaci\'on tipo \emph{Batch} y se muestran resultados obtenidos aplicando m\'etodos de multiresoluci\'on a los esquemas desarrollados por B\"urger \emph{et al.} \cite{BCS,BEK,BEKL,BK,BWC}. Se observa que el m\'etodo de multiresoluci\'on es de gran ayuda para reducir el costo computacional en este tipo de problemas sin afectar la calidad de la soluci\'on. 



\clearpage{\pagestyle{empty}\cleardoublepage}
\chapter{Multiresoluci\'on y compresi\'on de datos}
En este cap{\'i}tulo se presentan los conceptos y definiciones b\'asicas introducidas por Harten \cite{Harten1} para el an\'alisis de multiresoluci\'on. Se presentan adem\'as herramientas adicionales utilizadas por Kozakevicius y Santos \cite{Alice} para el desarrollo de m\'etodos con mallas generadas mediante an\'alisis de multiresoluci\'on.
\section{An\'alisis de multiresoluci\'on para valores puntuales}
Considerar $N_0=2^{n_0}$ valores
\begin{equation}
u^0=\{u_j^0\}_{j=1}^{N_0},
\end{equation}
correspondientes a los valores puntuales de una funci\'on $u(x)$ sobre una partici\'on uniforme de [-1,1]:
\begin{equation}
G^0=\{x_j^0\}_{j=0}^{N_0},\quad, x_j^0=-1+j\cdot h_L,\quad h_0=\frac{2}{N_0},\quad u_j^0=u(x_j^0),\quad 1\leqslant j\leqslant N_0.
\end{equation}
Se supone que $u(x)$ es 2-peri\'odica. Sus valores fuera de ]-1,1] son los de su extensi\'on peri\'odica: $u_0^0=u_{N_0}^0$, etc.

Considerar el conjunto de mallas anidadas di\'adicas $G^k,\ k=0,\ldots,L$:
\begin{equation}
G^k=\{x_j^k\}_{j=0}^{N_k},\quad, x_j^k=-1+j\cdot h_k,\quad h_k=2^{N_k+1}h_0,\quad N_k=\frac{N_0}{2^k},
\end{equation}
donde el nivel $k=0$ corresponde a la malla original, que es la m\'as fina; y $k=L$ corresponde a la malla m\'as gruesa. Notar que $G^k$ est\'a formada a partir de la malla m\'as fina $G^{k-1}$ eliminando las componentes de la malla con {\'i}ndice impar, es decir
\begin{equation}
G^{k-1}\setminus G^k=\{x^{k-1}_{2j-1}\}_{j=1}^{N_k},\qquad x_j^k=x_{2j}^{k-1},\quad 0\leqslant j\leqslant N_k.
\end{equation}
\begin{figure}[ht]
	\begin{center}
		\includegraphics[height=4.5cm]{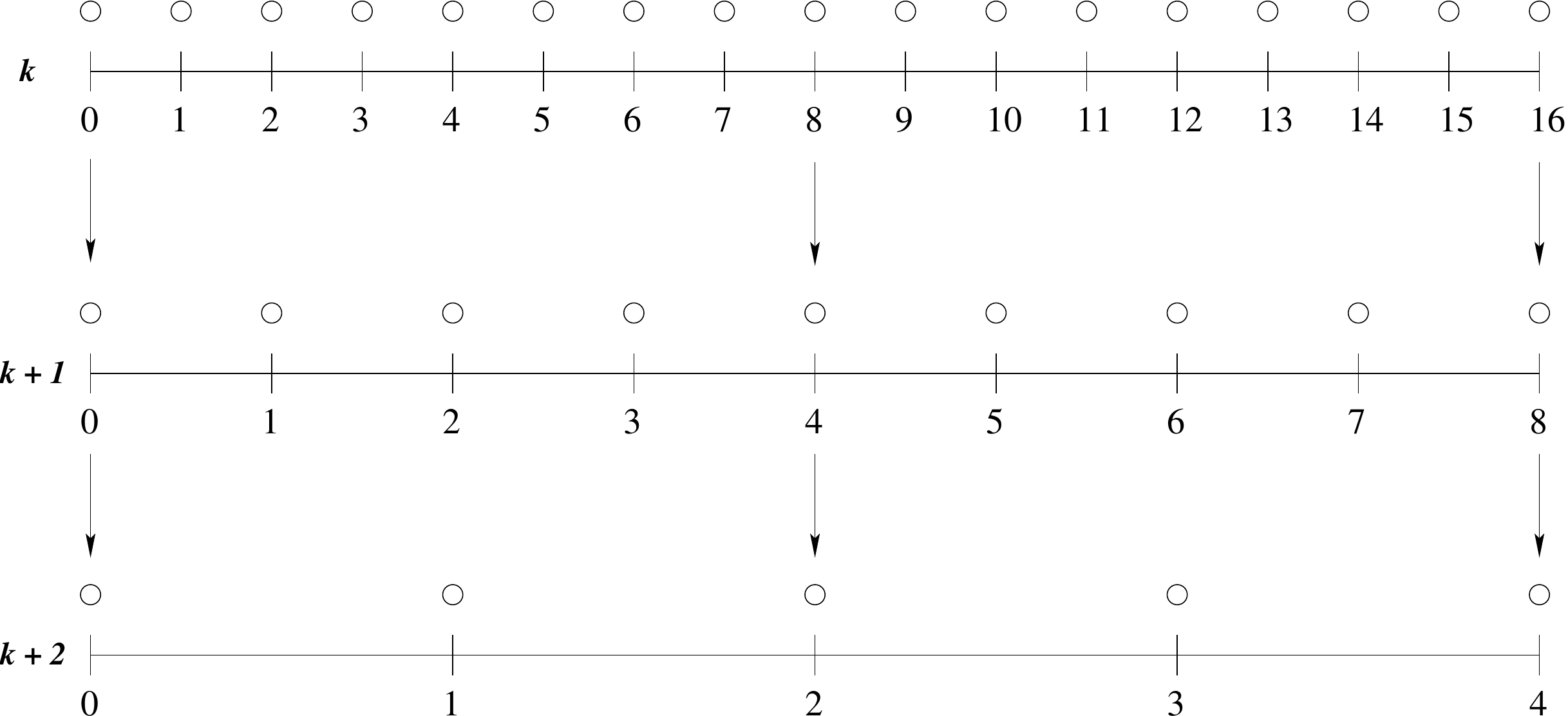}
		\caption{Diferentes escalas de valores puntuales}
		\label{fig:point}
	\end{center}
\end{figure}

Adem\'as se definen
\begin{equation}\label{aa1}
u_j^k=u(x_j^k)=u(x^0_{2^{k}j})=u^0_{2^kj},\qquad 0\leqslant j\leqslant N_k,
\end{equation}
por lo tanto este proceso (ver figura \ref{fig:point}) permite obtener $u^k$ a partir de $u^{k-1}$ mediante
\begin{eqnarray}
u_j^k&=&u_{2j}^{k-1},\quad 1\leqslant j\leqslant N_k,\\
u^{k-1}-u^k&=&\{u_{2j-1}^{k-1}\}_{j=1}^{N_k}.
\end{eqnarray}

Sea $\mathcal{I}(x,u^k)$ una funci\'on de interpolaci\'on de la malla $k$-\'esima, es decir,
\begin{equation}
\mathcal{I}(x_j^k,u^k)=u_j^k,\qquad 0\leqslant j\leqslant N_k,
\end{equation}
que puede utilizarse para obtener aproximaciones para los valores ausentes en la malla $k-1$\'esima
\begin{equation}
\tilde{u}_{2j-1}^{k-1}=\mathcal{I}(x_{2j-1}^{k-1},u^k),\qquad 0\leqslant j\leqslant N_k.
\end{equation}
Sea $D^k(u^0)=\{D_j^k\}_{j=1}^{N_k}$ la sucesi\'on de errores de interpolaci\'on al predecir los valores puntuales de cada nivel de resoluci\'on a partir del pr\'oximo nivel m\'as grueso
\begin{equation}
D_j^k=u_{2j-1}^{k-1}-\tilde{u}_{2j-1}^{k-1}=u_{2j-1}^{k-1}-\mathcal{I}(x_{2j-1}^{k-1},u^k),\ 1\leqslant j\leqslant N_k.
\end{equation}
Estos $D_j^k$ se conocen como \emph{coeficientes de ondelette} o \emph{detalles}. Es sencillo comprobar que los conjuntos de datos $(u^k,D^k)$ y $u^{k-1}$ contienen exactamente la misma informaci\'on,
\begin{equation}\label{equiv1}
u^{k-1}\leftrightarrow (u^k,D^k)
\end{equation}
en el sentido de que existe una transformaci\'on uno a uno entre ambos conjuntos (notar que la cardinalidad es la misma: $N_{k-1}=2N_k$).

Claramente utilizando (\ref{equiv1}) sucesivamente para $1\leqslant k\leqslant L$, se obtiene
\begin{eqnarray}
u^0\leftrightarrow (u^1,D^1)&\leftrightarrow &(D^1,(D^2,u^2))=(D^1,D^2,u^2)\leftrightarrow\cdots\\ \nonumber
&\leftrightarrow&(D^1,D^2,\ldots,D^L,u^L)=:(u_M)^T
\end{eqnarray}
donde $u_M=(D^1,D^2,\ldots,D^L,u^L)^T$ es la \emph{representaci\'on de multiresoluci\'on} de $u^0$, equivalente a la representaci\'on original. Esta permite extraer informaci\'on sobre la suavidad de la soluci\'on a partir de los errores de interpolaci\'on. La transformaci\'on uno a uno entre $u^0$ y $u_M$
\begin{equation}
u_M=Mu^0,\qquad u^0=M^{-1}u_M
\end{equation}
es lineal si $\mathcal{I}(\cdot,u^k)$ es independiente de los datos. En principio, puede utilizarse cualquier t\'ecnica de interpolaci\'on para $\mathcal{I}$. En este caso se utilizar\'a interpolaci\'on central polinomial
\begin{equation}\label{papa1}
\mathcal{I}(x,u^k)=q_j(x),\qquad x\in I_j=[x_{j-1},x_j],\ j=1,\ldots,N_k
\end{equation}
donde $q_j(x)$ es un plinomio de grado $r=2s$ un\'{\i}vocamente determinado por los datos $(u_{j-s}^k,\ldots,u_{j+s-1}^k)$ en los puntos $(x_{j-s}^k,\ldots,x_{j+s-1}^k)$; el valor en $x_{2j-1}^{k-1}$ se calcula a partir del polinomio de grado $r-1$ (es decir, cada est\'encil est\'a formado por $r$ puntos consecutivos de la malla) que interpola los puntos $(u_{j-s}^k,\ldots,u_{j+s-1}^k)$, por consiguiente
\begin{equation}
\tilde{u}_{2j-1}^{k-1}=\mathcal{I}(x_{2j-1}^{k-1},u^k)=\sum_{l=1}^s\beta_l(u_{j+l-1}^k+u_{j-l}^k),\quad r=2s,
\end{equation}
con
\begin{equation}\label{betas}
\left\{\begin{array}{ll}
r=2 \Rightarrow&\beta_1=1/2\\
r=4 \Rightarrow&\beta_1=9/16,\ \beta_2=-1/16\\
\end{array}\right.
\end{equation}
Ver detalles en el ap\'endice \ref{ap1}.

En este caso $M$ es un operador lineal que puede ser representado por una matriz de $N_0\times N_0$. Sin escribir la forma expl{\'i}cita de esta matriz, se sigue que $u_M=Mu^0$ puede ser calculado mediante el siguiente \emph{Algoritmo de Codificaci\'on}
\begin{equation}
\underline{u_M=Mu^0}\quad\left\{\begin{array}{l}
FOR\ k=1,2,\ldots,L\\
\qquad u_j^k=u_{2j}^{k-1},\quad 1\leqslant j\leqslant N_k,\\
\qquad D_j^k=u_{2j-1}^{k-1}-\sum_{l=1}^s\beta_l(u_{j+l-1}^k+u_{j-l}^k),\quad 1\leqslant j\leqslant N_k
\end{array}\right.
\end{equation}
y adem\'as $u^0=M^{-1}u_M$ puede ser  calculado mediante el siguiente \emph{Algoritmo de Decodificaci\'on}
\begin{equation}
\underline{u^0=M^{-1}u_M}\quad\left\{\begin{array}{l}
FOR\ k=L,L-1,\ldots,1\\
\qquad u_{2j}^{k-1}=u_j^k,\quad 1\leqslant j\leqslant N_k,\\
\qquad u_{2j-1}^{k-1}=\sum_{l=1}^s\beta_l(u_{j+l-1}^k+u_{j-l}^k)+D_j^k,\quad 1\leqslant j\leqslant N_k.
\end{array}\right.
\end{equation}

Notar que el algoritmo de Codificaci\'on va de fino a grueso mientras que el algoritmo de Decodificaci\'on va de grueso a fino; ambos son algoritmos cuyo costo computacional es de $O(N_0)$ operaciones ($(N_0-N_L)\cdot(s+1)$ sumas y $(N_0-N_L)\cdot s$ multiplicaciones).

Notar adem\'as que los algoritmos de Codificaci\'on y Decodificaci\'on representan una transformada de ondelette exacta, pues $u=\mathbf{M}^{-1}(\mathbf{M}u)$.
\section{An\'alisis de multiresoluci\'on por medias en celda}
En esta secci\'on se considera la sucesi\'on de $N_0$ valores
\begin{equation}
\bar{u}^0=\{\bar{u}^0_j\}_{j=1}^{N_0}
\end{equation}
que se interpretar\'an como medias en celda (\emph{cell-averages}) de cierta funci\'on $u(x)$ sobre la malla fina $G^0$:
\begin{equation}
\bar{u}^0_j=\frac{1}{h_0}\int_{x^0_{j-1}}^{x^0_j}u(x)dx,\quad 1\leqslant j\leqslant N_0.
\end{equation}

\begin{figure}[h]
	\begin{center}
		\includegraphics[height=3.2cm]{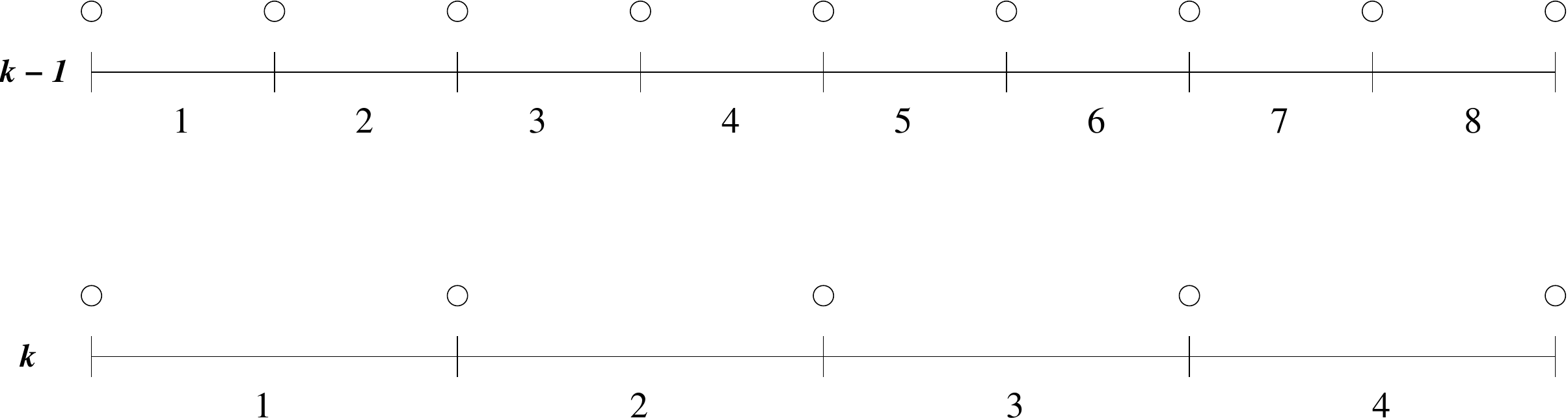}
		\caption{Diferentes escalas de medias en celda}
		\label{fig:cell}
	\end{center}
\end{figure}

Se consideran las mallas anidadas $G^k$, $1\leqslant k\leqslant L$ y se definen
\begin{equation}
\bar{u}^k_j=\frac{1}{h_k}\int_{x^k_{j-1}}^{x^k_j}u(x)dx,\quad 1\leqslant j\leqslant N_k.
\end{equation}
Se sigue de esta definici\'on y de (\ref{aa1}) que
\begin{eqnarray*}
\bar{u}^k_j&=&\frac{1}{h_k}\int_{x^k_{j-1}}^{x^k_j}u(x)dx\\
&=&\frac{1}{2h_{k-1}}\left(\int_{x^{k-1}_{2j-2}}^{x^{k-1}_{2j-1}}u(x)dx+\int_{x^{k-1}_{2j-1}}^{x^{k-1}_{2j}}u(x)dx\right)\\
&=&\frac{1}{2}(\bar{u}^{k-1}_{2j-1}+\bar{u}^{k-1}_{2j})
\end{eqnarray*}
Por lo tanto $\{\bar{u}^k_j\}_{j=1}^{N_k}$, $1\leqslant j\leqslant N_k$, puede ser calculado en forma directa del dato inicial $\bar{u}^0$, y sin ning\'un conocimiento expl{\'i}cito de la funci\'on $u(x)$, mediante el algoritmo
\begin{equation}
\left\{\begin{array}{l}
\textrm{FOR } k=1,2,\ldots,L\\
\quad \textrm{FOR } j=1,\ldots,N_k\\
\quad \quad \bar{u}^k_j=\frac{1}{2}(\bar{u}^{k-1}_{2j-1}+\bar{u}^{k-1}_{2j}).
\end{array}\right.
\end{equation}
Considerar la primitiva de $u(x)$
\begin{equation}
U(x)=\int_0^xu(y)dy,
\end{equation}
y observar que conocer las medias en celda $\bar{u}^k$ es equivalente al conocimiento de los valores puntuales $U^k$ de la funci\'on primitiva, es decir,
$$U^k=\{U^k_j\}_{j=1}^{N_k}\quad \leftrightarrow \quad \bar{u}^k=\{\bar{u}^k_j\}_{j=1}^{N_k},$$
lo cual es evidente de las siguientes dos relaciones:
\begin{equation}
U_j^k=U(x_j^k)=\int_0^{x_j^k}u(y)dy=\sum_{i=1}^j\int_{x_{i-1}^k}^{x_j^k}u(y)dy=h_k\sum_{i=1}^j\bar{u}^k_i,
\end{equation}
\begin{equation}
\bar{u}^k_j=\frac{U(x_j^k)-U(x_{j-1}^k)}{h_k}=\frac{U_j^k-U_{j-1}^k}{h_k}.
\end{equation}
En consecuencia conociendo los valores de $\bar{u}^k$ es posible calcular $U^k$ y utilizar una funci\'on de interpolaci\'on para aproximar el valor ausente $U_{2j-1}^{k-1}$, $1\leqslant j\leqslant N_k$ por $\tilde{U}_{2j-1}^{k-1}$, es decir,
\begin{equation}
\tilde{U}_{2j-1}^{k-1}=\mathcal{I}(x_{2j-1}^{k-1};U^k).
\end{equation}
Con esto, y teniendo en cuenta que $U_j^k=U^{k-1}_{2j}$, es posible lograr una aproximaci\'on $\tilde{u}^{k-1}$ para $\bar{u}^{k-1}$ mediante
\begin{equation}
\tilde{u}^{k-1}_{2j-1}=\frac{\tilde{U}_{2j-1}^{k-1}-\tilde{U}_{j-1}^k}{h_{k-1}},\quad \tilde{u}^{k-1}_{2j}=\frac{\tilde{U}_j^k-\tilde{U}_{2j-1}^{k-1}}{h_{k-1}}.
\end{equation}
Notar que
\begin{equation}
\frac{1}{2}(\tilde{u}^{k-1}_{2j-1}+\tilde{u}^{k-1}_{2j})=\frac{U_j^k-U_{j-1}^k}{2h_{k-1}}=\bar{u}^k_j,
\end{equation}
por lo tanto $\tilde{u}^{k-1}_{2j}$ puede calcularse a partir de $\bar{u}^k_j$ y $\tilde{u}^{k-1}_{2j-1}$ mediante
\begin{equation}
\tilde{u}^{k-1}_{2j}=2\bar{u}^k_j-\tilde{u}^{k-1}_{2j-1}.
\end{equation}
Se denota por $d^k(\bar{u}^0)=\{d_j^k\}_{j=1}^{N_k}$ a la sucesi\'on de errores de aproximaci\'on cometidos al predecir $\{\bar{u}^{k-1}_{2j-1}\}_{j=1}^{N_k}$ desde $\bar{u}^k$
\begin{eqnarray}
d_j^k&=&\bar{u}^{k-1}_{2j-1}-\tilde{u}^{k-1}_{2j-1}\\
&=&\bar{u}^{k-1}_{2j-1}-\frac{\mathcal{I}(x^{k-1}_{2j-1},U^k)-U^k_{j-1}}{h_{k-1}}.\label{cell12}
\end{eqnarray}
An\'alogamente al caso de valores puntuales, puede concluirse que existe una transformaci\'on uno a uno entre $\bar{u}^0$ y su representaci\'on de multiresoluci\'on
\begin{equation}\label{celda4}
\bar{u}_M=(d^1,\ldots,d^L,\bar{u}^L)^T,
\end{equation}
que se denota por
\begin{equation}
\bar{u}_M=\bar{M}\bar{u}^0,\qquad  \bar{u}^0=\bar{M}^{-1}\bar{u}_M.
\end{equation}
En (\ref{cell12}) el valor en $x^{k-1}_{2j-1}$ se calcula a partir de la funci\'on polinomial que interpola los puntos $(U_{j-s}^k,\ldots,U_{j+s-1}^k)$. Utilizando lo visto anteriormente para el caso de valores puntuales, y como $2h_{k-1}=h_k$, se obtiene
\begin{equation}
d_j^k=\bar{u}_{2j-1}^{k-1}-\frac{\sum_{l=1}^{s-1}\beta_l(U_{j+l-1}^k+U_{j-l}^k)-U_{j-1}^k}{2h_k},
\end{equation}
con los $\beta_l$ calculados en (\ref{betas}).

De este modo, los coeficientes de ondelette est\'an dados por
\begin{equation}
d_j^k=\bar{u}_{2j-1}^{k-1}-\bar{u}_j^k-\sum_{l=1}^{s-1}\gamma_l(\bar{u}_{j+l}^k-\bar{u}_{j-l}^k),\quad 1\leqslant j\leqslant N_k.
\end{equation}
Notar que se utiliza el est\'encil $(\bar{u}_{j-s+1}^k,\ldots,\bar{u}_{j+s-1}^k)$ y por lo tanto el orden de precisi\'on correspondiente es $\bar{r}=2s-1$, y los coeficientes correspondientes $\gamma_l$ son
$$\left\{\begin{array}{ll}
r=3 \Rightarrow&\gamma_1=-1/8\\
r=5 \Rightarrow&\gamma_1=-22/128,\ \gamma_2=3/128
\end{array}\right.
$$
Ver detalles en el ap\'endice \ref{ap2}.

Cuando se utiliza interpolaci\'on central (o cualquier interpolaci\'on independiente de los datos), se tiene que $\bar{M}$ es un operador lineal que puede ser expresado por una matriz de $N_0\times N_0$. En el caso de que $\mathcal{I}(\cdot,U^k)$ sea el especificado en la secci\'on anterior, con $r$ y $s$ dados, las transformaciones en (\ref{celda4}) pueden ser llevadas a cabo sin escribir la forma expl{\'i}cita de la matriz, calculadas mediante los algoritmos siguientes:

\noindent\emph{Algoritmo de Codificaci\'on}
\begin{equation}
\underline{\bar{u}_M=\bar{M}\bar{u}^0}\quad\left\{\begin{array}{l}
FOR\ k=1,2,\ldots,L\\
\qquad \bar{u}_j^k=\frac{1}{2}(\bar{u}_{2j-1}^{k-1}+\bar{u}_{2j}^{k-1}),\quad 1\leqslant j\leqslant N_k,\\
\qquad d_j^k=\bar{u}_{2j-1}^{k-1}-\bar{u}_j^k-\sum_{l=1}^{s-1}\gamma_l(\bar{u}_{j+l}^k-\bar{u}_{j-l}^k),\ 1\leqslant j\leqslant N_k
\end{array}\right.
\end{equation}

\noindent\emph{Algoritmo de Decodificaci\'on}
\begin{equation}
\underline{u^0=M^{-1}u_M}\quad\left\{\begin{array}{l}
FOR\ k=L,L-1,\ldots,1\\
\qquad FOR\ j=1,\ldots,N_k\\
\qquad \qquad \Delta= \sum_{l=1}^{s-1}\gamma_l(\bar{u}_{j+l}^k-\bar{u}_{j-l}^k)+d_j^k,\\
\qquad \qquad \bar{u}_{2j-1}^{k-1}=\bar{u}_j^k+\Delta,\quad \bar{u}_{2j}^{k-1}=\bar{u}_j^k-\Delta.
\end{array}\right.
\end{equation}

Ambos algoritmos poseen un costo computacional de $O(N_0)$ operaciones ($(N_0-N_L)\cdot(s+2)$ sumas en ambos algoritmos y $(N_0-N_L)\cdot s$ multiplicaciones en ambos algoritmos).

Es interesante observar que dado que $\bar{u}^0$ es equivalente a $U^0$, tambi\'en $\bar{u}_M$ es equivalente a $U_M$, la representaci\'on de multiresoluci\'on de los valores puntuales de la funci\'on primitiva $U(x)$
\begin{equation*}
(d^1,d^2,\ldots,d^L,\bar{u}^L)^T=\bar{u}_M\leftrightarrow U_M=(D^1,D^2,\ldots,D^L,U^L)^T.
\end{equation*}
Adem\'as la transformaci\'on entre $d_j^k(\bar{u}^0)$ y $D_j^k(U^0)$ est\'a dada por
\begin{equation}\label{relacion4}
d_j^k(\bar{u}^0)=D_j^k(U^0)/h_{k-1}.
\end{equation}
\section{An\'alisis de regularidad}
El an\'alisis de multiresoluci\'on ser\'a de gran utilidad para obtener un algoritmo de compresi\'on de datos de las medias en celda. Luego se estudiar\'a su aplicaci\'on a la soluci\'on num\'erica $v^n$ del esquema conservativo
\begin{equation}
v_j^{n+1}=v_j^n-\lambda(\bar{f}_j-\bar{f}_{j-1}),\qquad \lambda=\tau/h.
\end{equation}
Utilizando resultados de interpolaci\'on est\'andar y notando que $U(x)$ es m\'as suave que $u(x)$, se obtiene de (\ref{relacion4}) la siguiente caracterizaci\'on cualitativa de $d_j^k(\bar{u}^0)$ (ver \cite{Harten1}):

\begin{Th} Si la funci\'on $u(x)$ en $x=x^*$ posee $p-1$ derivadas continuas y una discontinuidad de salto en la derivada $p-$\'esima, entonces en los puntos $x_j^k$ cercanos a $x^*$ se tiene
\begin{equation}\label{teor1}
d_j^k(\bar{u}^0)\sim\left\{\begin{array}{ll}
(h_k)^p[u^{(p)}],& si\  0\leqslant p\leqslant \bar{r},\\
(h_k)^pu^{(p)},& si\ p>\bar{r},
\end{array}\right.
\end{equation}
donde $\bar{r}$ es el orden de precisi\'on de la aproximaci\'on ($\bar{r}=r-1$), $p\leqslant 1$ y $[\ ]$ denota el salto en la discontinuidad.
\end{Th}
\dem Sea $\mathcal{I}(x,U^{k-1})$ como en (\ref{papa1}). Se tiene que
\begin{equation}\label{papa2}
U(x)=\mathcal{I}(x,U^{k-1})+U[x_{j-s}^{k-1},\ldots,x_{j+s-1}^{k-1},x]\prod_{i=j-s}^{j+s-1}(x-x_i^{k-1}),
\end{equation}
con $x\in[x_{j-1}^{k-1},x_j^{k-1}]$. Notar que si $u(x)$ tiene $p-1$ derivadas continuas en $x^*$ y una discontinuidad de salto en $u^{(p)}$ cerca de $x^*$, entonces $U(x)$ tiene $p$ derivadas continuas en $x^*$ y una discontinuidad de salto en $U^{(p+1)}$ cerca de $x^*$. Con esto, de \cite{Arandiga} se deduce que
\begin{equation}
U[x_l^{k-1},\ldots,x_{l+t}^{k-1}]=\left\{\begin{array}{ll}
\frac{O([U^{(p+1)}])}{h_k^{t-(p+1)}},&\textrm{ si } 0\leqslant p+1\leqslant t\\
O(\|U^{(t)}\|),&\textrm{ si }t<p+1.
\end{array}\right.\end{equation}
Dado que $D_j^k=U^{k-1}_{2j-1}-\mathcal{I}(x_{2j-1}^{k-1},U^k)$, la relaci\'on (\ref{papa2})
conduce a
\begin{equation}
D_j^k(U)=U[x_{j-s}^k,\ldots,x_{j+s-1}^k,x_{2j-1}^{k-1}]\prod_{i=j-s}^{j+s-1}(x_{2j-1}^{k-1}-x_i^k),
\end{equation}
y teniendo en cuenta que $x_{2j-1}^{k-1}-x_i^k$ es aproximadamente del orden de $h_k$, con $i \in \{j-s,\ldots,j+s-1\}$, se obtiene que
\begin{equation}\label{teor3}
D_j^k(U)\sim\left\{\begin{array}{ll}
\frac{[U^{(p+1)}]}{h_k^{r-(p+1)}}h^r_k,&\textrm{ si } 0\leqslant p+1\leqslant t\\
\|U^{(r)}\|h^r_k,&\textrm{ si }t<p+1.
\end{array}\right.
\end{equation}
Finalmente, de (\ref{relacion4}), (\ref{teor3}) y remarcando que $U^{(n+1)}(x)\equiv u^{(n)}(x)$, se obtiene (\ref{teor1}).

\qedo

Ahora, la ecuaci\'on (\ref{teor1}) en el nivel $k-1$ corresponde a
\begin{equation}\label{teor2}
d_j^{k-1}\sim\left\{\begin{array}{ll}
(h_{k-1})^p[u^{(p)}],& si\  0\leqslant p\leqslant r-1,\\
(h_{k-1})^{r-1}u^{(r-1)},& si\ p>r-1,
\end{array}\right.
\end{equation}
y como $h_k=2h_{k-1}$, entonces
\begin{equation}
d_{2j}^k\sim\left\{\begin{array}{ll}
2^{-p}(h_{k-2})^p[u^{(p)}],& si\  0\leqslant p\leqslant r-1,\\
2^{-r+1}(h_{k-2})^{r-1}u^{(r-1)},& si\ p>r-1,
\end{array}\right..
\end{equation}
Por lo tanto
\begin{equation}
|d_{2j}^{k-1}|\approx 2^{-\bar{p}}|d_j^k|,\quad \bar{p}=\min(p,\bar{r}).
\end{equation}

Pueden obtenerse entonces algunas conclusiones \'utiles
\begin{itemize}
\item Lejos de las discontinuidades, los coeficientes $d_j^k$ decrecen a medida que se va a niveles m\'as finos.
\item La tasa de decaimiento de los coeficientes $d_j^k$ es determinada por la regularidad local de la funci\'on y el orden de precisi\'on de la aproximaci\'on.
\item En la vecindad de una irregularidad de $u(x)$, los coeficientes $d_j^k$ permanecen del mismo orden $O([u])$, independiente del nivel de refinamiento.
\end{itemize}
Por lo tanto el an\'alisis de multiresoluci\'on de $\bar{u}^0$ puede verse como un estudio de la regularidad local de $u(x)$.

Puede hacerse un an\'alisis de regularidad similar si se considera el caso de valores puntuales en vez de medias en celdas. De forma an\'aloga, Kozakevicius (ver \cite{Alice}) propone que dependiendo de la regularidad de la funci\'on, un gran n\'umero de coeficientes de  ondelette pueden ser extremadamente peque\~nos, y por lo tanto podr\'{\i}an ser descartados de la representaci\'on de multiresoluci\'on.

\section{Compresi\'on de datos}
La idea principal es reducir la cantidad de datos mediante una t\'ecnica de truncamiento, que consiste en hacer ceros los coeficientes que est\'an por debajo de una tolerancia prescrita.

\noindent Sea $\mathbf{tr_{\varepsilon_k}}$ el operador de truncamiento definido por
\begin{equation}
\hat{d_j^k}=\mathbf{tr_{\varepsilon_k}}(d_j^k)=\left\{\begin{array}{ll}
0,& \textrm{ si }|d_j^k|\leqslant \varepsilon_k\\
d_j^k,& \textrm{ en otro caso}.
\end{array}\right.,
\end{equation}
Sea $\hat{u_M}$ el resultado de la operaci\'on de truncamiento aplicada a $u_M$
\begin{equation}
\hat{u_M}=(\hat{d^1},\hat{d^2},\ldots,\hat{d^L},u^L).
\end{equation}
Si se aplica el algoritmo de decodificaci\'on al dato truncado $\hat{u_M}$, se obtiene una aproximaci\'on $\tilde{u^0}=\mathbf{M}^{-1}\hat{u_M}$, que por \cite{Harten1} se sabe que permanece \emph{cerca} del dato inicial $u^0$.

Dado que se est\'a en el caso de multiresoluci\'on por valores puntuales de $u$, cada coeficiente de ondelette est\'a relacionado con una posici\'on espec\'{\i}fica en la malla fina uniforme y por lo tanto los procesos de codificaci\'on y decodificaci\'on pueden ser simplificados. Los coeficientes $d_j^k$ se calculan entonces s\'olo para decidir si $x_j^k$ seguir\'a o no en la malla y se evita as\'{\i} construir la representaci\'on de multiresoluci\'on completa \cite{Alice}. Esto quiere decir, que en estos puntos, la informaci\'on sobre la funci\'on puede ser obtenida mediante interpolaci\'on.

La representaci\'on de $u^0$ al cabo de este proceso, contendr\'a s\'olo los valores puntuales en las posiciones asociadas a coeficientes de ondelette significativos, y los puntos en el nivel m\'as grueso. Esto se conoce como \emph{representaci\'on puntual esparsa} de $u$, y se denota por $u_S$.

La elecci\'on de $\varepsilon_k$ puede variar de acuerdo a las propiedades de los espacios funcionales \cite{Holm}, o suavidad de la funci\'on \cite{Harten1}.
En este caso, con $\varepsilon$ fijo, los niveles de tolerancia en cada nivel estar\'an dados por $\varepsilon_k=\varepsilon/2^{L-k}$. Notar que a escalas m\'as finas, $\varepsilon_k$ es m\'as peque\~no; esto con el fin de preservar la informaci\'on asociada a la parte regular del dato inicial y descartar perturbaciones de alta frecuencia (pues una se\~nal regular posee mayores coeficientes de ondelette en escalas m\'as gruesas y una se\~nal perturbada, o una funci\'on con singularidades, posee mayores coeficientes de ondelette en escalas m\'as finas). Adem\'as esta elecci\'on de $\varepsilon_k$ es \'optima en el sentido que mantiene la mejor relaci\'on entre compresi\'on de datos y disipaci\'on de informaci\'on durante la evoluci\'on temporal de la soluci\'on.

\begin{figure}[ht]
\begin{center}
\includegraphics[height=5cm]{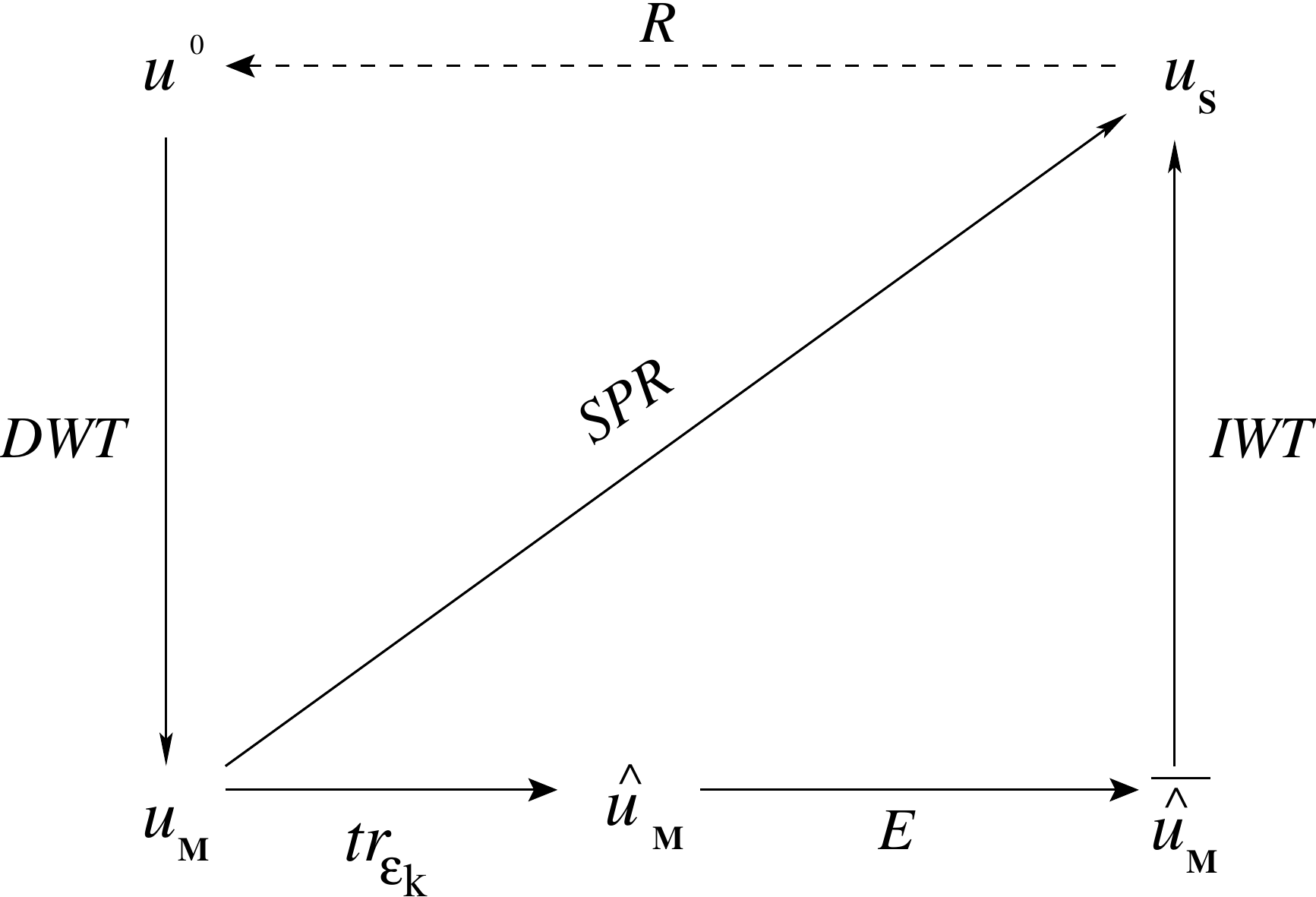}
\caption{\small Secuencia de operaciones para obtener la representaci\'on puntual esparsa de una funci\'on. $DWT$: transformada del dato inicial, $\mathbf{tr}_{\varepsilon_k}$: operador de truncamiento, $E$: inclusi\'on de \emph{safety points}, $IWT$: transformada de ondelette inversa y $R$: reconstrucci\'on de malla uniforme.}
\end{center}
\end{figure}

La representaci\'on puntual esparsa $u_S$ tambi\'en incluir\'a algunos \emph{safety points} necesarios para evitar la disipaci\'on num\'erica; este corresponde al operador de extensi\'on $\mathbf{E}$. Los safety points ser\'an incluidos en las vecindades de puntos cuyos coeficientes de ondelette son significativos \cite{Alice}. Se incluir\'an dos tipos de safety points: Puntos en el mismo nivel de multiresoluci\'on que el coeficiente de ondelette respectivo (con el fin de mantener la calidad del transporte de informaci\'on desde un punto a su vecino en la malla) y puntos en un nivel de multiresoluci\'on m\'as fino que el nivel del coeficiente de ondelette (s\'olo si el detalle es mayor que una tolerancia adicional $2\varepsilon_k$, esto con el fin de mejorar la captura de choques).

\section{Estructura de datos}
Dado las caracter\'{\i}sticas de los problemas hiperb\'olicos que poseen discontinuidades que se propagan, el n\'umero de puntos en la representaci\'on puntual esparsa es mucho menor que el n\'umero de puntos en la malla fina uniforme. Luego, ser\'a de gran utilidad almacenar la informaci\'on relevante en alg\'un tipo de estructura que saque provecho de ello, tal como se hace en \cite{Alice} (MORSE, SPARSE, etc.)

\begin{figure}[h]
\begin{center}
\includegraphics[height=5cm]{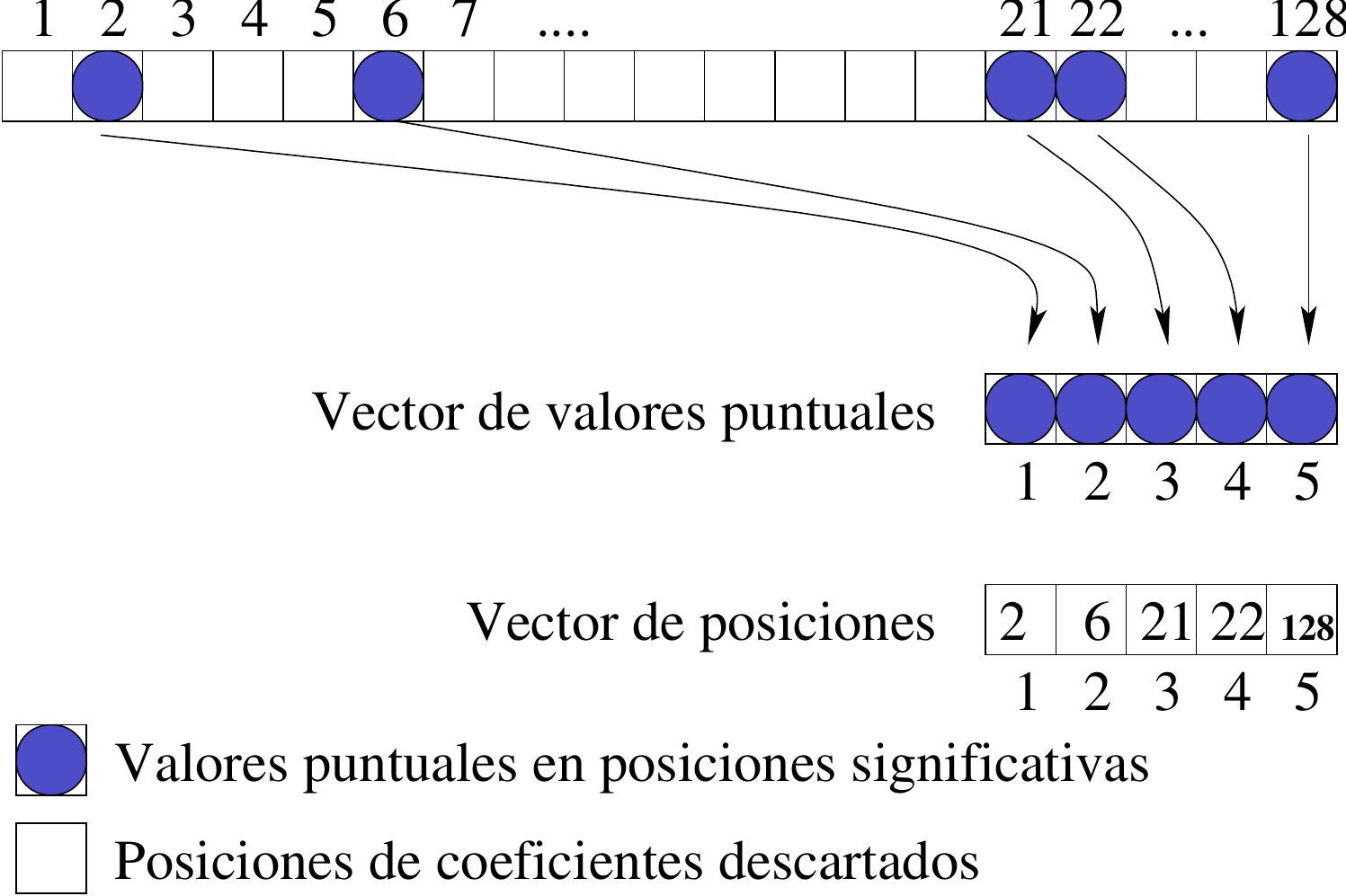}
\caption{\small Ejemplo de almacenamiento de datos s\'olo para posiciones significativas de la representaci\'on truncada (MORSE o SPARSE).}
\end{center}
\end{figure}

\clearpage{\pagestyle{empty}\cleardoublepage}
\chapter{Caso hiperb\'olico}
En esta secci\'on se presenta una forma eficiente de resolver leyes de conservaci\'on hiperb\'olicas mediante un m\'etodo de alta resoluci\'on en mallas generadas por ondelettes desarrollado por Kozakevicius y Santos \cite{Alice}. La eficiencia de este m\'etodo se basa en la asociaci\'on de dos t\'ecnicas independientes: mallas adaptativas generadas por una transformaci\'on de ondelettes \cite{Harten1,Cohen,Holm} y m\'etodos de alta resoluci\'on basados en interpolaciones ENO para el c\'alculo de los flujos \cite{ENO,Alice}.
\section{Esquema ENO Lax-Friedrichs}
Se necesitan esquemas conservativos para la parte espacial del operador (forma semi-discreta)
$$\frac{d}{dt}(u_j(t))=\frac{-1}{\Delta x_j}\left(\hat{f}_{j+1/2}-\hat{f}_{j-1/2}\right),$$
donde $\hat{f}_{j+1/2}=\hat{f}(u_{j-r},\ldots,u_{j-s})$ es el flujo num\'erico, en que la primera posici\'on del est\'encil $j-r$ es elegida mediante un algoritmo ENO, manteniendo la relaci\'on $j-r<j+1/2<j-s$. Esta funci\'on de flujo num\'erico es Lipschitz continua en sus argumentos y es consistente con el flujo exacto, es decir, $\hat{f}(u,\ldots,u)=f(u)$.

Para lograr un alto orden de aproximaci\'on para $\frac{\partial f}{\partial x_j}$, se utilizar\'an posiciones escalonadas auxiliares $\{x_{j+1/2}\}_j$ \cite{GR} con respecto a la malla gruesa esparsa. El flujo num\'erico evaluado en estas posiciones se obtiene mediante interpolaci\'on ENO.

\begin{figure}[ht]
\begin{center}
\includegraphics[height=5cm]{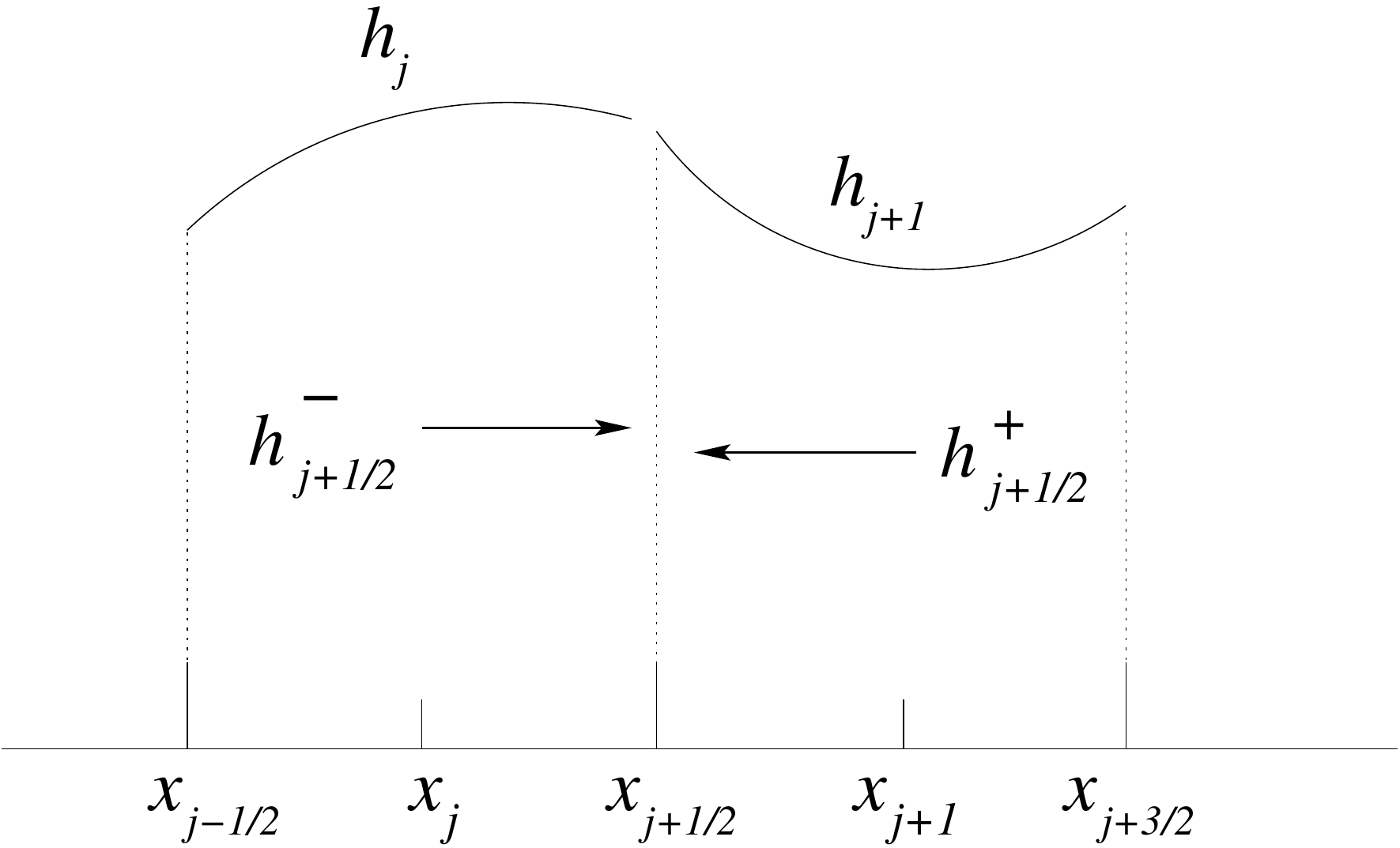}
\caption{Componentes de la separaci\'on del flujo num\'erico en la frontera. $h_j$ es el interpolador ENO para la celda $]x_{j-1/2},x_{j+1/2}[$ y $h_{j+1}$ es el interpolador ENO para la celda $]x_{j+1/2},x_{j+1+1/2}[$.}
\end{center}
\end{figure}

Es necesario considerar esquemas upwind en la construcci\'on del flujo num\'erico con el objetivo de mantener la estabilidad del esquema. Para ello se utilizar\'a la forma m\'as sencilla, m\'as robusta y menos costosa de obtener esquemas upwind sin violar condiciones de entrop{\'i}a de la soluci\'on. Esta es, la separaci\'on de flujo de \emph{Lax-Friedrichs}:
$$f(u)=f^+(u)+f^-(u),\quad f^+(u)=\frac{1}{2}(f(u)+\alpha u),\quad f^-(u)=\frac{1}{2}(f(u)-\alpha u),$$
donde
\begin{equation}\label{alpha}
\alpha=\max_u|f'(u)|.
\end{equation}
El n\'umero de puntos escogidos para la reconstrucci\'on depende del orden de la interpolaci\'on. En este caso, se utilizar\'a interpolaci\'on c\'ubica.

El flujo num\'erico en las posiciones de la malla auxiliar corresponde a la suma de las aproximaciones generadas para cada parte de la separaci\'on de flujos
\begin{equation}
\hat{f}_{j+1/2}=\hat{f}^+_{j+1/2}+\hat{f}^-_{j+1/2}.
\end{equation}
Notar que $\hat{f}^+_{j+1/2}$ y $\hat{f}^-_{j+1/2}$ son aproximaciones para el mismo borde $x_{j+1/2}$ del volumen de control $]x_{j-1/2},x_{j+1/2}[$, obtenidas mediante interpoladores distintos.

Notar adem\'as, que una vez que se elige el n\'umero de puntos en el est\'encil, este permanece igual para todos los puntos de la malla $G^k$. Esta forma de construir predictores para la transformada de ondelette no considera la suavidad local de la funci\'on a ser interpolada. Si la funci\'on es suave a trozos, una aproximaci\'on de est\'encil fijo puede comportarse inadecuadamente cerca de las irregularidades, generando oscilaciones en las celdas correspondientes. Estas oscilaciones (conocidas como \emph{Fen\'omeno de Gibbs} en m\'etodos espectrales) ocurren debido a que los est\'enciles contienen una \emph{celda discontinua} (volumen de control que contiene una irregularidad), es decir, poseen un punto $x_j$ bastante cerca de una irregularidad. Adem\'as, cada vez que el est\'encil cruza una singularidad, la calidad de la interpolaci\'on se ve reducida. Cuanto mayor es el grado del interpolador, mayor es la regi\'on afectada por la singularidad.

La idea es entonces utilizar interpolaci\'on ENO (\emph{Essentially non oscillatory}), que aumenta la regi\'on de precisi\'on para el interpolador, eligiendo un est\'encil diferente, para evitar las oscilaciones cerca de las discontinuidades.

Se presenta a continuaci\'on la forma en que se prepara la reconstrucci\'on ENO. Inicialmente se conocen los valores de los flujos en la malla esparsa $S$. Se define $V(x_{j+1/2})$, la primitiva de la componente de separaci\'on de flujo en la malla auxiliar con respecto a $S$. Se construir\'a un polinomio interpolador por partes de $V$, en la variable $x$: $H(x,V)$, sobre la malla auxiliar, es decir,
\begin{eqnarray*}
H(x_{j+1/2},V)=V_{j+1/2}&=&V(x_{j+1/2})=\sum_{k=0}^jf(x_k),\\
H(x,V)=q_m(x,V),& & x_{j-1/2}\leqslant x\leqslant x_{j+1/2},
\end{eqnarray*}
donde $q_m$ es el \'unico polinomio interpolador de grado $m$, que utiliza $m+1$ puntos consecutivos $(x_{i_m(j)},\ldots,x_{i_m(j)+m})$, incluyendo a $x_{j-1/2}$ y $x_{j+1/2}$.

Notar que dependiendo de la elecci\'on del primer punto del est\'encil $i_m(j)$, existen $m$ polinomios interpoladores posibles. >Cu\'al elegir? El est\'encil asociado a $[x_{j-1/2},x_{j+1/2}]$ ser\'a aquel tal que $V(x)$ es m\'as suave (en un sentido asint\'otico) y el valor $x$ donde se evaluar\'a el interpolador, ser\'a $x_{j-1/2}$ o $x_{j+1/2}$.

La informaci\'on de la suavidad de $V$ puede obtenerse de las diferencias divididas:
\begin{eqnarray*}
w[x_{j-1/2}]&=&V(x_{j-1/2})\\
w[x_{j-1/2},\ldots,x_{j-1/2+k}]&=&\frac{w[x_{j-1/2+1},\ldots,x_{j-1/2+k}]-w[x_{j-1/2},\ldots,x_{j-1/2+k-1}]}{x_{j-1/2+k}-x_{j-1/2}}.
\end{eqnarray*}
El siguiente Teorema (ver \cite{Alice,Donat}) entrega un criterio para medir asint\'oticamente la suavidad de una funci\'on.
\begin{Th} Si $f(x)$ es $C^\infty([x_i,x_{i+k}])$, entonces
\begin{equation}
w[x_i,\ldots,x_{i+k}]=\frac{1}{k!}\frac{d^k}{dx^k}f(\xi_{i,k}),\quad x_i\leqslant \xi_{i,k}\leqslant x_{i+k};
\end{equation}
pero si $f(x)$ tiene una discontinuidad de salto en su $p-$\'esima derivada, $0\leqslant p \leqslant k$, entonces
\begin{equation}
w[x_i,\ldots,x_{i+k}]=O(d^{-k+p})[f^{(p)}],\quad d=|x_{i+k}-x_i|.
\end{equation}
\end{Th}
Luego, utilizando $|w[x_{j-1/2},\ldots,x_{j-1/2+k}]|$ es posible medir asint\'oticamente la suavidad de $f(x)$ en $[x_{j-1/2},x_{j-1/2+k}]$: El mejor est\'encil, ser\'a aquel asociado a la diferencia dividida m\'as peque\~na \cite{Donat}.

La cuesti\'on es ahora, c\'omo hallar $i_m(j)$. Para ello se seguir\'a el siguiente procedimiento (\cite{Alice}):
\begin{enumerate}
\item $i_1(j)=j-\frac{1}{2}$, donde $q_1$ es el polinomio interpolador para $V$ en $x_{j-1/2}$ y $x_{j+1/2}$.
\item Suponer un polinomio interpolador de grado $n$, $q_n$ para $V$ en $x_{i_n(j)},\ldots,x_{i_n(j)+n}$.
\item De acuerdo con la diferencia dividida m\'as peque\~na, $q_{n+1}$ comenzar\'a con $i_{n+1}(j)$ $=$ $i_n(j)-1$ (si el siguiente punto elegido est\'a a la izquierda del \'ultimo punto en el est\'encil) o con $i_{n+1}(j)$$=i_n(j)$ (si el siguiente punto elegido est\'a a la derecha del \'ultimo punto en el est\'encil).
\end{enumerate}
Las aproximaciones $h_{j-1/2}$ y $h_{j+1/2}$ para cada componente de la separaci\'on de flujo ser\'an entonces la derivada de $H$ evaluada en $x_{j-1/2}$ y $x_{j+1/2}$ respectivamente. Una vez calculado el flujo num\'erico en las posiciones auxiliares, se obtiene una aproximaci\'on de alto orden para el t\'ermino de la derivada espacial en las posiciones de la malla esparsa $S$.
\section{Evoluci\'on temporal}\label{Etemp}
Notar que se est\'a frente a un proceso de discretizaci\'on en dos etapas, primero se ha discretizado s\'olo espacio, dejando el problema continuo en tiempo. Esto conduce a las llamadas \emph{ecuaciones semi-discretas}. La discretizaci\'on puede hacerse utilizando un m\'etodo num\'erico est\'andar para sistemas de ecuaciones diferenciales ordinarias. Este mecanismo es particularmente ventajoso en el desarrollo de m\'etodos con orden de precisi\'on \emph{mayor a dos}, ya que permite alcanzar de forma relativamente sencilla la misma precisi\'on espacial y temporal.

Los experimentos realizados en \cite{ENO} indican que las formulaciones semi-discretas con discretizaci\'on temporal Runge-Kutta TVD desarrollados por Shu y Osher no generan oscilaciones para $CFL \leqslant 0.5$ aproximadamente, y son \'optimas en el sentido de que permiten el mayor $CFL$ para esquemas expl{\'i}citos, $CFL=1$.

Se utilizar\'an entonces m\'etodos Runge-Kutta TVD de segundo o tercer orden.

\underline{R-K TVD \'optimo de segundo orden:}
\begin{eqnarray*}
u^{(1)}&=&u^n+\Delta t\mathcal{L}(u^n)\\
u^{n+1}&=&\frac{1}{2}u^n+\frac{1}{2}u^{(1)}+\Delta t\mathcal{L}(u^{(1)}),
\end{eqnarray*}

\underline{R-K TVD \'optimo de tercer orden orden:}
\begin{eqnarray*}
u^{(1)}&=&u^n+\Delta t\mathcal{L}(u^n)\\
u^{(2)}&=&\frac{3}{4}u^n+\frac{1}{4}u^{(1)}+\Delta t\mathcal{L}(u^{(1)})\\
u^{n+1}&=&\frac{1}{3}u^n+\frac{2}{3}u^{(2)}+\frac{2}{3}\Delta t\mathcal{L}(u^{(2)}),
\end{eqnarray*}
con $\mathcal{L}(u)=-(\Delta x)^{-1}(\hat{f}_j(u)-\hat{f}_{j-1}(u))$.

La alta resoluci\'on (asociada a discretizaciones espaciales TVB, ENO o TVD) es necesaria para asegurar estabilidad. En los pasos intermedios del esquema de evoluci\'on temporal, se conserva la malla esparsa del paso $n$.
\section{Adaptatividad de la representaci\'on esparsa}\label{remesh}
Con el fin de actualizar la malla esparsa, es necesario aplicar el operador de reconstrucci\'on $R$ para reconstruir la soluci\'on en la malla uniforme.  Una vez aplicada la transformada de ondelette, el operador de truncamiento y el operador de extensi\'on, puede llevarse a cabo la  evoluci\'on temporal.

Dado que recalcular la malla es costoso, puede utilizarse la misma malla para varios pasos temporales. Para problemas donde la velocidad de la onda es baja (en el sentido $CFL$), es posible utilizar la misma representaci\'on puntual esparsa para 5 o m\'as pasos temporales sin aumentar la disipaci\'on num\'erica, y luego realizar la actualizaci\'on de la configuraci\'on. Para problemas con una alta velocidad de onda, la reconstrucci\'on de la malla puede hacerse cada dos pasos temporales, sin afectar la calidad de la soluci\'on \cite{Alice}.

Cuando se trabaja con ecuaciones multivariadas se construye una malla esparsa ``unificada''. Es la uni\'on de las posiciones significativas de la representaci\'on esparsa de cada componente y todos los safety points necesarios para la evoluci\'on, en cada componente. El criterio para la malla unificada es bastante simple. Una vez que una posici\'on tiene asociado un coeficiente de ondelette significativo en cualquier componente del vector de cantidades, tal posici\'on debe permanecer en la malla unificada, y todas las componentes del vector de cantidades deben tener sus valores puntuales en esta posici\'on. Lo mismo sucede con el operador de extensi\'on.

Notar que como cada variable del vector de cantidades desarrolla discontinuidades bastante localizadas, la malla unificada \emph{seguir\'a siendo esparsa} \cite{Harten1,Harten2,Alice}.

La actualizaci\'on de la malla es an\'aloga al caso escalar. Los mismos operadores deben ser aplicados a cada componente del vector de cantidades para obtener la siguiente configuraci\'on de la malla unificada y realizar la evoluci\'on temporal.

\section{M\'etodo adaptativo de alta resoluci\'on}
Dado el n\'umero de puntos en la malla fina, $N_0$, el n\'umero de niveles de multiresoluci\'on, $L$, el grado $r$ del predictor intermallas y del interpolador ENO, el nivel de truncamiento $\varepsilon_k$; dadas adem\'as las condiciones de contorno e inicial de la ley de conservaci\'on, el algoritmo del m\'etodo descrito puede ser resumido como sigue:
\begin{enumerate}
\item \begin{itemize}
\item Transformada de ondelette discreta ($DWT$) (u operador de codificaci\'on $\mathbf{M}$) aplicada al dato inicial.
\item Representaci\'on puntual esparsa ($SPR$) de la soluci\'on. Esta incluye truncamiento, extensi\'on, y transformada inversa de ondelette ($IWT$) (u operador de decodificaci\'on $\mathbf{M}^{-1}$).
\end{itemize}
\item \begin{itemize}
\item C\'alculo del flujo exacto en malla esparsa (correspondiente al nivel m\'as fino de multiresoluci\'on).
\item C\'alculo del valor global de $\alpha$ (\ref{alpha}).
\item C\'alculo de $\Delta t$ para la evoluci\'on temporal: $\Delta t=\frac{CFL\cdot h_0}{\alpha}$, donde $h_0$ es el paso espacial en la malla fina.
\item Factorizaci\'on Lax-Friedrichs del flujo exacto: $f^+$ y $f^-$.
\item C\'alculo del flujo num\'erico $\hat{f}_{j+1/2}$:
\begin{itemize}
\item para $f^+$, construir la aproximaci\'on ENO $h^-$.
\item para $f^-$, construir la aproximaci\'on ENO $h^+$.
\item $\hat{f}_{j+1/2}=\hat{f}^+_{j+1/2}+\hat{f}^-_{j+1/2}$.
\end{itemize}

\end{itemize}

\item \begin{itemize}
\item Evoluci\'on temporal: Runge-Kutta TVD de segundo o tercer orden. Se necesitan pasos intermedios.
\item Repetir 2. para la soluci\'on intermedia necesaria para 3.
\item Evoluci\'on temporal de la soluci\'on intermedia (El m\'etodo Runge-Kutta TVD de segundo orden completa el paso temporal, el m\'etodo Runge-Kutta TVD de tercer orden necesita otro paso intermedio).
\end{itemize}
\item Aplicaci\'on del operador de reconstrucci\'on de la soluci\'on en malla fina $R$.
\item Volver a 1., aplicar $DWT$ a la soluci\'on obtenida y repetir (ver \cite{Alice}).
\end{enumerate}
\newpage
\section{Resultados num\'ericos}\label{sec:hip-num}
En esta secci\'on se reproducir\'an algunos resultados obtenidos por Harten \cite{Harten1}. Para ello se aplicar\'a el algoritmo de multiresoluci\'on a la soluci\'on num\'erica de una ley de conservaci\'on, tomando como modelo la ecuaci\'on de Burgers (caso escalar y unidimensional)
\begin{equation}
u_t+(u^2/2)_x=0
\end{equation}
asociada a la condici\'on inicial
\begin{equation}\label{techo}
u(x,0)=\left\{\begin{array}{cc}
1,& \textrm{ si }|x|\leqslant 1/2\\
0,& \textrm{ si }1/2<|x|\leqslant 1.
\end{array}\right.
\end{equation}
Se utilizan condiciones peri\'odicas en $x=-1$ y $x=1$. Se opera hasta antes de que las discontinuidades alcancen las fronteras del dominio.

El primer objetivo es mostrar la relaci\'on existente entre la capacidad de compresi\'on de este m\'etodo de multiresoluci\'on y las propiedades de aproximaci\'on de las t\'ecnicas de reconstrucci\'on utilizadas. La localizaci\'on de los coeficientes de ondelette que est\'an por sobre una tolerancia prescrita, ayuda a visualizar esta conexi\'on. Recordar que los coeficientes de ondelette $d_j^k$ representan los errores cometidos en el proceso de predicci\'on y est\'an directamente relacionados a errores de interpolaci\'on, los cuales son peque\~nos en regiones de suavidad. En las proximidades de las singularidades el proceso de reconstrucci\'on podr{\'i}a conducir a regiones de exactitud pobre, por lo tanto, se examina el efecto del esquema de compresi\'on basado en la multiresoluci\'on.
Como una medida de la mejora en velocidad alcanzada mediante la utilizaci\'on del an\'alisis de multiresoluci\'on, se presenta la \emph{tasa de compresi\'on} o \emph{eficiencia} $\mu$ \cite{Bihari,Harten1} definida por $\mu=\frac{N_0}{N_0/2^L+|D^n|}$, donde $D^n$ es el conjunto de coeficientes de ondelette significativos, en todos los niveles de multiresoluci\'on, en el paso temporal $n$.

\begin{figure}[tp]
\begin{center}
\includegraphics[width=6in,height=3in]{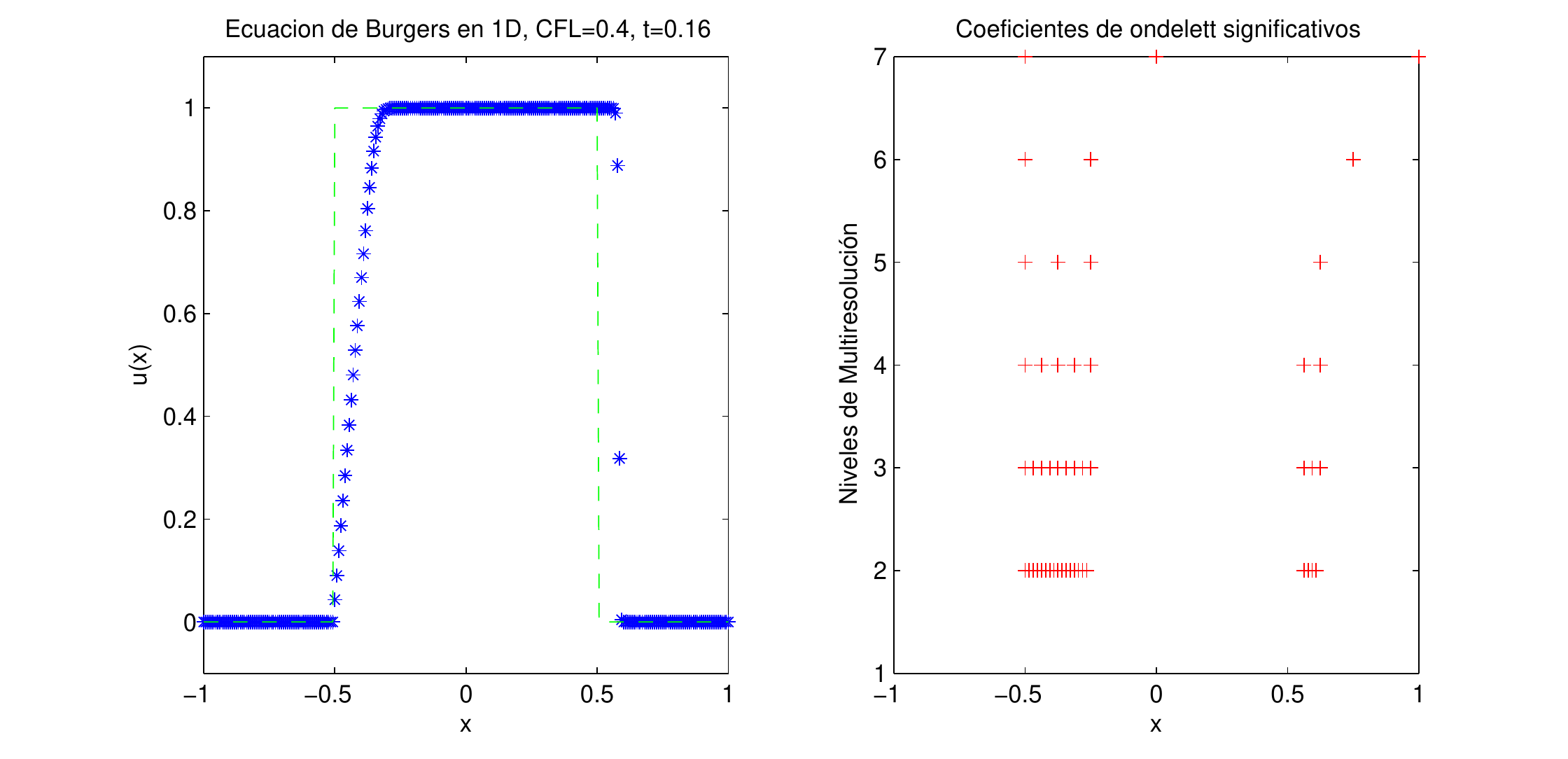}		
\caption[Soluci\'on num\'erica de multiresoluci\'on en el tiempo $t=0.16$ para la ec. de Burgers en 1D, $N_0=257$]{\small Izquierda: Soluci\'on inicial \emph{(rayas)} y soluci\'on num\'erica de multiresoluci\'on \emph{(asteriscos)} en el tiempo $t=0.16$ para la ec. de Burgers en 1D asociada a la condici\'on inicial (\ref{techo}), con $\varepsilon=10^{-5}$, $N_0=257$ y $L=7$. Derecha: Estructura de coeficientes de ondelette significativos correspondientes.}
\label{fig:hiper1}
\end{center}
\end{figure}

\begin{figure}[bp]
\begin{center}
\includegraphics[width=6in,height=3in]{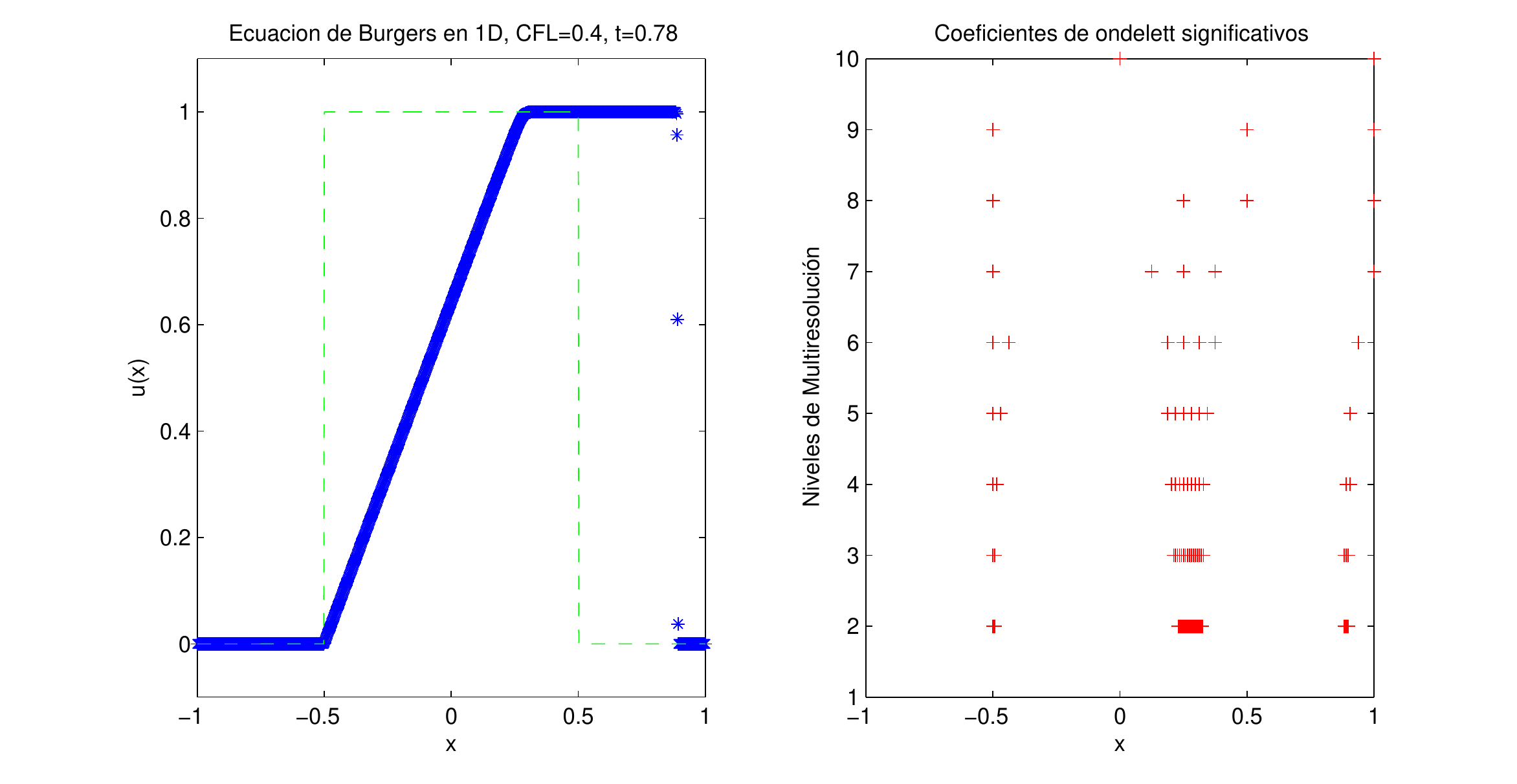}		
\caption[Soluci\'on num\'erica de multiresoluci\'on en el tiempo $t=0.78$ para la ec. de Burgers en 1D, $N_0=1025$]{\small Izquierda: Soluci\'on inicial \emph{(rayas)} y soluci\'on num\'erica de multiresoluci\'on \emph{(asteriscos)} en el tiempo $t=0.78$ para la ec. de Burgers en 1D asociada a la condici\'on inicial (\ref{techo}), con $N_0=1025$, $L=10$, $\varepsilon=10^{-3}$. Derecha: Estructura de coeficientes de ondelette significativos.}
\label{fig:hiper8}
\end{center}
\end{figure}

Las figuras \ref{fig:hiper1} a \ref{fig:hiper8} y las tablas \ref{tabla:hiper1} y \ref{tabla:hiper2} resumen el resultado de los test num\'ericos realizados. En cada figura, la parte izquierda representa a la soluci\'on num\'erica con asteriscos. La parte derecha muestra el conjunto de los coeficientes de ondelette significativos en el plano $x-k$, dibujando un + alrededor de cada $(x_j^k,k)$. Cada tabla muestra resultados de multiresoluci\'on para la soluci\'on num\'erica de la ecuaci\'on de Burgers para diferentes tiempos $t$. Se muestra la tasa de compresi\'on $\mu$, proporci\'on $V$ (entre el tiempo total de CPU de la soluci\'on num\'erica sin multiresoluci\'on y el tiempo total de CPU de la soluci\'on num\'erica con multiresoluci\'on) y los errores $e_p=\|u^n-u^n_{MR}\|_p$, $p=1,2,\infty$, donde
$$e_\infty=\max|u^n_j-u^n_{MR_j}|,\quad 1\leqslant j\leqslant N_0$$
y
$$e_p=\left(\frac{1}{N_0}\sum_{j=1}^{N_0}|u^n_j-u^n_{MR_j}|^p\right)^{1/p},\quad p=1,2.$$
En ambas tablas se ver\'a que el error obtenido es menor que la tolerancia prescrita. La norma $\mathcal{L_1}$ obtiene el menor error principalmente en funciones discontinuas \cite{Harten1}. Es importante precisar que los errores son calculados entre la representaci\'on puntual esparsa y la soluci\'on en malla fina, a\'un cuando la longitud de estos vectores no coincide (ya que existen posiciones en la malla fina para los cuales no corresponde ning\'un punto en la representaci\'on puntual esparsa). 

Se presentan los resultados correspondientes para el caso de $N_0=257$ puntos en la malla fina con $L=7$ niveles de multiresoluci\'on y el caso de $N_0=1025$ puntos en la malla fina con $L=10$ niveles de multiresoluci\'on. En ambos casos se utiliza una tolerancia de truncamiento $\varepsilon_k=\varepsilon/2^{L-k}$, condici\'on $CFL=0.5$, multiresoluci\'on con interpolador c\'uadr\'atico, flujos num\'ericos calculados mediante reconstrucci\'on ENO de segundo orden (ver secci\'on \ref{sec:eno2}) y evoluci\'on temporal Runge-Kutta de orden 2 (\ref{RK2}).

\begin{table}[h]
\begin{center}
\begin{tabular}{lccccc}
\hline
$t$ & $V$    & $\mu$  &            $e_1$   &        $e_2$     &     $e_\infty$     \\
\hline
0.16 &1.9330 & 19.7633 & 8.89$\times10^{-7}$&1.92$\times10^{-5}$&1.80$\times10^{-4}$\\
0.47 &1.8334 & 19.8122 & 1.99$\times10^{-6}$&3.15$\times10^{-5}$&6.14$\times10^{-5}$\\
0.62 &1.7696 & 19.4591 & 2.46$\times10^{-5}$&3.58$\times10^{-5}$&5.91$\times10^{-5}$\\
0.78 &1.6881 & 19.7633 & 2.92$\times10^{-5}$&3.96$\times10^{-5}$&5.77$\times10^{-5}$\\
\hline
\end{tabular}
\end{center}
\caption{Soluci\'on num\'erica de la Ecuaci\'on de Burgers, condici\'on inicial (\ref{techo}). Tolerancia prescrita $\varepsilon=10^{-5}$, $N_0=257$ puntos en la malla fina y $L=7$ niveles de multiresoluci\'on.}
\label{tabla:hiper1}
\end{table}

\begin{figure}[bp]
\begin{center}
\includegraphics[width=6in,height=3in]{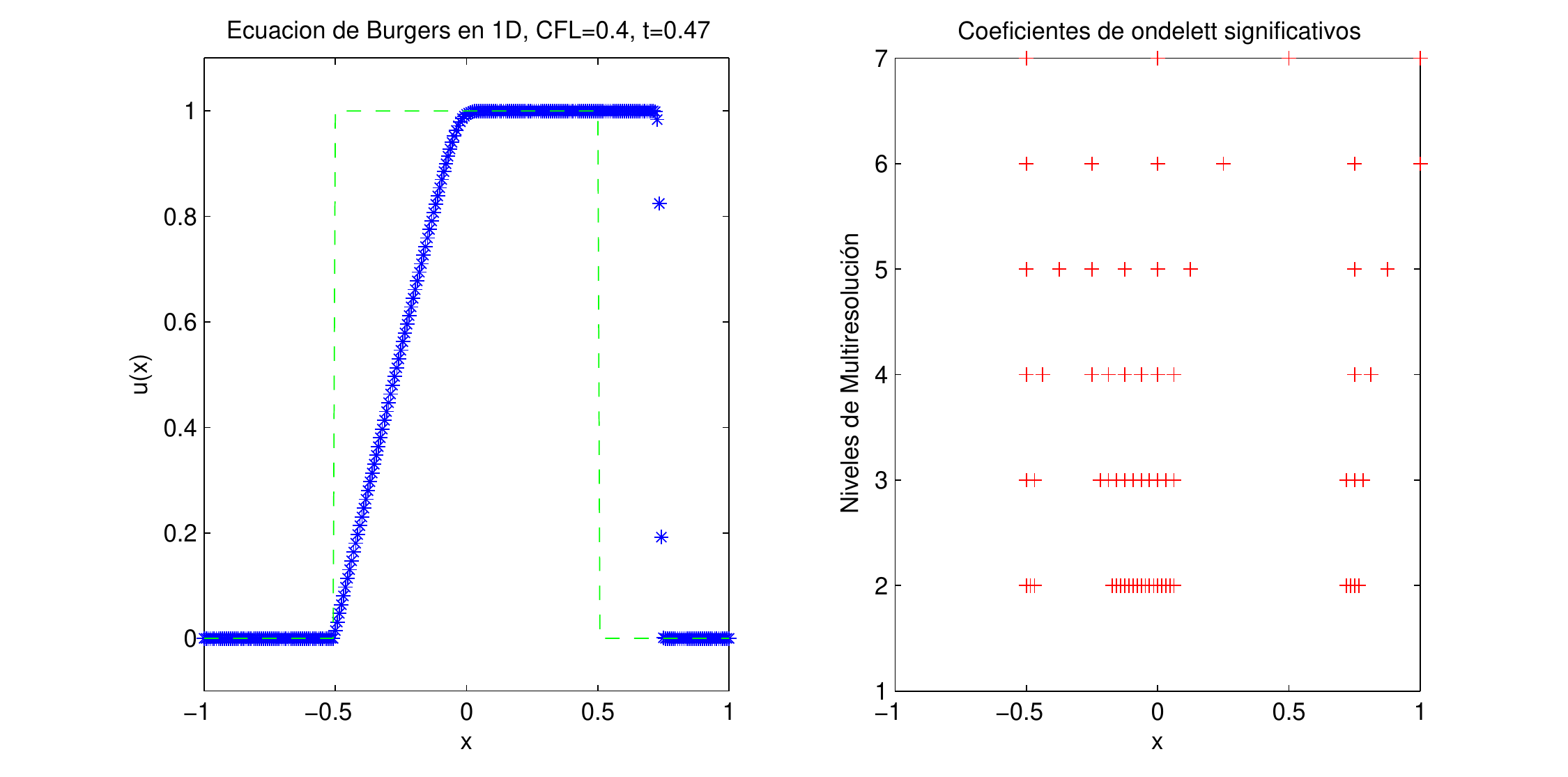}		
\caption[Soluci\'on num\'erica de multiresoluci\'on en el tiempo $t=0.47$ para la ec. de Burgers en 1D, $N_0=257$]{\small Izquierda: Soluci\'on inicial \emph{(rayas)} y soluci\'on num\'erica de multiresoluci\'on \emph{(asteriscos)} en el tiempo $t=0.47$ para la ec. de Burgers en 1D asociada a la condici\'on inicial (\ref{techo}), con $\varepsilon=10^{-5}$, $N_0=257$ y $L=7$. Derecha: Estructura de coeficientes de ondelette significativos.}
\label{fig:hiper2}
\end{center}
\end{figure}

\begin{figure}[tp]
\begin{center}
\includegraphics[width=6in,height=3in]{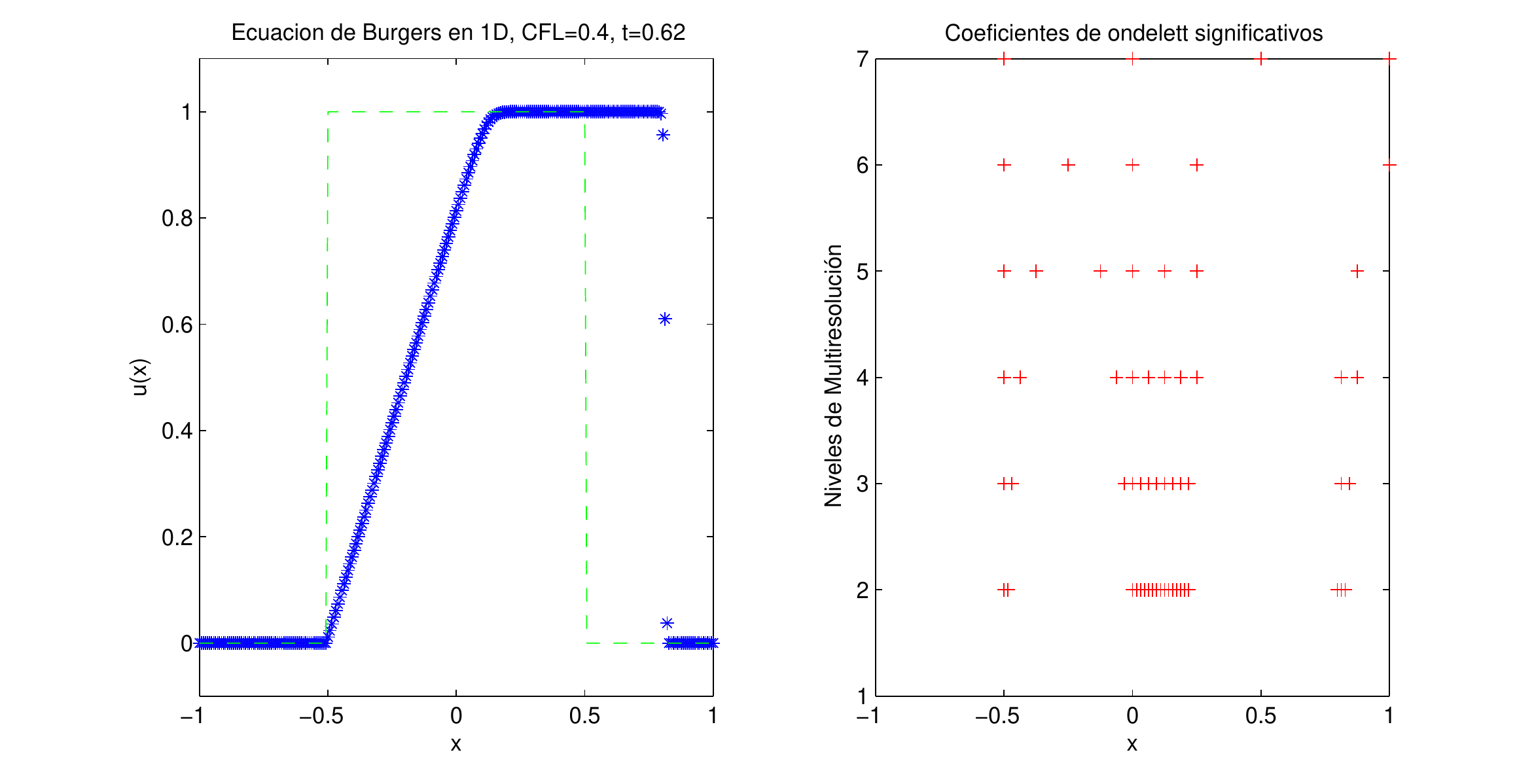}		
\caption[Soluci\'on num\'erica de multiresoluci\'on en el tiempo $t=0.62$ para la ec. de Burgers en 1D, $N_0=257$]{\small Izquierda: Soluci\'on inicial \emph{(rayas)} y soluci\'on num\'erica de multiresoluci\'on \emph{(asteriscos)} en el tiempo $t=0.62$ para la ec. de Burgers en 1D asociada a la condici\'on inicial (\ref{techo}), con $\varepsilon=10^{-5}$, $N_0=257$ y $L=7$. Derecha: Estructura de coeficientes de ondelette significativos.}
\label{fig:hiper3}
\end{center}
\end{figure}

\begin{figure}[bp]
\begin{center}
\includegraphics[width=6in,height=3in]{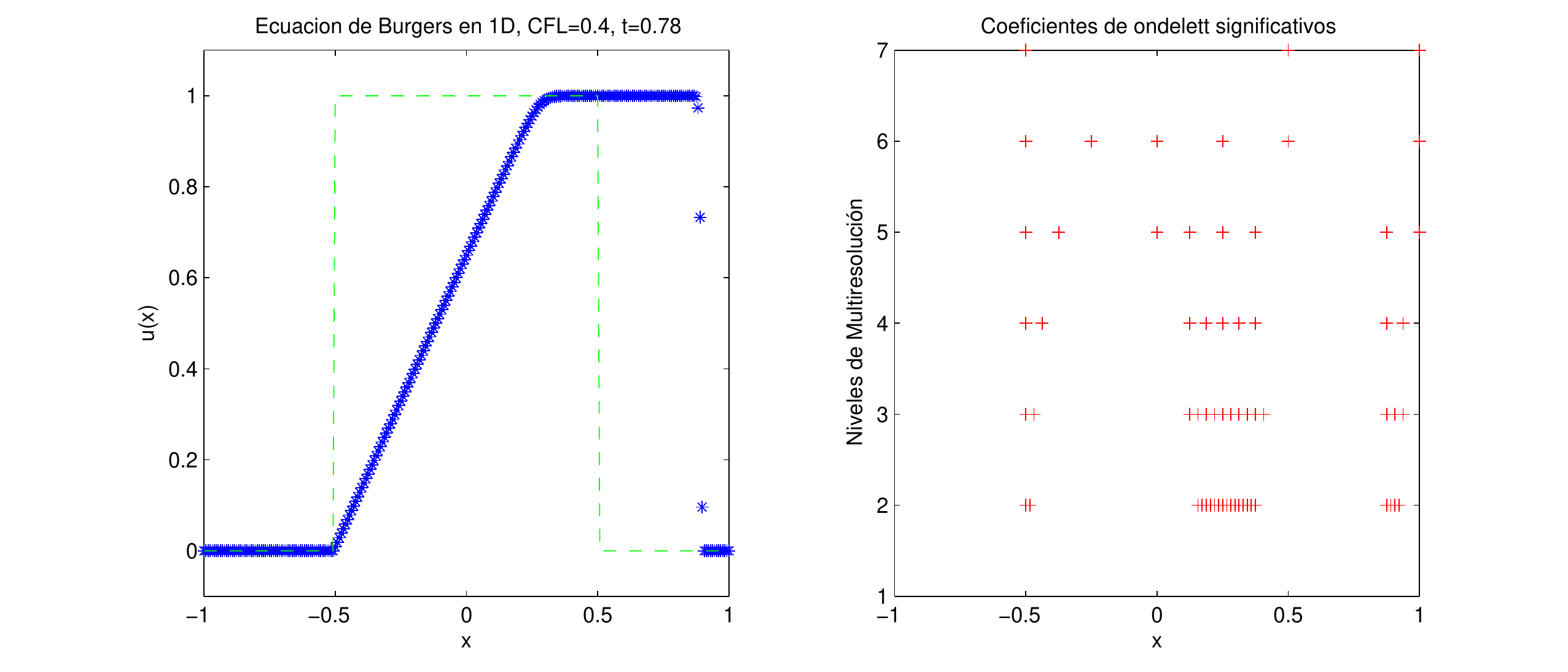}		
\caption[Soluci\'on num\'erica de multiresoluci\'on en el tiempo $t=0.78$ para la ec. de Burgers en 1D, $N_0=257$]{\small Izquierda: Soluci\'on inicial \emph{(rayas)} y soluci\'on num\'erica de multiresoluci\'on \emph{(asteriscos)} en el tiempo $t=0.78$ para la ec. de Burgers en 1D asociada a la condici\'on inicial (\ref{techo}), con $\varepsilon=10^{-5}$, $N_0=257$ y $L=7$. Derecha: Estructura de coeficientes de ondelette significativos.}
\label{fig:hiper4}
\end{center}
\end{figure}


\begin{figure}[tp]
\begin{center}
\includegraphics[width=6in,height=3in]{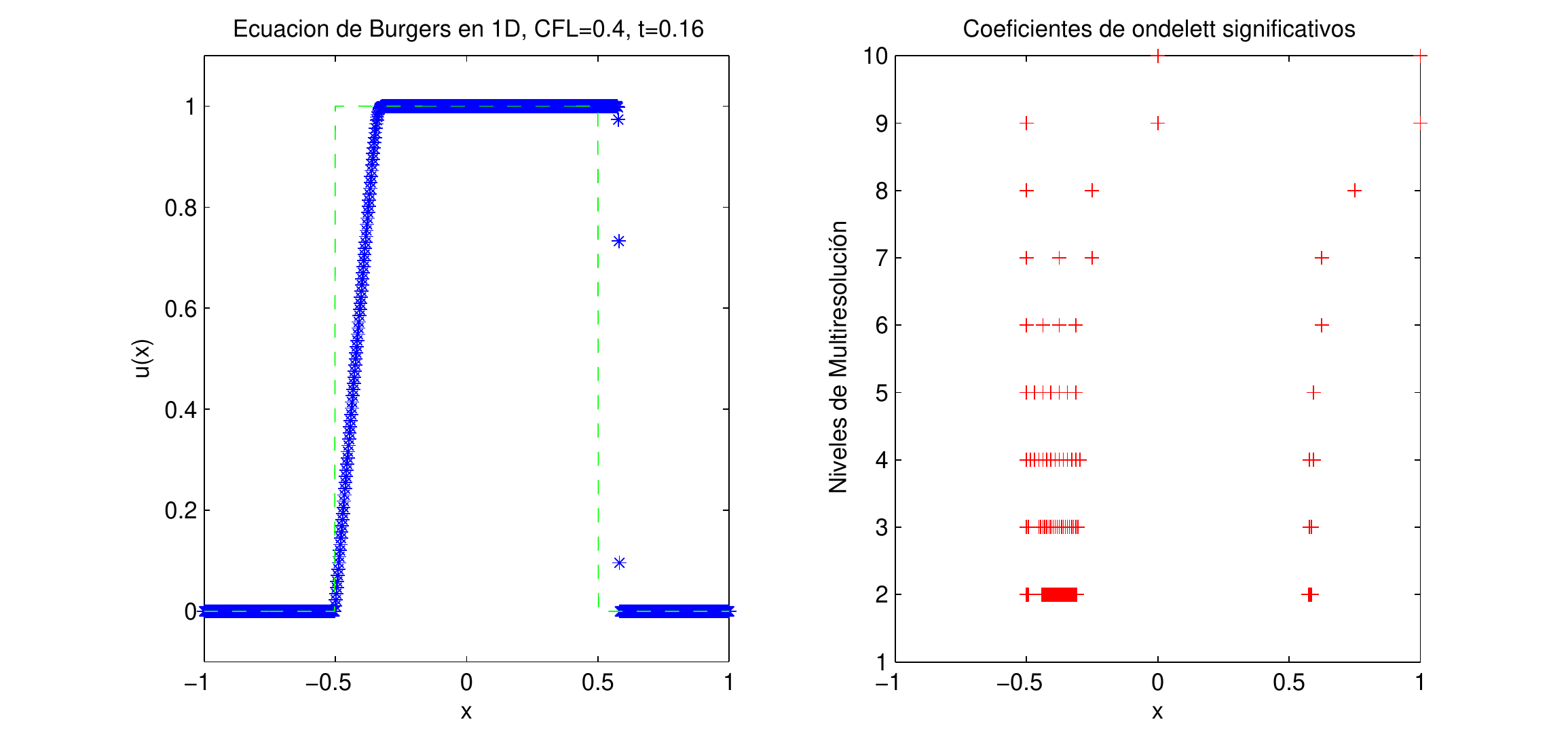}		
\caption[Soluci\'on num\'erica de multiresoluci\'on en el tiempo $t=0.16$ para la ec. de Burgers en 1D, $N_0=1025$]{\small Izquierda: Soluci\'on inicial \emph{(rayas)} y soluci\'on num\'erica de multiresoluci\'on \emph{(asteriscos)} en el tiempo $t=0.16$ para la ec. de Burgers en 1D asociada a la condici\'on inicial (\ref{techo}), con $N_0=1025$, $L=10$, $\varepsilon=10^{-3}$. Derecha: Estructura de coeficientes de ondelette significativos.}
\label{fig:hiper5}
\end{center}
\end{figure}

\begin{figure}[bp]
\begin{center}
\includegraphics[width=6in,height=3in]{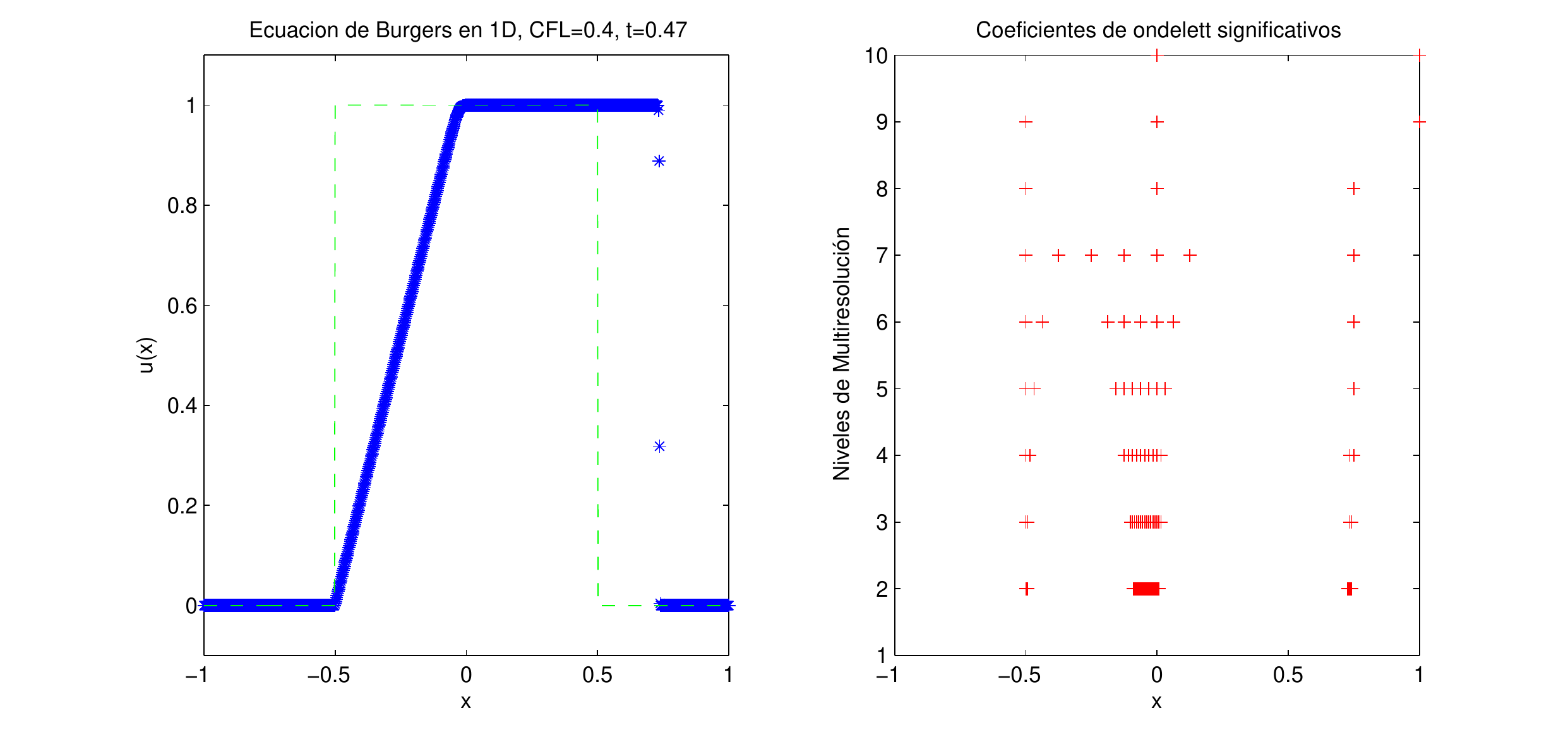}		
\caption[Soluci\'on num\'erica de multiresoluci\'on en el tiempo $t=0.47$ para la ec. de Burgers en 1D, $N_0=1025$]{\small Izquierda: Soluci\'on inicial \emph{(rayas)} y soluci\'on num\'erica de multiresoluci\'on \emph{(asteriscos)} en el tiempo $t=0.47$ para la ec. de Burgers en 1D asociada a la condici\'on inicial (\ref{techo}), con $N_0=1025$, $L=10$, $\varepsilon=10^{-3}$. Derecha: Estructura de coeficientes de ondelette significativos.}
\label{fig:hiper6}
\end{center}
\end{figure}

\begin{figure}[tp]
\begin{center}
\includegraphics[width=6in,height=3in]{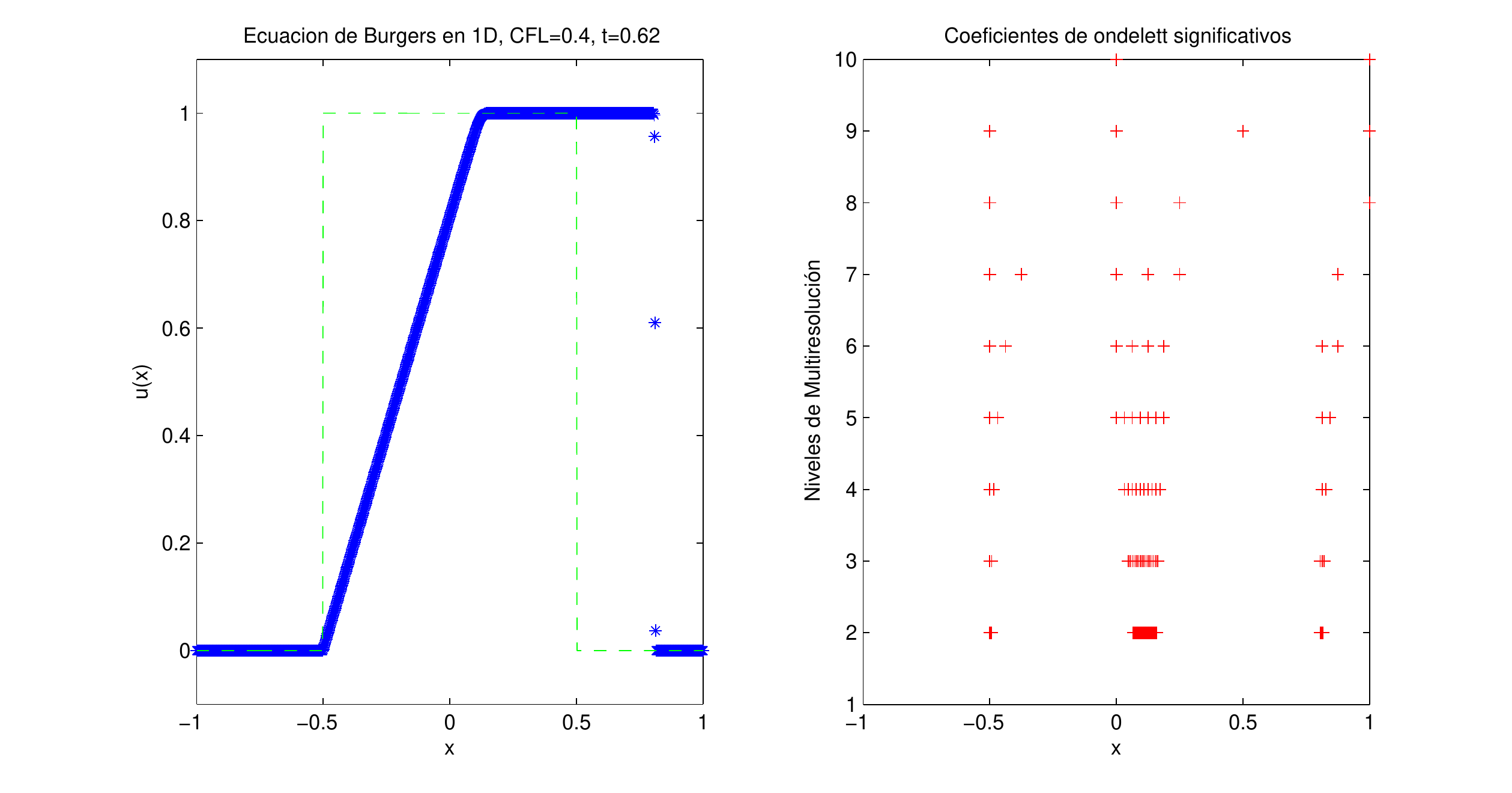}		
\caption[Soluci\'on num\'erica de multiresoluci\'on en el tiempo $t=0.62$ para la ec. de Burgers en 1D, $N_0=1025$]{\small Izquierda: Soluci\'on inicial \emph{(rayas)} y soluci\'on num\'erica de multiresoluci\'on \emph{(asteriscos)} en el tiempo $t=0.62$ para la ec. de Burgers en 1D asociada a la condici\'on inicial (\ref{techo}), con $N_0=1025$, $L=10$, $\varepsilon=10^{-3}$. Derecha: Estructura de coeficientes de ondelette significativos.}
\label{fig:hiper7}
\end{center}
\end{figure}

\begin{table}[h]
\begin{center}
\begin{tabular}{lccccc}
\hline
$t$ & $V$    & $\mu$  &            $e_1$   &        $e_2$     &     $e_\infty$     \\
\hline
0.16 &2.7872 & 53.9446 & 8.94$\times10^{-6}$&1.99$\times10^{-6}$&4.79$\times10^{-6}$\\
0.47 &2.5986 & 53.0172 & 2.09$\times10^{-5}$&3.01$\times10^{-6}$&5.59$\times10^{-6}$\\
0.62 &2.6170 & 53.5019 & 2.49$\times10^{-5}$&3.99$\times10^{-6}$&7.26$\times10^{-6}$\\
0.78 &2.5029 & 53.2874 & 2.97$\times10^{-5}$&4.26$\times10^{-6}$&1.88$\times10^{-5}$\\
\hline
\end{tabular}
\end{center}
\caption{Soluci\'on num\'erica de la Ecuaci\'on de Burgers, condici\'on inicial (\ref{techo}). Tolerancia prescrita $\varepsilon=10^{-3}$, $N_0=1025$ puntos en la malla fina y $L=10$ niveles de multiresoluci\'on.}
\label{tabla:hiper2}
\end{table}

\clearpage{\pagestyle{empty}\cleardoublepage}
\chapter{Caso parab\'olico}
En este cap{\'i}tulo se aplicar\'a el algoritmo de multiresoluci\'on a ecuaciones parab\'olicas. Se reproducir\'an los experimentos num\'ericos realizados por Roussel \emph{et al.} \cite{Sch}, Bihari \cite{Bihari}, Liandrat y Tchamitchian\cite{Liandrat}.
\section{M\'etodo num\'erico}
A continuaci\'on se presenta un m\'etodo general de vol\'umenes finitos para ecuaciones hiperb\'olicas, incluyendo la descripci\'on de los esquemas utilizados para la discretizaci\'on espacial y evoluci\'on temporal \cite{Sch}.
\subsection{Leyes de conservaci\'on parab\'olicas}
Se considera el problema de valores iniciales para una ecuaci\'on parab\'olica en $(x,t)\in \Omega\times[0,\infty[$, $\Omega\subset \RR^d$ de la forma
\begin{equation}\label{parab}
\begin{split}
\frac{\partial u}{\partial t}+\nabla\cdot F(u,\nabla u)&=S(u),\\
u(x,0)=u_0(x)&
\end{split}
\end{equation}
asociada a condiciones de borde apropiadas.

Se considerar\'a la restricci\'on al caso en que el flujo difusivo se define por un operador gradiente, suponiendo difusividad constante $\nu>0$, es decir,
\begin{equation*}
F(u,\nabla u)=f(u)-\nu\nabla u.
\end{equation*}
Para la ecuaci\'on de convecci\'on-difusi\'on en 1D, se tiene ($c>0$)
\begin{eqnarray*}
f(u)&=&cu,\\
S(u)&=&0.
\end{eqnarray*}
En el caso de la ecuaci\'on viscosa de Burgers en 1D, se tiene
\begin{eqnarray*}
f(u)&=&\frac{u^2}{2},\\
S(u)&=&0,
\end{eqnarray*}
y para la ecuaci\'on de reacci\'on-difusi\'on ($\alpha>0,\beta>0$),
\begin{eqnarray*}
f(u)&=&0,\\
S(u)&=&\frac{\beta^2}{2}(1-u)\exp\frac{\beta(1-u)}{\alpha(1-u)-1}.
\end{eqnarray*}

Se define el \emph{t\'ermino fuente y de divergencia} por
$$
\mathcal{D}(u,\nabla u)=-\nabla\cdot F(u,\nabla u)+S(u).
$$
Luego (\ref{parab}) puede escribirse como
\begin{equation}\label{parab2}
\frac{\partial u}{\partial t}=\mathcal{D}(u,\nabla u).
\end{equation}
\subsection{Discretizaci\'on}
Para discretizar (\ref{parab2}), se utiliza una formulaci\'on de vol\'umenes finitos en la forma conservativa est\'andar. En el caso general, consid\'erese el dominio computacional $\Omega$ y una partici\'on de \'el en vol\'umenes de control $(\Omega_i)_{i\in\Lambda}$, $\Lambda=\{1,\ldots,i_{max}\}$. Se denota entonces por $\bar{q}_i(t)$ al promedio de cierta cantidad $q$ sobre $\Omega_i$ en el instante $t$,
\begin{equation}
\bar{q}_i(t)=\frac{1}{|\Omega_i|}\int_{\Omega_i}q(x,t)dx.
\end{equation}
Integrando (\ref{parab2}) y promediando sobre $\Omega_i$,
\begin{equation}
\frac{1}{|\Omega_i|}\int_{\Omega_i}\frac{\partial u}{\partial t}(x,t)dx=\frac{1}{|\Omega_i|}\int_{\Omega_i}\mathcal{D}(u(x,t),\nabla u(x,t))dx.
\end{equation}
Luego
\begin{equation}\label{parab3}
\frac{\partial \bar{u}_i}{\partial t}(t)=\bar{\mathcal{D}}_i(t).
\end{equation}
Si se aplica el teorema de la divergencia, se obtiene
\begin{equation}\label{parab4}
\bar{\mathcal{D}}_i=-\frac{1}{|\Omega_i|}\int_{\partial\Omega_i}F(u,\nabla u)\cdot \sigma_i(x)\, dx + \bar{S}_i(t),
\end{equation}
donde $\sigma_i(x)$ es el vector normal hacia $\Omega_i$. La conservatividad en el c\'alculo del flujo se garantiza si y s\'olo si, para dos vol\'umenes de control adyacentes $\Omega_{i_1}$ y $\Omega_{i_2}$, el flujo que va de $\Omega_{i_1}$ a $\Omega_{i_2}$ se equilibra con el flujo que va de $\Omega_{i_2}$ a $\Omega_{i_1}$.
\subsection{Integraci\'on temporal}
Notar que se est\'a frente a un proceso de discretizaci\'on en dos etapas, debido a la adaptatividad de la discretizaci\'on espacial. Primero se discretiza s\'olo en tiempo, y luego en espacio. Esto conduce a las ya mencionadas ecuaciones semi-discretas (ver secci\'on \ref{Etemp}). La discretizaci\'on puede hacerse utilizando un m\'etodo num\'erico est\'andar para sistemas de ecuaciones diferenciales ordinarias. Este mecanismo es particularmente ventajoso en el desarrollo de m\'etodos con orden de precisi\'on mayor a dos, ya que permite alcanzar de forma relativamente sencilla la misma precisi\'on espacial y temporal. Los experimentos realizados en \cite{ENO} indican que las formulaciones semi-discretas con discretizaci\'on temporal Runge-Kutta TVD desarrollados por Shu y Osher no generan oscilaciones para $CFL \leqslant 0.5$ aproximadamente, y son \'optimas en el sentido de que permiten la mayor $CFL$ para esquemas expl{\'i}citos, $CFL=1$ \cite{ENO,ENO2,Alice}.

Se utilizar\'a entonces un m\'etodo expl{\'i}cito Runge-Kutta TVD de segundo orden que en este caso se expresa por
\begin{eqnarray}
\bar{u}_i^{n+1/2}&=&\bar{u}_i^n+\Delta t\bar{\mathcal{D}}_i^n,\nonumber\\
\bar{u}_i^{n+1}&=&\frac{1}{2}\left[\bar{u}_i^n+\bar{u}_i^{n+1/2}+\Delta t\bar{\mathcal{D}}_i^{n+1/2}\right].\label{RK2}
\end{eqnarray}
Notar que (\ref{RK2}) tambi\'en se conoce como \emph{M\'etodo de Heun} \cite{Bihari}.

Si se denota por $\bar{u}^n$ al vector $(\bar{u}_i^n)_{i\in\Lambda}$, entonces el operador de evoluci\'on temporal discreto $\bar{E}(\Delta t)$ est\'a definido por
\begin{equation}\label{time2}
\bar{u}^{n+1}=\bar{E}(\Delta t)\cdot \bar{u}^n,
\end{equation}
donde
\begin{equation}\label{time3}
\bar{E}(\Delta t)=\mathbf{I}+\frac{\Delta t}{2}\left[\bar{\mathcal{D}}+\bar{\mathcal{D}}(\mathbf{I}+\Delta t\bar{\mathcal{D}})\right].
\end{equation}
La discretizaci\'on del operador $\bar{\mathcal{D}}$ se describe en la siguiente secci\'on.

\subsection{Flujo num\'erico}\label{section:fn}
Consid\'erese ahora un tiempo fijo $t^n$. Para el caso unidimensional general, $\Omega_i$ es el intervalo $[x_{i-1/2},x_{i+1/2}]$ de longitud $\Delta x_i=x_{i+1/2}-x_{i-1/2}$. Mediante una discretizaci\'on de vol\'umenes finitos est\'andar, la ecuaci\'on (\ref{parab4}) puede escribirse como
\begin{equation}\label{flux1}
\bar{\mathcal{D}}_i=-\frac{1}{\Delta x_i}\left(\bar{F}_{i+\frac{1}{2}}-\bar{F}_{i-\frac{1}{2}}\right)+\bar{S}_i,
\end{equation}
donde
\begin{equation}\label{flux2}
\bar{F}_{i+\frac{1}{2}}=f^R\left(\bar{u}^-_{i+\frac{1}{2}},\bar{u}^+_{i+\frac{1}{2}}\right)-\nu\frac{\bar{u}_{i+1}-\bar{u}_i}{\Delta x_{i+\frac{1}{2}}},
\end{equation}
con $\Delta x_{i+\frac{1}{2}}=\frac{1}{2}(\Delta x_i+\Delta x_{i+1})$. El t\'ermino $f^R$ denota, para la parte advectiva, la soluci\'on aproximada de Roe para el problema de Riemann \cite{GR}, dados los estados de derecha e izquierda de $u$. La versi\'on escalar correspondiente es
\begin{equation}\label{roe}
f^R(u^-,u^+)=\frac{1}{2}[f(u^-)+f(u^+)-|a(u^-,u^+)|(u^+-u^-)],
\end{equation}
donde
\begin{equation*}
a(u^-,u^+)=\left\{\begin{array}{cl}
\frac{f(u^+)-f(u^-)}{u^+-u^-}, & \textrm{ si } u^+\neq u^-,\\
f'(u^+),  & \textrm{ si } u^+= u^-.
\end{array}\right.
\end{equation*}
Los valores de izquierda y derecha $\bar{u}^-_{i+\frac{1}{2}}$ y $\bar{u}^+_{i+\frac{1}{2}}$, respectivamente, son obtenidos mediante interpolaci\'on ENO de segundo orden (ver secci\'on \ref{sec:eno2}).

Notar de (\ref{flux2}) que los t\'erminos advectivo y difusivo son aproximados de diferente forma. Para la parte advectiva, se utiliza el esquema de Roe cl\'asico con una interpolaci\'on ENO de segundo orden; mientras que para la parte difusiva, se escoge un esquema centrado en $\bar{u}_i$ de segundo orden.

En \cite{Bihari} se prueba que el esquema global resultante, que es no lineal,
\small\begin{equation}\label{flux3}
\bar{\mathcal{D}}_i=-\frac{1}{\Delta x_i}\left(f^R\left(\bar{u}^-_{i+\frac{1}{2}},\bar{u}^+_{i+\frac{1}{2}}\right)-f^R\left(\bar{u}^-_{i-\frac{1}{2}},\bar{u}^+_{i-\frac{1}{2}}\right)-\nu\frac{\bar{u}_{i+1}-2\bar{u}_i+\bar{u}_{i-1}}{\Delta x_{i+\frac{1}{2}}} \right)+\bar{S}_i,
\end{equation}
\normalsize
es de segundo orden (en espacio).

El t\'ermino fuente es aproximado por $\bar{S}_i\approx S(\bar{u}_i)$. Para un t\'ermino fuente no lineal, esta elecci\'on tambi\'en implica una precisi\'on de orden dos \cite{Sch}.
\subsection{Reconstrucci\'on ENO de segundo orden}\label{sec:eno2}
Para obtener los valores de la funci\'on $u$ en las fronteras de los vol\'umenes de control, se utiliza una reconstrucci\'on lineal a trozos de $u$ a partir de los valores de las medias en celda. Es decir, los t\'erminos de izquierda y derecha $\bar{u}^-_{i+\frac{1}{2}}$ y $\bar{u}^+_{i+\frac{1}{2}}$, respectivamente, son obtenidos mediante interpolaci\'on ENO de segundo orden \cite{Alice,ENO,ENO2,Sch}. Este tipo de m\'etodos utiliza una construcci\'on adaptativa del est\'encil a fin de evitar la generaci\'on de oscilaciones esp\'ureas cerca de las discontinuidades. Se puede generar oscilaciones, pero del orden del error local de truncamiento en la parte suave de la soluci\'on. En este caso particular, se tiene
\begin{eqnarray}
\bar{u}^-_{i+\frac{1}{2}}&=&\bar{u}_i+\frac{1}{2}M\left(\bar{u}_{i+1}-\bar{u}_i,\bar{u}_i- \bar{u}_{i-1}\right),\\
\bar{u}^+_{i+\frac{1}{2}}&=&\bar{u}_{i+1}+\frac{1}{2}M\left(\bar{u}_{i+2}-\bar{u}_{i+1},\bar{u}_{i+1}- \bar{u}_i\right),
\end{eqnarray}
donde $M$ es el limitador \emph{Min-Mod}, que escoge la pendiente m{\'i}nima entre los extremos izquierdo y derecho, es decir,
\begin{equation*}
M(a,b)=\left\{\begin{array}{ll}
a,& \textrm{ si } |a|\leqslant |b|,\\
b,& \textrm{ si } |a|>|b|.
\end{array}\right.
\end{equation*}

Notar que (\ref{flux2}) es la forma \emph{semi-discreta} de (\ref{time2}). (\ref{flux2}) se resuelve utilizando una actualizaci\'on temporal Runge-Kutta de segundo orden; por lo tanto se obtiene un esquema de segundo orden tanto en tiempo como en espacio.

Mediante un argumento de producto tensorial, puede llevarse a cabo la extensi\'on natural de la reconstrucci\'on a 2D y 3D en geometr{\'i}as cartesianas \cite{Sch}.

\subsection{Soluci\'on exacta de la onda viajera}
Para formar una idea cualitativa de la estructura del choque, consid\'erese la soluci\'on  $u(x,t)=u(\psi)$, $\psi=(x-st)/\nu$ del problema de la onda viajera
\begin{eqnarray}
u_t+f(u)_x&=&\nu u_{xx},\\
u(x,0)&=&\left\{\begin{array}{ll}
u_L,&\textrm{ si } x<0,\\
u_R<u_L,&\textrm{ si } x\geqslant 0.
\end{array}\right.
\end{eqnarray}
La ecuaci\'on diferencial ordinaria resultante en $\psi$ puede integrarse para obtener
\begin{equation}\label{exact2}
-su+f(u)+c=u',
\end{equation}
donde $s$ y $c$ pueden ser determinadas de las ``condiciones de borde''
$$\lim_{\psi\to-\infty}u(\psi)=u_L,\quad  \lim_{\psi\to\infty}u(\psi)=u_R$$
como sigue
\begin{eqnarray}
c&=&su_L-f(u_L),\\
s&=&\frac{f(u_R)-f(u_L)}{u_R-u_L},
\end{eqnarray}
donde la velocidad de la onda $s$ puede ser identificada como la velocidad del choque (asume la misma expresi\'on que en el caso puramente hiperb\'olico). Una nueva integraci\'on de (\ref{exact2}) entrega una f\'ormula impl{\'i}cita para $u$:
\begin{equation}
\int\frac{du}{f(u)-su+c}=\psi+c_1.
\end{equation}
En el caso particular de la ecuaci\'on de Burgers viscosa, es decir, $f(u)=\frac{1}{2}u^2$, se obtiene
\begin{equation}\label{exact3}
u(\psi)=u_L\tanh\frac{u_L-u_R}{4}\psi.
\end{equation}
Ver detalles en \cite{Bihari}.
\subsection{Estabilidad num\'erica}\label{estab1}
Como el paso temporal es el mismo para todas las escalas de multiresoluci\'on, la condici\'on de estabilidad es la correspondiente al esquema de vol\'umenes finitos en la malla fina. Si denotamos por $\Delta x$ al menor paso espacial, el n\'umero CFL $\sigma$ est\'a dado por
\begin{equation}
\sigma=u_{\max}\frac{\Delta t}{\Delta x}.
\end{equation}
Para el caso lineal (ecuaci\'on de convecci\'on-difusi\'on), si $c$ es la velocidad,
\begin{equation}\label{defisig}
\sigma=\frac{c\Delta t}{\Delta x}
\end{equation}
y el n\'umero de Reynolds \emph{Re} est\'a dado por
\begin{equation}\label{defire}
Re=\frac{c\Delta x}{\nu}.
\end{equation}
En \cite{Bihari} y \cite{GR} se muestra que una condici\'on suficiente para asegurar la estabilidad del esquema de vol\'umenes finitos es
\begin{equation}
\sigma\leqslant\min\left(\frac{Re}{2},\frac{6}{Re}\right).
\end{equation}
A\'un m\'as, una condici\'on suficiente \emph{para que el esquema sea TVD} (ver ap\'endice \ref{anexob}), es
\begin{equation}\label{sigma}
\sigma\leqslant\frac{Re}{Re+4}.
\end{equation}
La mayor ventaja de utilizar un esquema expl{\'i}cito para el t\'ermino difusivo, es que no se necesita resolver un sistema lineal. Sin embargo, esto generalmente implica que $\Delta t=O(\Delta x^2)$. S\'olo para el caso $Re>>1$ se puede esperar $\Delta t=O(\Delta x)$ \cite{FV}.

A continuaci\'on se analizar\'a un \emph{Esquema de multiresoluci\'on conservativo completamente adaptativo} dise\~nado por Roussel \emph{et al.} \cite{Sch}.

\subsection{\'Arbol graduado din\'amico}
El principio del an\'alisis de multiresoluci\'on es representar un conjunto de datos dados en malla fina como valores en la malla m\'as gruesa y un conjunto de detalles a diferentes escalas de mallas anidadas. Se propone organizar la estructura de datos como un \emph{\'arbol graduado} din\'amico, que posee una capacidad mayor de compresi\'on que la estructura MORSE o SPARSE de la representaci\'on puntual esparsa.

En la terminolog{\'i}a de las ondelettes, una estructura de \'arbol graduado corresponde a una aproximaci\'on adaptativa en la que est\'a garantizada la conectividad para la estructura de \'arbol.

\begin{figure}[ht]
	\begin{center}
	\includegraphics[width=5.2in,height=1.22in]{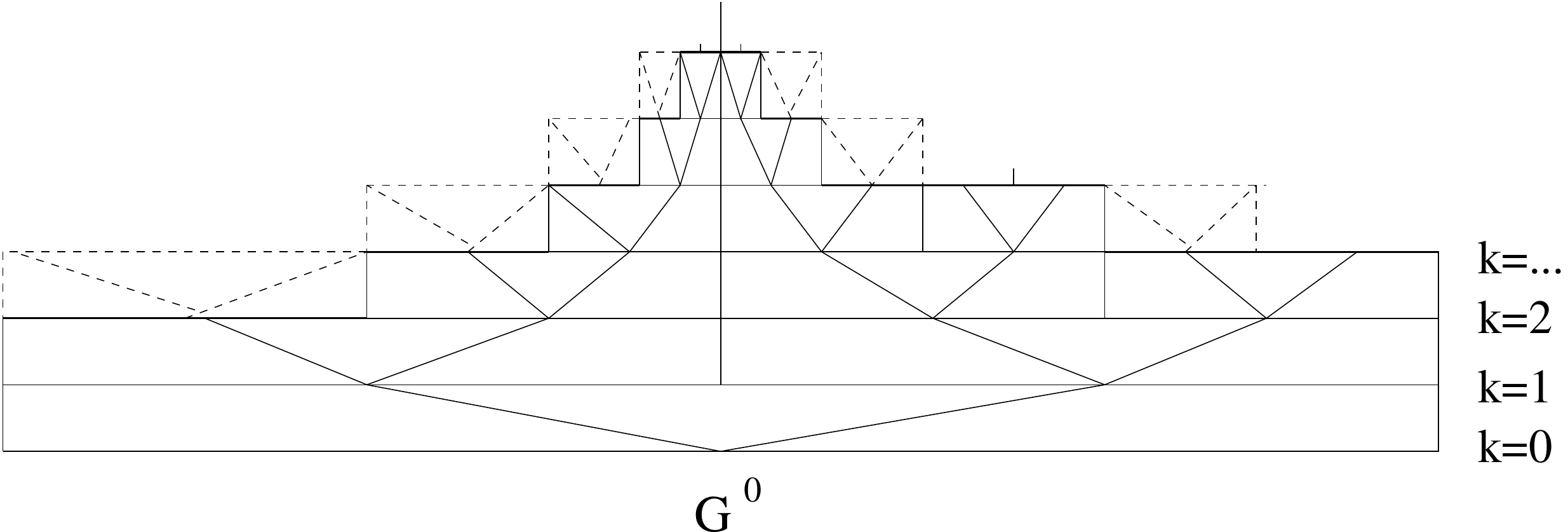}
	\caption{Estructura de datos tipo \'arbol graduado din\'amico unidimensional.}
	\label{fig:arbol}
	\end{center}
\end{figure}

Para definir la estructura de \'arbol, se introduce la terminolog{\'i}a utilizada por Cohen \cite{Cohen,Sch} :
\begin{itemize}
\item La \emph{ra{\'i}z} es la base del \'arbol.
\item Un \emph{nodo} es un elemento del \'arbol. Cada volumen de control ser\'a considerado un nodo.
\item Un nodo \emph{padre} tiene 2 nodos \emph{hijos}; los nodos hijos de un mismo nodo padre son llamados \emph{hermanos}.
\item Un nodo tiene vecinos cercanos en cada direcci\'on, llamados \emph{primos cercanos}. Los nodos hermanos son tambi\'en considerados como primos cercanos.
\item Un nodo es llamado \emph{hoja} cuando no tiene hijos.
\item Para calcular los flujos entrantes y salientes de cada hoja, se necesitan los primos cercanos. Cuando alguno de ellos no existe, se crea una \emph{hoja virtual} (representada por rayas en la figura \ref{fig:arbol}). Esta no se considera como un nodo existente, sino s\'olo se utiliza para calcular flujos.
\end{itemize}
Un \emph{\'arbol din\'amico} es un \'arbol que cambia en el tiempo. Si es necesario, algunos nodos pueden ser agregados o quitados. Para permanecer \emph{graduado}, el \'arbol debe respetar las condiciones siguientes:
\begin{itemize}
\item Cuando un hijo es creado, todos sus hermanos son creados en el mismo tiempo;
\item Un nodo tiene siempre dos primos cercanos en cada direcci\'on. Si no existe, debe ser creado como hoja virtual.
\item Un nodo puede ser quitado s\'olo si son quitados todos sus hermanos y s\'olo si no es el primo cercano de un nodo existente.
\end{itemize}

\subsection{An\'alisis del error}
El error global entre los valores puntuales de la soluci\'on exacta en el nivel $L$, $u^L_{ex}$, y los valores de la soluci\'on num\'erica por multiresoluci\'on con un nivel m\'aximo $L$, $u^L_{MR}$, puede ser descompuesto en dos errores \cite{Harten1,Sch}:
\begin{equation}\label{err1}
\|u^L_{ex}-u^L_{MR}\|\leqslant\|u^L_{ex}-u^L_{FV}\|+\|u^L_{FV}-u^L_{MR}\|,
\end{equation}
donde $\|\cdot\|$ es la norma $\mathcal{L}^1$, $\mathcal{L}^2$, o $\mathcal{L}^\infty$. El primer error del lado derecho de (\ref{err1}), llamado \emph{error de discretizaci\'on}, es el error del esquema de vol\'umenes finitos en malla fina, para un nivel m\'aximo $L$. Puede ser acotado por
\begin{equation}
\|u^L_{ex}-u^L_{FV}\| \leqslant C2^{-\xi L},\quad C>0,
\end{equation}
donde $\xi$ es el orden de convergencia del esquema de vol\'umenes finitos. En este caso, se utilizar\'an esquemas de segundo orden (en tiempo y espacio). Luego $\xi=2$.

El segundo error del lado derecho de (\ref{err1}) es llamado \emph{error de perturbaci\'on}. En \cite{Cohen} se prueba que si los detalles en un nivel de multiresoluci\'on $k$ son truncados bajo cierta tolerancia prescrita $\varepsilon_k$, si el operador de evoluci\'on temporal discreto $\mathbf{\bar{E}}$ es contractivo en la norma correspondiente, y si la tolerancia prescrita en el nivel $k$ es
$$ \varepsilon_k=2^{(k-L)}\varepsilon,$$
entonces la diferencia entre la soluci\'on por vol\'umenes finitos en la malla fina y la soluci\'on obtenida mediante un algoritmo de multiresoluci\'on, se acumula en el tiempo y verifica
\begin{equation}
\|u^L_{MR}-u^L_{FV}\| \leqslant Cn\varepsilon,\quad C>0,
\end{equation}
donde $n$ es el n\'umero de pasos temporales. Considerando un tiempo fijo $T=n\Delta t$, esto es
\begin{equation*}
\|u^L_{MR}-u^L_{FV}\| \leqslant C\frac{T}{\Delta t}\varepsilon,\quad C>0.
\end{equation*}
Para la ecuaci\'on lineal de convecci\'on-difusi\'on, de (\ref{sigma}), el paso temporal $\Delta t$ debe verificar
\begin{equation*}
\Delta t\leqslant\frac{\Delta x^2}{4\nu+c\Delta x}.
\end{equation*}
Si se denota por $X$ a la longitud del dominio, $\Delta x$ al paso espacial en la malla fina, y en el caso de que la ra{\'i}z del \'arbol graduado contenga s\'olo un nodo, se tiene $\Delta x=X2^{-L}$. Luego
$$\Delta t=C\frac{(\Delta x)^2}{4\nu+c\Delta x}=C\frac{X2^{-2L}}{4\nu+cX2^{-L}},\quad 0<C<1.$$
Si se quiere que el error de perturbaci\'on sea del mismo orden que el error de discretizaci\'on,
$$\varepsilon/\Delta t\propto2^{-\xi L}.$$
Por lo tanto,
$$\varepsilon2^{2L}(4\nu+ cX2^{-L})\propto 2^{-\xi L},$$
y si se define el n\'umero de Peclet como $Pe=\frac{cX}{\nu}$,
\begin{equation}\label{propto}
\varepsilon\propto \frac{2^{-(\xi+1)L}}{Pe+2^{(L+2)}}.
\end{equation}
Para el caso inv{\'i}scido ($Pe\to\infty$), \ref{propto} es equivalente a los resultados obtenidos en \cite{Cohen}:
$$\varepsilon\propto 2^{-(\xi+1)L}.$$
Con esto, elegiremos una \emph{tolerancia de referencia}:
\begin{equation}
\varepsilon_R=C\frac{2^{-(\xi+1)L}}{Pe+2^{(L+2)}}.
\end{equation}

\subsection{C\'alculo del flujo conservativo}
\begin{figure}[ht]
	\begin{center}
	\includegraphics[height=1in]{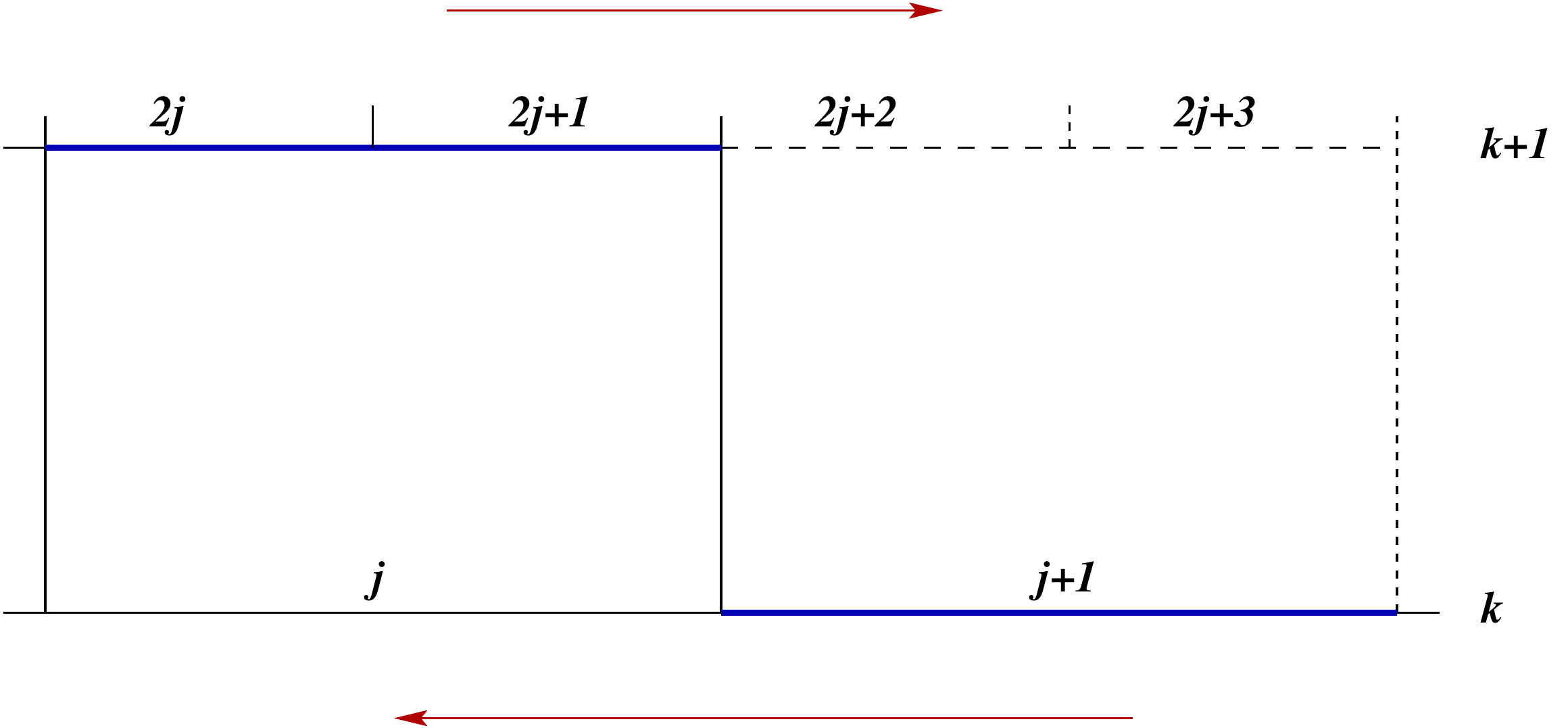}
	\caption{Flujo saliente y entrante para dos niveles diferentes.}
	\label{fig:flux2}
	\end{center}
\end{figure}

Consid\'erese una hoja $\Omega_{k+1,2j+1}$ con primos virtuales $\Omega_{k+1,2j+2}$ y $\Omega_{k+1,2j+3}$ a la derecha. Su padre $\Omega_{k,j+1}$ es una hoja.
Como se ve en la figura \ref{fig:flux2}, el flujo que sale de $\Omega_{k+1,2j+1}$ hacia la derecha $F_{k+1,2j+1\to k+1,2j+2}$ no est\'a en equilibrio con el flujo que sale de $\Omega_{k,j+1}$ hacia la izquierda  $F_{k,j+1\to k,j}$. Es posible calcular directamente los flujos que salen de $\Omega_{k+1,2j+1}$ hacia $\Omega_{k,j+1}$ o pueden calcularse s\'olo los flujos en el nivel $k+1$ y para determinar el flujo entrante a la hoja en el nivel $k$, \'este ser\'a igual a la suma de los flujos salientes de las hojas en el nivel $k+1$.

Esta elecci\'on asegura una conservatividad estricta en el c\'alculo de los flujos entre vol\'umenes de control de niveles diferentes, sin un aumento significativo de las evaluaciones (generalmente costosas) de los flujos.
\subsection{Implementaci\'on del algoritmo}
A continuaci\'on se presenta la estrategia a seguir por el algoritmo. En primer lugar, dependiendo de la condici\'on inicial dada, se crea un \'arbol graduado inicial. Luego se realiza la evoluci\'on temporal sobre las hojas y finalmente se actualiza el \'arbol graduado.
\small
\begin{itemize}
\item Inicializaci\'on de par\'ametros: tiempo de simulaci\'on, tama\~no del dominio, niveles de multiresoluci\'on, n\'umero de puntos en la malla fina, condici\'on $CFL$, etc.
\item Creaci\'on de la estructura de \'arbol graduado inicial: C\'alculo de detalles mediante transformada de multiresoluci\'on, obtenci\'on de la representaci\'on puntual esparsa.
\item Evoluci\'on temporal: C\'alculo del operador discreto de divergencia en todas las hojas, c\'alculo de la evoluci\'on temporal Runge-Kutta.
\item Si alg\'un valor resulta overflow, el proceso se considera num\'ericamente inestable.
\item Actualizaci\'on de la estructura de \'arbol.
\item Estudio de distintos indicadores de error. C\'alculo de la tasa de compresi\'on.
\end{itemize}
\normalsize
Notar que el algoritmo puede resumirse esquem\'aticamente por
\begin{equation}\label{alg}
u^{n+1}=\mathbf{\bar{E}}(\Delta t)\cdot\mathbf{\bar{M}}^{-1}\cdot\mathbf{Tr}(\varepsilon)\cdot\mathbf{\bar{M}}\cdot u^n,
\end{equation}
donde $\mathbf{\bar{M}}$ es el operador de multiresoluci\'on (Codificaci\'on), $\mathbf{Tr}(\varepsilon)$ es el operador de truncamiento con la tolerancia prescrita $\varepsilon$, y $\mathbf{\bar{E}}(\Delta t)$ es el operador discreto de evoluci\'on temporal.
\newpage
\section{Resultados num\'ericos}
En esta secci\'on se reproducen los resultados num\'ericos en 1D obtenido por Roussell \emph{et al.} \cite{Sch} y Bihari \cite{Bihari}, utilizando para la evoluci\'on temporal un m\'etodo expl{\'i}cito Runge-Kutta TVD de segundo orden; para la discretizaci\'on del t\'ermino advectivo se utiliza el esquema cl\'asico de Roe, con los estados de izquierda y derecha obtenidos mediante interpolaci\'on ENO de segundo orden y para la discretizaci\'on de la parte difusiva, se utiliza un esquema centrado de segundo orden. Se utiliza un orden de precisi\'on para el interpolador de multiresoluci\'on de $r=2$. Se utilizan mallas finas de 256, 512, 1024, 2048 y 4096 vol\'umenes de control, tolerancias prescritas de $\varepsilon=5\times10^{-3}$ y $\varepsilon=10^{-3}$, niveles de multiresoluci\'on hasta $L=13$ y una estrategia para el operador de truncamiento $\varepsilon_k=\frac{\varepsilon}{2^{L-k}}$, $1\leqslant k\leqslant L$.
\subsection{Ecuaci\'on de convecci\'on-difusi\'on en 1D}
En el caso de que el flujo sea lineal, se considera la ecuaci\'on lineal de convecci\'on-difusi\'on para $(x,t)\in\ [-1,1]\times[0,\infty[$, $c>0$, $\nu>0$,
\begin{equation}\label{convecc}
\frac{\partial u}{\partial t}+c\frac{\partial u}{\partial x}=\nu\frac{\partial^2u}{\partial x^2}.
\end{equation}
Si se considera como escala espacial caracter{\'i}stica al largo del dominio $X$ y como escala temporal caracter{\'i}stica a $T=c/X$, (\ref{convecc}) puede escribirse en la forma adimensional siguiente
\begin{equation}\label{convecc2}
\frac{\partial u}{\partial t}+\frac{\partial u}{\partial x}=\frac{1}{Pe}\frac{\partial^2u}{\partial x^2},
\end{equation}
donde $Pe$ denota el \emph{n\'umero de Peclet} $Pe=cX/\nu$.
Se estudia (\ref{convecc2}) asociada a la condici\'on inicial

\begin{equation}\label{discont1}
u(x,0)=u_0(x)=\left\{\begin{array}{ll}
1,&\textrm{ si }x\leqslant 0,\\
0,&\textrm{ si }x> 0
\end{array}\right.
\end{equation}

y condiciones de Dirichlet en la frontera
\begin{eqnarray*}
u(-1,t)&=&1,\\
u(1,t)&=&0.
\end{eqnarray*}
La soluci\'on anal{\'i}tica est\'a dada por Hirsch \cite{Bihari}
\begin{equation}
u_{ex}(x,t)=\frac{1}{2}\textrm{erfc}\left(\frac{x-t}{2}\sqrt{\frac{Pe}{t}}\right).
\end{equation}

\begin{figure}[h]
\begin{center}
\includegraphics[width=6in,height=3in]{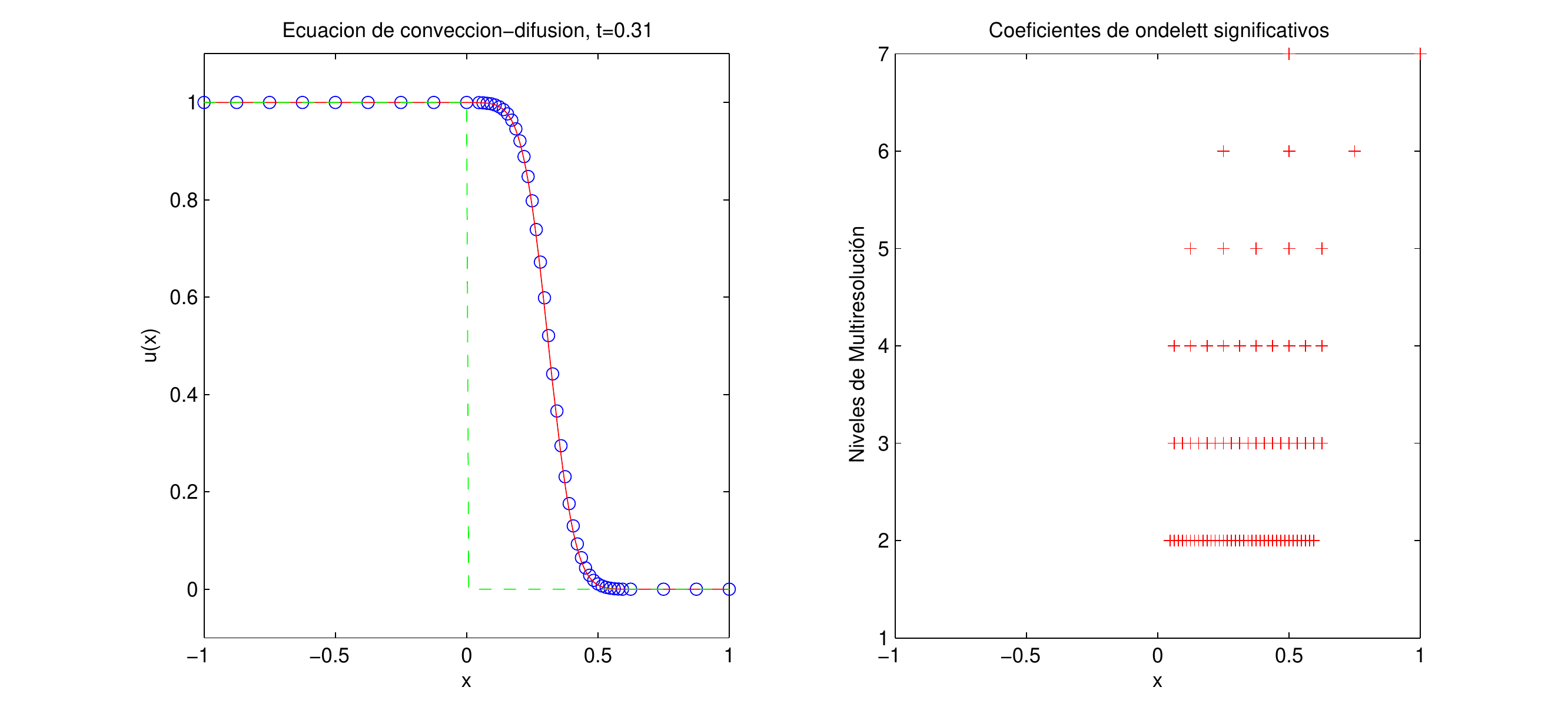}		
\caption[Soluci\'on num\'erica de multiresoluci\'on en el tiempo $t=0.31$ para la ec. de convecci\'on-difusi\'on en 1D]{\small Izquierda: Soluci\'on inicial \emph{(rayas)}, soluci\'on anal{\'i}tica \emph{(linea)}, y soluci\'on num\'erica de multiresoluci\'on \emph{(c{\'i}rculos)} en el tiempo $t=0.31$ para la ec. de convecci\'on-difusi\'on en 1D asociada a la condici\'on inicial (\ref{discont1}), con $Pe=100$, $L=7$, $\varepsilon=10^{-3}$ y $N_0=257$. Derecha: Estructura de coeficientes de ondelette significativos correspondientes.}
\label{fig:cd01}
\end{center}
\end{figure}

\begin{figure}[tp]
\begin{center}
\includegraphics[width=6in,height=3in]{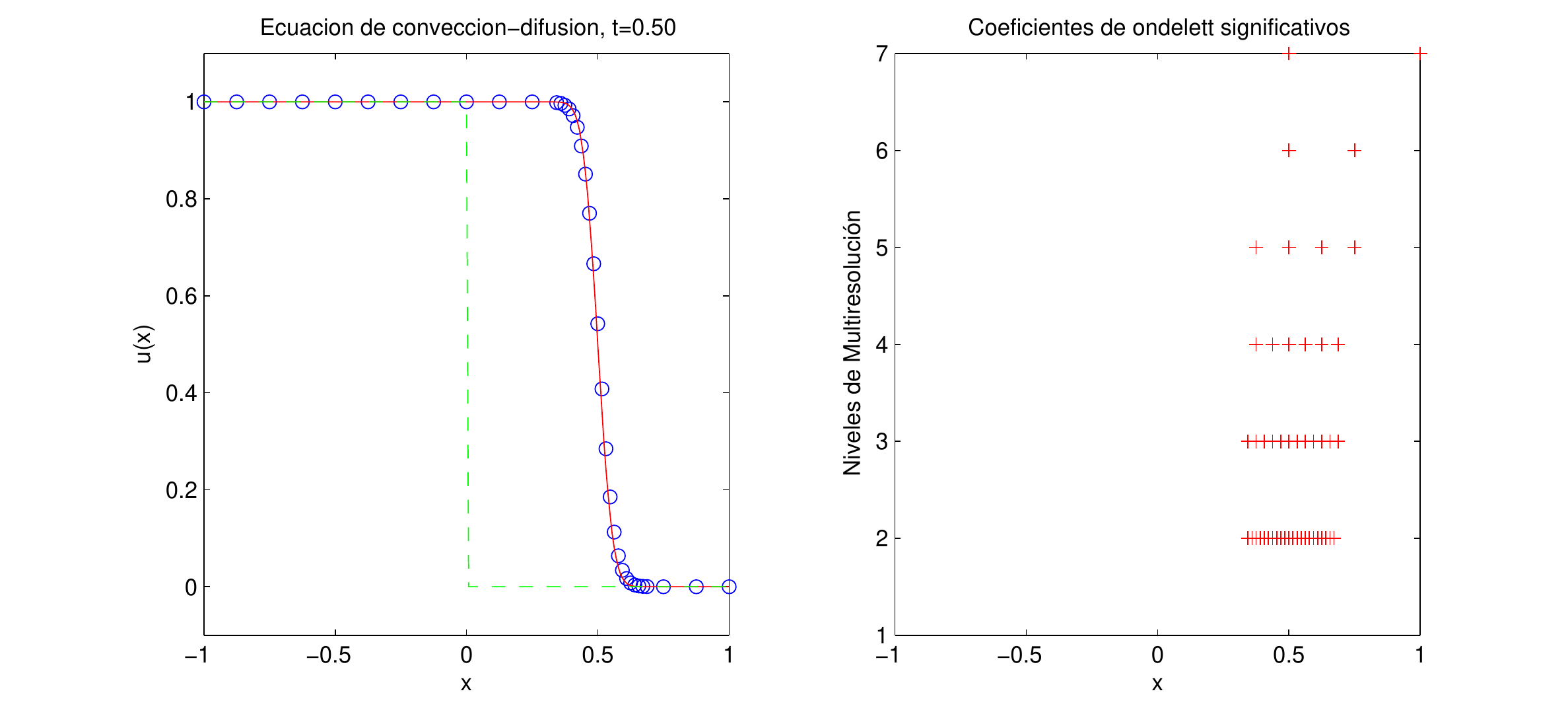}
\caption[Soluci\'on num\'erica de multiresoluci\'on en el tiempo $t=0.50$ para la ec. de convecci\'on-difusi\'on en 1D]{\scriptsize Izquierda: Soluci\'on inicial \emph{(rayas)}, soluci\'on anal{\'i}tica \emph{(linea)}, y soluci\'on num\'erica de multiresoluci\'on \emph{(c{\'i}rculos)} en el tiempo $t=0.50$ para la ec. de convecci\'on-difusi\'on en 1D asociada a la condici\'on inicial (\ref{discont1}), con $Pe=1000$, $L=7$, $\varepsilon=10^{-3}$ y $N_0=257$. Derecha: Estructura de coeficientes de ondelette significativos.}
\label{fig:rouss8}
\end{center}
\end{figure}

\begin{figure}[p]
\begin{center}
\includegraphics[width=6in,height=3in]{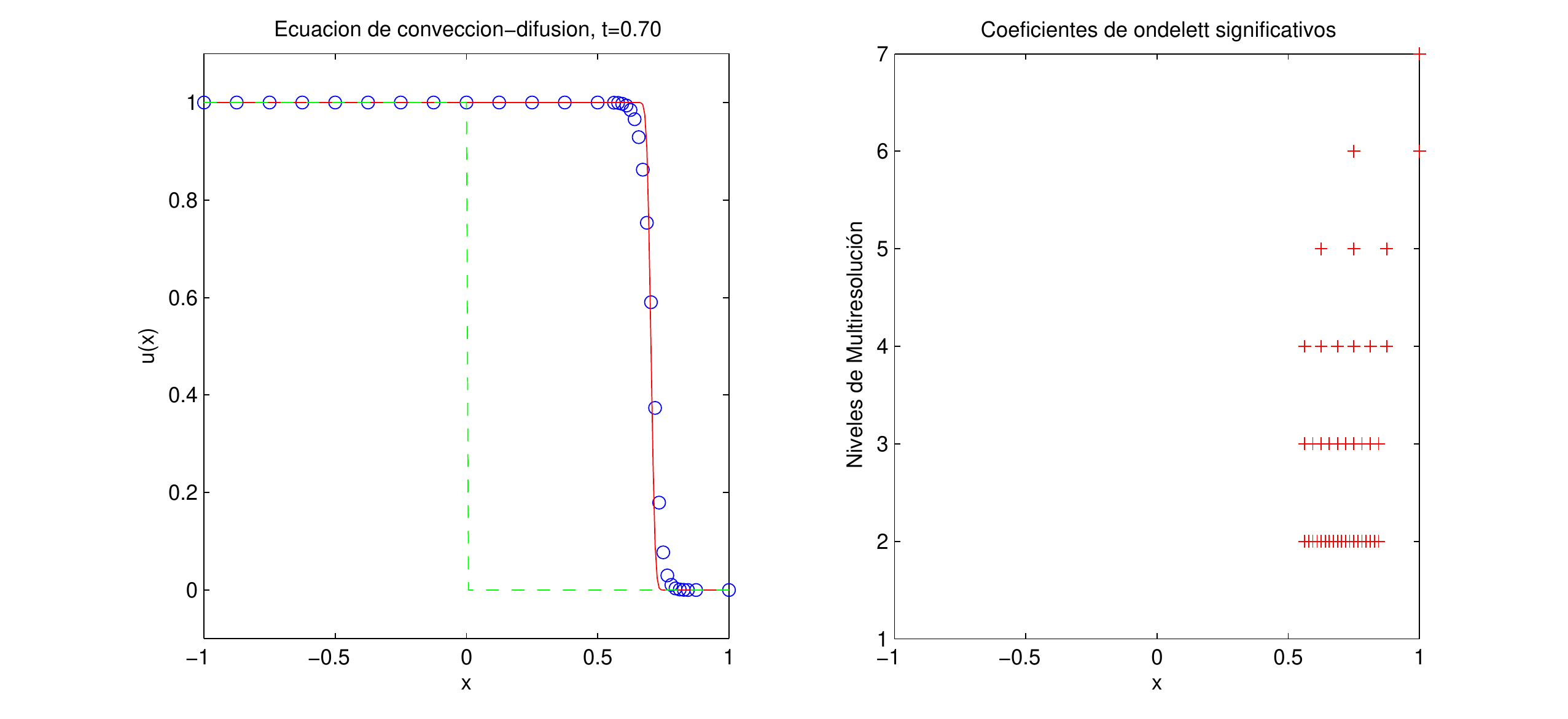}
\caption{\scriptsize Izquierda: Soluci\'on inicial \emph{(rayas)}, soluci\'on anal{\'i}tica \emph{(linea)}, y soluci\'on num\'erica de multiresoluci\'on \emph{(c{\'i}rculos)} en el tiempo $t=0.70$ para la ec. de convecci\'on-difusi\'on en 1D asociada a la condici\'on inicial (\ref{discont1}), con $Pe=10000$, $L=7$, $\varepsilon=10^{-3}$ y $N_0=257$. Derecha: Estructura de coeficientes de ondelette significativos.}
\label{fig:cd03}
\end{center}
\end{figure}

Se testearon tres casos en que el par\'ametro de control es el n\'umero de Peclet $Pe$:
\begin{itemize}
\item[i)] $Pe=100$. En la figura \ref{fig:cd01} (izquierda) se muestra la soluci\'on num\'erica de (\ref{convecc2}) en el tiempo $t=0.3125$. Se observa el fen\'omeno de propagaci\'on lineal de la discontinuidad hacia la derecha. Notar de la tabla \ref{tabla:bihari1}, que los errores al comparar la soluci\'on obtenida mediante multiresoluci\'on y la soluci\'on obtenida sin aplicar el proceso de multiresoluci\'on, son bastante peque\~nos, pero se acumulan con el paso del tiempo.
\item[ii)] $Pe=1000$. En la figura \ref{fig:rouss8} (izquierda) se muestra la soluci\'on num\'erica de (\ref{convecc2}) en el tiempo $t=0.5$. La suavidad de la soluci\'on se debe principalmente a la difusividad.
\item[iii)] $Pe=10000$. Este caso es cercano al caso l{\'i}mite en que la viscosidad es baja en extremo, y el efecto ``suavizante'' es bastante lento. Este caso (y se ver\'a lo mismo para el caso no lineal), es un ejemplo de que la soluci\'on inv{\'i}scida puede obtenerse haciendo $\nu\to0$. En la figura \ref{fig:cd03} (izquierda) se muestra la soluci\'on num\'erica de (\ref{convecc2}) en el tiempo $t=0.7031$. Notar de la
tabla \ref{tabla:bihari1}, la tasa de compresi\'on es considerablemente alta.
\end{itemize}

\begin{table}[tp]
\begin{center}
\begin{tabular}{lccccc}
\hline
$Pe$ & $n$& $\mu$  &            $e_1$   &        $e_2$     &     $e_\infty$     \\
\hline
100  & 10 & 24.0963 & 6.00$\times10^{-4}$&1.88$\times10^{-4}$&9.01$\times10^{-4}$\\
     &100 & 23.9254 & 1.90$\times10^{-4}$&1.02$\times10^{-4}$&6.14$\times10^{-4}$\\
     &200 & 23.6491 & 4.31$\times10^{-4}$&6.58$\times10^{-4}$&8.65$\times10^{-4}$\\
     &600 & 24.1358 & 8.29$\times10^{-4}$&7.17$\times10^{-4}$&9.61$\times10^{-4}$ (*)\\
&&&&&\\
1000 & 10 & 28.2134 & 7.56$\times10^{-4}$&2.02$\times10^{-4}$&7.90$\times10^{-3}$\\
     &100 & 27.6779 & 8.38$\times10^{-6}$&6.80$\times10^{-5}$&6.56$\times10^{-4}$ (*)\\
     &200 & 28.7502 & 1.45$\times10^{-5}$&2.72$\times10^{-4}$&9.77$\times10^{-4}$ \\
     &600 & 28.0683 & 4.04$\times10^{-4}$&5.65$\times10^{-4}$&1.00$\times10^{-3}$ \\
&&&&&\\
10000& 10 & 32.0937 & 1.77$\times10^{-6}$&2.37$\times10^{-5}$&5.82$\times10^{-5}$\\
     &100 & 32.0901 & 1.93$\times10^{-5}$&2.70$\times10^{-4}$&2.22$\times10^{-4}$\\
     &200 & 32.0949 & 1.82$\times10^{-4}$&5.72$\times10^{-4}$&4.43$\times10^{-4}$ (*)\\
     &600 & 32.1005 & 2.94$\times10^{-4}$&7.16$\times10^{-4}$&9.23$\times10^{-4}$\\
\hline
\end{tabular}
\end{center}
\caption{Soluci\'on num\'erica de la Ecuaci\'on de Convecci\'on-difusi\'on en 1D, con condici\'on inicial (\ref{discont1}), $L=7$, $\varepsilon=10^{-3}$ y $N_0=257$. Se adjuntaron figuras para los casos marcados con (*).}
\label{tabla:bihari1}
\end{table}
\newpage
\subsection{Ecuaci\'on de Burgers viscosa en 1D}
Se llevaron a cabo experimentos con la ecuaci\'on de Burgers viscosa, la que contiene un t\'ermino convectivo no lineal, para la cual se conoce soluci\'on anal{\'i}tica.
Para $(x,t)\in\ [-1,1]\times[0,\infty[$, la ecuaci\'on puede ser escrita en su forma adimensional:
\begin{equation}\label{adimens}
\frac{\partial u}{\partial t}+\frac{\partial}{\partial x}\left(\frac{u^2}{2}\right)=\frac{1}{Re}\frac{\partial^2u}{\partial x^2},
\end{equation}
donde $Re=\frac{\Delta x}{\nu}$ es el n\'umero de Reynolds.
\subsubsection{Dato inicial suave}
Asociada a la ecuaci\'on (\ref{adimens}), consid\'erese la condici\'on inicial
\begin{equation}\label{smooth}
u(x,0)=u_0(x)=\sin(\pi x),\quad -1\leqslant x<1
\end{equation}
y condiciones de borde peri\'odicas.
Excepto en el caso l{\'i}mite cuando $Re$ es muy grande, nunca existe un choque completamente discontinuo. Como se muestra en los resultados siguientes, se obtienen tasas de compresi\'on cercanas a 4.

\begin{figure}[pt]
\begin{center}
\includegraphics[height=3in,width=6in]{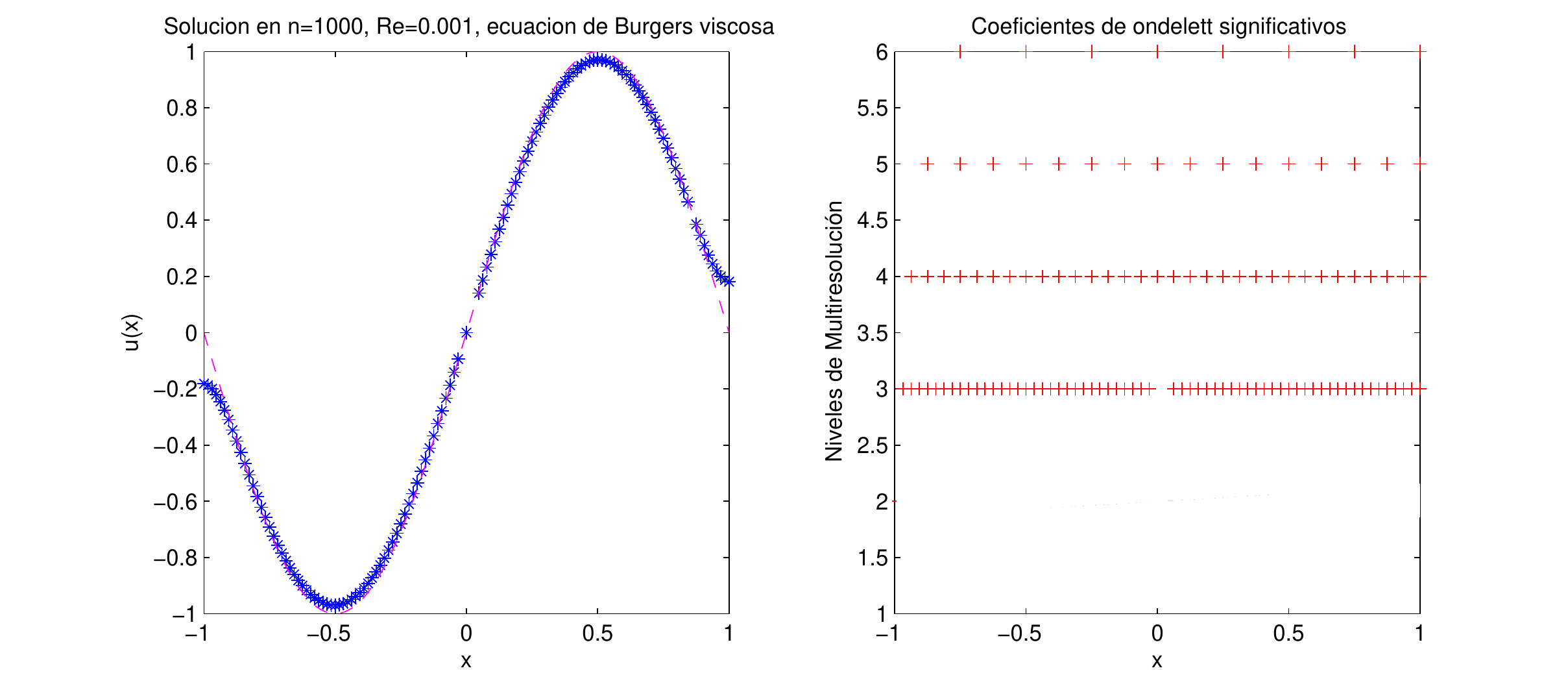}
\caption{\scriptsize Izquierda: Soluci\'on \emph{(rayas)} y soluci\'on num\'erica de multiresoluci\'on \emph{(asteriscos)} en el paso temporal $n=1000$ para la ec. de Burgers viscosa, con $Re=0.001$, $L=7$, $N_0=257$ y $\varepsilon=10^{-3}$. Derecha: Estructura de coeficientes de ondelette significativos correspondientes.}
\label{fig:bihari11}
\end{center}
\end{figure}

\begin{figure}[pb]
	\begin{center}
		\includegraphics[height=3in,width=6in]{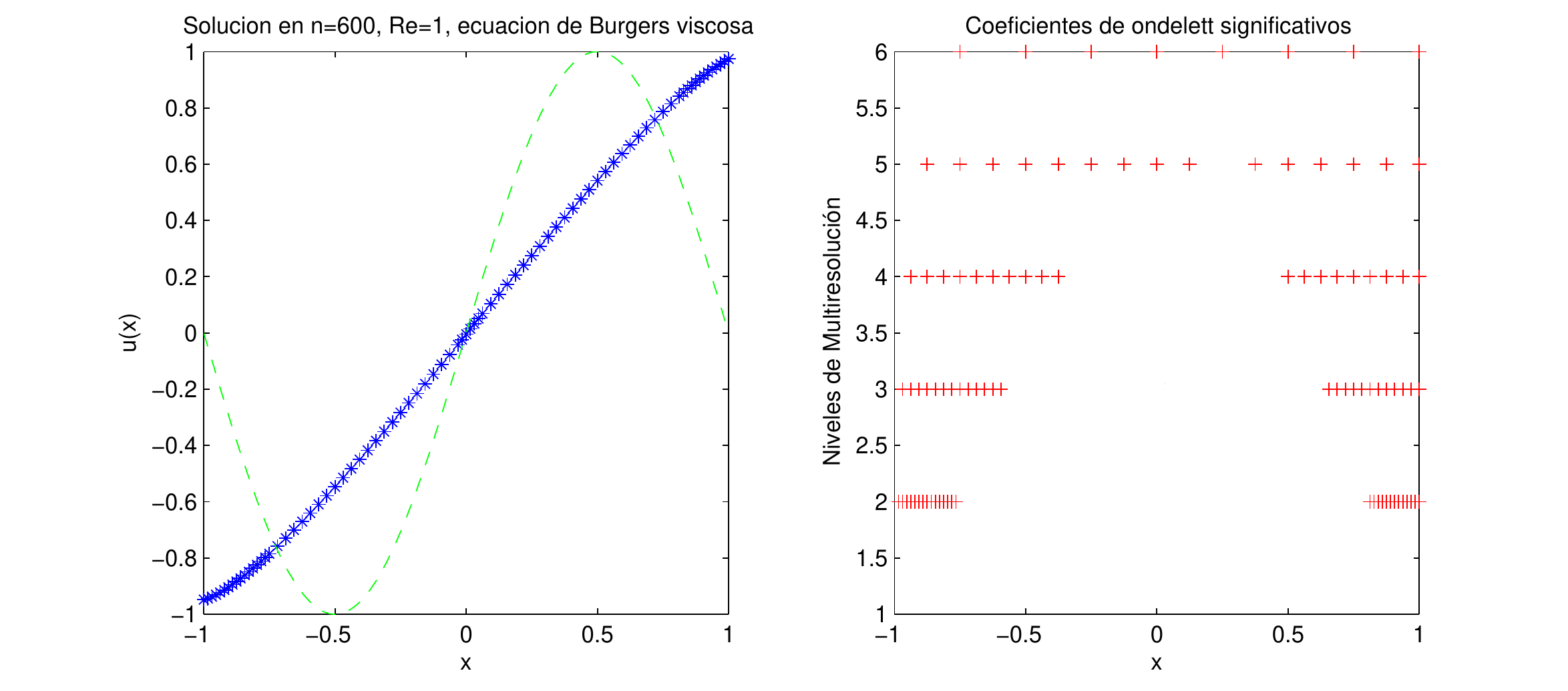}
		\caption{\scriptsize Izquierda: Soluci\'on \emph{(rayas)} y soluci\'on num\'erica de multiresoluci\'on \emph{(asteriscos)} en el paso temporal $n=600$ para la ec. de Burgers viscosa, con $Re=1$, $L=7$, $N_0=257$ y $\varepsilon=10^{-3}$. Derecha: Estructura de coeficientes de ondelette significativos correspondientes.}
		\label{fig:bihari12}
	\end{center}
\end{figure}

\begin{figure}[ht]
	\begin{center}
		\includegraphics[height=3in,width=6in]{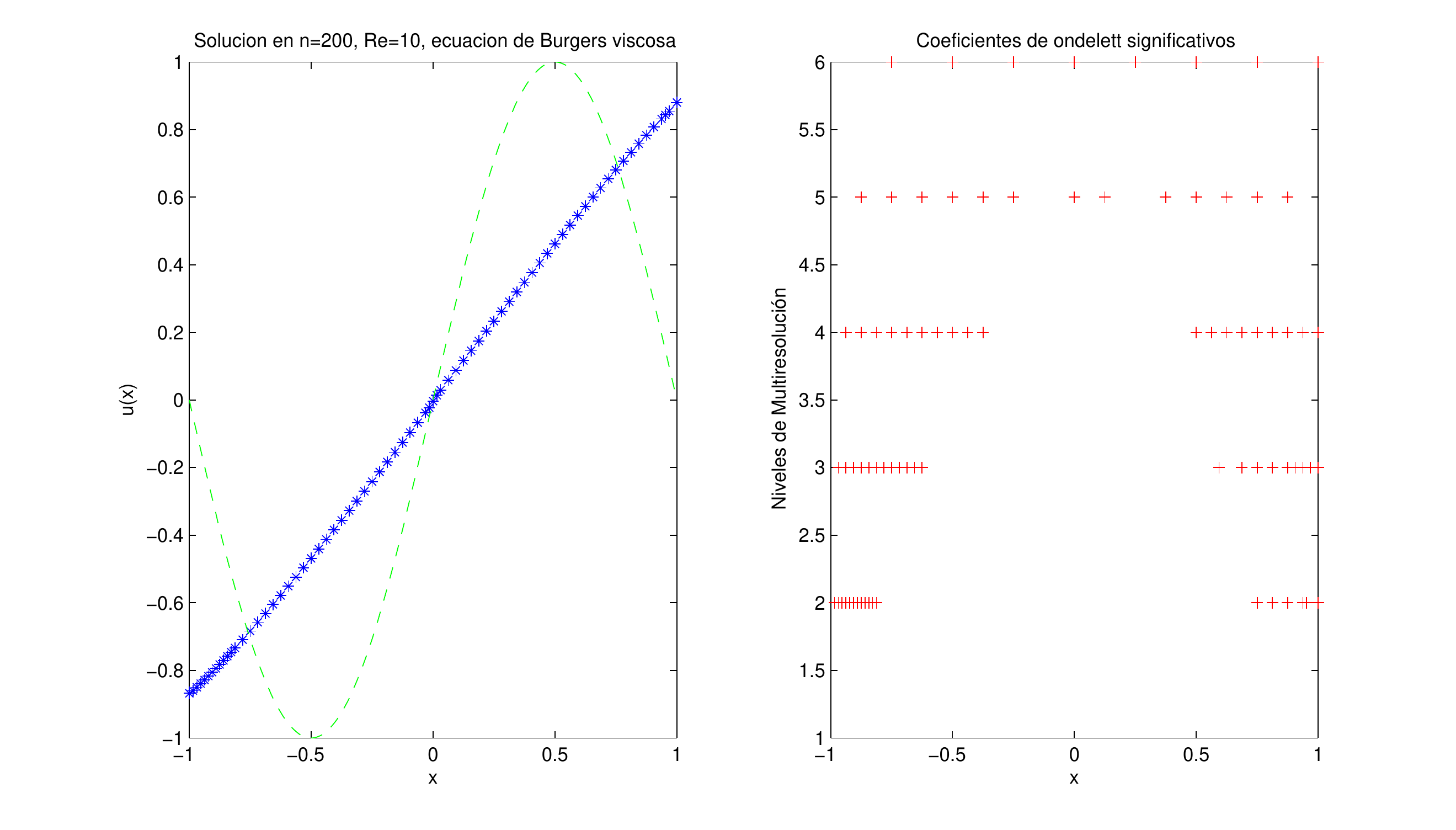}
		\caption{\small Izquierda: Soluci\'on \emph{(rayas)} y soluci\'on num\'erica de multiresoluci\'on \emph{(asteriscos)} en el paso temporal $n=200$ para la ec. de Burgers viscosa, con $Re=10$, $L=7$, $N_0=257$ y $\varepsilon=10^{-3}$. Derecha: Estructura de coeficientes de ondelette significativos correspondientes.}
		\label{fig:bihari13}
	\end{center}
\end{figure}

Se presentan resultados para $Re=0.001$, $Re=1$, y $Re=10$:
\begin{itemize}
\item[ i)] $Re=0.001$. Este caso corresponde a una difusividad grande, lo que provoca que el dato inicial se mantenga suave para todo tiempo $t$. Ver resultados en la tabla \ref{tabla:bihari} y figura \ref{fig:bihari11}.
\item[ ii)] $Re=1$. En $n=600$ se advierte la creaci\'on de una N-onda y el diagrama de coeficientes de multiresoluci\'on es similar al obtenido en el caso inv{\'i}scido. Ver resultados en la tabla \ref{tabla:bihari} y figura \ref{fig:bihari12}.
\item[ iii)] $Re=10$. Este caso produce resultados similares a los obtenidos en el caso inv{\'i}scido. Debido a la capacidad del algoritmo de mantener perfiles afilados, la tasa de compresi\'on se mantiene bastante alta. La figura \ref{fig:bihari13} muestra que el choque se encuentra en un estado casi estacionario. Ver resultados en la tabla \ref{tabla:bihari} y figura \ref{fig:bihari13}.
\end{itemize}

Notar que en todos los casos, los errores son bastante peque\~nos; por lo tanto la calidad de la soluci\'on no se ve comprometida al aplicar el proceso de multiresoluci\'on.

En la secci\'on siguiente se ver\'a que para un n\'umero de Reynolds bastante grande, el problema viscoso no necesita un tratamiento especial, y puede utilizarse el proceso de multiresoluci\'on desarrollado para leyes de conservaci\'on hiperb\'olicas.

\begin{table}[tp]
\begin{center}
\begin{tabular}{lccccc}
\hline
$Re$ & $n$& $\mu$  &            $e_1$   &        $e_2$     &     $e_\infty$     \\
\hline
10   & 10 & 3.0963 & 1.54$\times10^{-4}$&2.24$\times10^{-4}$&5.12$\times10^{-4}$\\
     &100 & 3.8254 & 5.19$\times10^{-4}$&6.69$\times10^{-4}$&1.01$\times10^{-3}$\\
     &200 & 4.6491 & 6.23$\times10^{-4}$&6.72$\times10^{-4}$&7.49$\times10^{-4}$ (*)\\
    &600 & 5.1358 & 7.47$\times10^{-4}$&6.56$\times10^{-4}$&6.05$\times10^{-4}$\\
     &1000& 5.1358 & 8.17$\times10^{-4}$&5.71$\times10^{-4}$&2.91$\times10^{-3}$\\
&&&&&\\
1    & 10 & 4.0198 & 9.86$\times10^{-5}$&1.62$\times10^{-4}$&4.21$\times10^{-3}$\\
     &100 & 3.9876 & 2.31$\times10^{-4}$&2.68$\times10^{-4}$&1.61$\times10^{-4}$\\
     &200 & 3.9902 & 2.92$\times10^{-4}$&2.74$\times10^{-4}$&3.48$\times10^{-5}$\\
     &600 & 4.0299 & 3.47$\times10^{-4}$&2.64$\times10^{-4}$&5.55$\times10^{-5}$ (*)\\
     &1000& 4.3742 & 3.71$\times10^{-4}$&2.48$\times10^{-4}$&1.06$\times10^{-4}$\\
&&&&&\\
0.001& 10 & 4.0279 & 1.15$\times10^{-5}$&1.84$\times10^{-5}$&4.52$\times10^{-5}$\\
     &100 & 4.0198 & 5.71$\times10^{-5}$&6.19$\times10^{-5}$&3.30$\times10^{-4}$\\
     &200 & 4.0198 & 7.42$\times10^{-5}$&1.67$\times10^{-4}$&4.79$\times10^{-4}$\\
     &600 & 4.0021 & 1.24$\times10^{-4}$&2.26$\times10^{-4}$&6.85$\times10^{-4}$\\
     &1000& 4.0021 & 4.07$\times10^{-4}$&4.71$\times10^{-4}$&9.02$\times10^{-4}$ (*)\\
\hline
\end{tabular}
\end{center}
\caption{Soluci\'on num\'erica de la Ecuaci\'on de Burgers viscosa en 1D, condici\'on inicial (\ref{smooth}), $L=7$, $\varepsilon=10^{-3}$ y $N_0=257$. Se adjuntan figuras para los casos marcados con (*).}
\label{tabla:bihari}
\end{table}

\newpage

\subsubsection{Dato inicial discontinuo}
\begin{figure}[ht]
	\begin{center}
		\includegraphics[width=5.6in,height=3in]{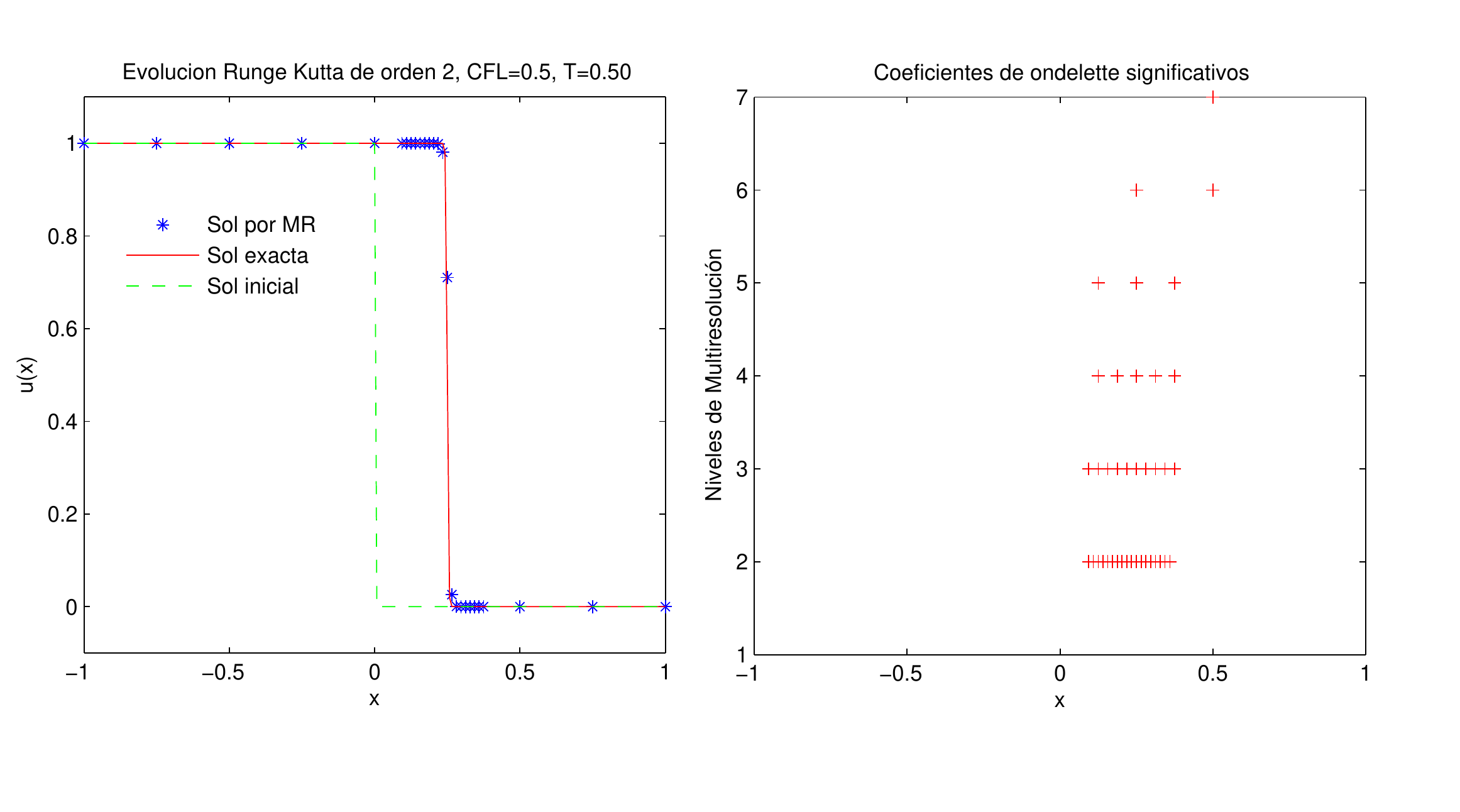}
		\caption{\small Izquierda: Condici\'on inicial \emph{(rayas)}, soluci\'on anal{\'i}tica \emph{(linea)}, y soluci\'on con multiresoluci\'on \emph{(asteriscos)} en el tiempo $t=0.5$, ec. de Burgers viscosa, $Re=1000$, $L=7$, $N_0=257$ y $\varepsilon=10^{-3}$. Derecha: Estructura de coeficientes de ondelette significativos.}
		\label{fig:burg1a}
	\end{center}
\end{figure}

Asociada a la ecuaci\'on (\ref{adimens}), consid\'erese la condici\'on inicial
\begin{equation}\label{discont}
u(x,0)=u_0(x)=\left\{\begin{array}{ll}
1,&\textrm{ si }x\leqslant 0,\\
0,&\textrm{ si }x> 0
\end{array}\right.
\end{equation}
y condiciones de Dirichlet en la frontera
\begin{eqnarray*}
u(-1,t)&=&1,\\
u(1,t)&=&0.
\end{eqnarray*}
La soluci\'on anal{\'i}tica est\'a dada por (\ref{exact3})
\begin{equation}
u_{ex}(x,t)=\frac{1}{2}\left[1-\tanh\left(\left( x-\frac{t}{2}\right)\frac{Re}{4}\right)\right].
\end{equation}
La soluci\'on num\'erica de (\ref{adimens}) en el tiempo $t=0.5$ se muestra en la parte izquierda de la figura \ref{fig:burg1a} para $Re=1000$, $\varepsilon=10^{-3}$ y $L=7$ escalas de multiresoluci\'on, correspondientes a un m\'aximo de 512 vol\'umenes de control en la malla fina. En la parte derecha de la figura \ref{fig:burg1a} se representan los coeficientes de ondelette significativos. Es posible notar el efecto de una propagaci\'on no lineal del choque hacia la derecha, adem\'as puede notarse la difusividad cerca de la discontinuidad.

\begin{figure}[tp]
	\begin{center}
		\includegraphics[width=5.9in,height=3.05in]{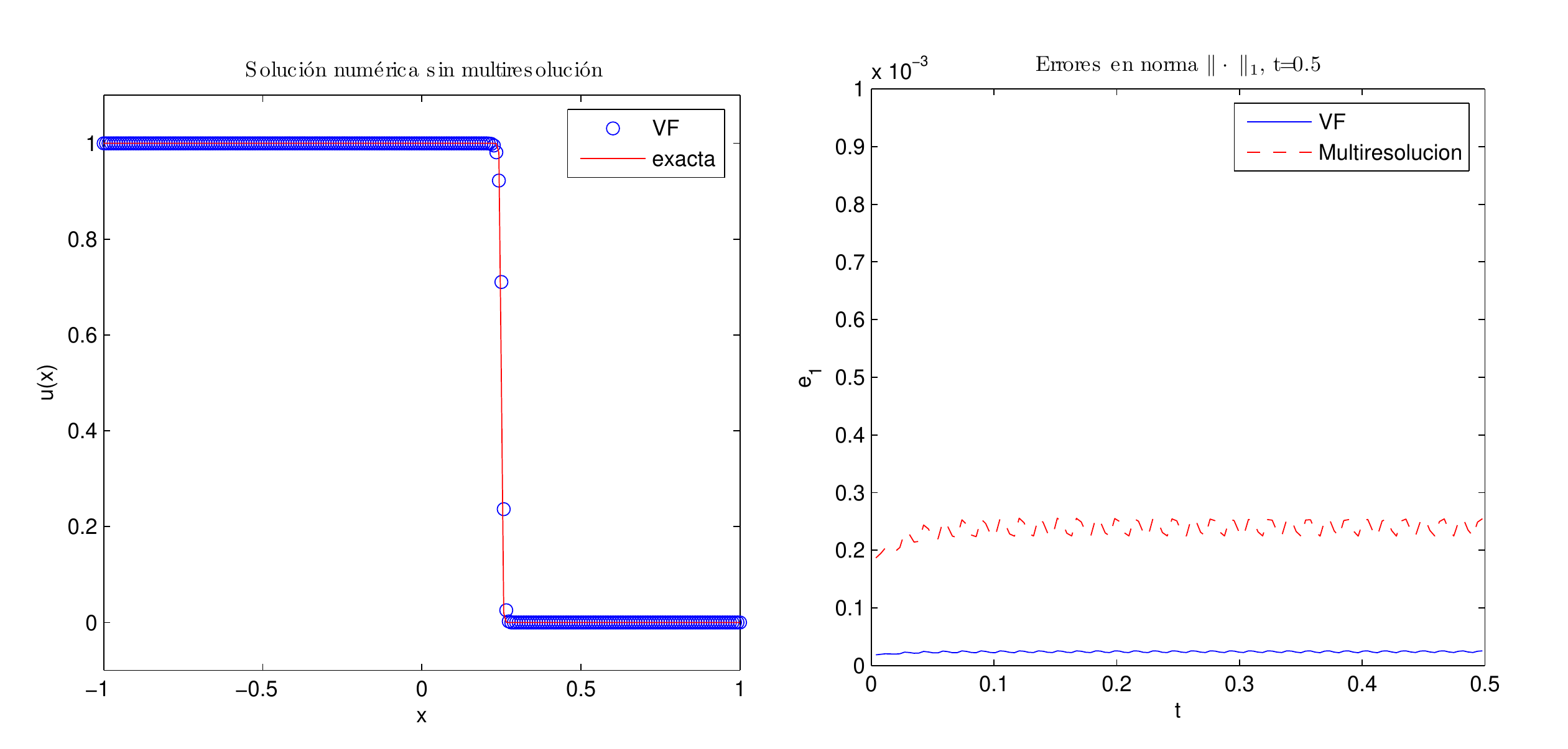}
		\caption{\scriptsize Izquierda: Soluci\'on anal{\'i}tica \emph{(linea)}, y soluci\'on num\'erica sin multiresoluci\'on \emph{(c{\'i}rculos)} en el tiempo $t=0.5$ para la ec. de Burgers viscosa, con $Re=1000$, $L=7$, $N_0=257$ y $\varepsilon=10^{-3}$. Derecha: Errores entre las soluciones anal{\'i}tica y de vol\'umenes finitos con y sin multiresoluci\'on.}
		\label{fig:burg1b}
	\end{center}
\end{figure}

Tambi\'en se presenta la soluci\'on num\'erica obtenida mediante el esquema ENO de segundo orden, con Runge-Kutta de segundo orden (ENO2-RK2) pero sin aplicar multiresoluci\'on (parte izquierda de la figura \ref{fig:burg1b}). La evoluci\'on temporal de los errores entre las soluciones anal{\'i}tica y calculada mediante vol\'umenes finitos con y sin multiresoluci\'on se presenta en la parte derecha de la figura \ref{fig:burg1b}. Notar que los errores est\'an bajo la tolerancia prescrita de $\varepsilon=10^{-3}$.
Como una medida de la mejora en velocidad alcanzada mediante la utilizaci\'on del an\'alisis de multiresoluci\'on, se utiliza la tasa de compresi\'on definida por
\begin{equation}
\mu=\frac{N_0}{N_0/2^L+|D^n|},
\end{equation}
donde $D^n$ es el conjunto de coeficientes de ondelette significativos, en todos los niveles de multiresoluci\'on, en el paso temporal $n$.

En las tablas \ref{table:cpu1} y \ref{table:cpu2} se muestra para diferentes tiempos la constante de proporci\'on $V$ entre el tiempo de CPU total para calcular la soluci\'on num\'erica sin multiresoluci\'on y el tiempo de CPU total para calcular la soluci\'on num\'erica con multiresoluci\'on. N\'otese que de los resultados de las tablas se concluye que la soluci\'on num\'erica tarda alrededor de 1.6 veces el tiempo de CPU que la soluci\'on de multiresoluci\'on.

\vspace{-.5cm}
\begin{figure}[p]
	\begin{center}
		\includegraphics[width=6.2in,height=2.6in]{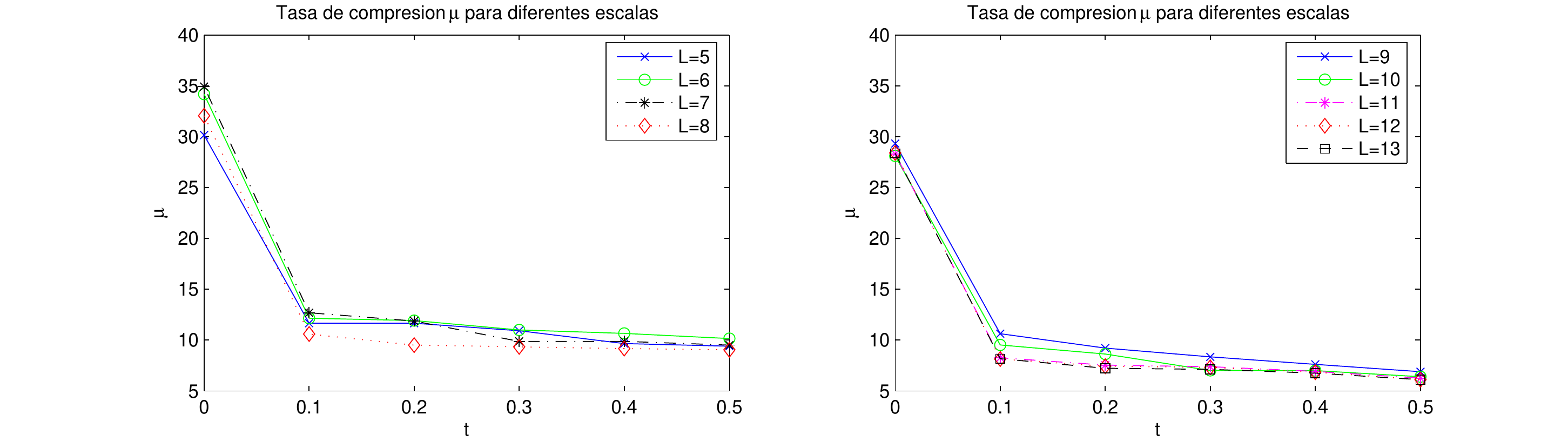}
		\caption{\small Tasa de compresi\'on para distintos niveles m\'aximos de multiresoluci\'on, a distintos tiempos hasta $t=0.5$.}
		\label{fig:tasa}
	\end{center}
\end{figure}

\begin{table}[tp]
{\small
\begin{center}
\begin{tabular}{ccccccccc}
\hline
$L$ & $t$ & $\mu$ & $L$ & $t$ & $\mu$ & $L$ & $t$ & $\mu$\\
\hline
 5 & 0.0019 & 30.1488& 6 & 0.0019 & 34.1822 & 7 & 0.0019 & 34.9112\\
   & 0.1    & 12.6550&   & 0.1    & 12.1503 &   & 0.1    & 12.6866\\
   & 0.2    & 11.9113&   & 0.2    & 11.9131 &   & 0.2    & 11.8701\\
   & 0.3    & 10.1484&   & 0.3    & 10.9858 &   & 0.3    & 9.8646\\
   & 0.4    & 9.9543&    & 0.4    & 10.6509 &   & 0.4    & 9.8646\\
   & 0.5    & 9.3960&    & 0.5    & 10.0013 &   & 0.5    & 9.4993 \\
\hline
 8 & 0.0019 & 32.0586& 9 & 0.0019 & 29.3126 & 10& 0.0019 & 28.1088\\
   & 0.1    & 10.5891&   & 0.1    & 10.6217 &   & 0.1    & 9.5145\\
   & 0.2    & 9.4997&    & 0.2    & 9.2067  &   & 0.2    & 8.6282\\
   & 0.3    & 9.3264&    & 0.3    & 8.3413  &   & 0.3    & 7.0034\\
   & 0.4    & 9.1604&    & 0.4    & 7.5994  &   & 0.4    & 6.9844\\
   & 0.5    & 9.0326&    & 0.5    & 6.8858  &   & 0.5    & 6.3925\\
\hline
 11& 0.0019 & 28.3031& 12& 0.0019 & 28.3012 & 13& 0.0019 & 28.1505\\
   & 0.1    & 8.2575 &   & 0.1    & 8.1348 &   & 0.1    & 8.1003\\
   & 0.2    & 7.5302 &   & 0.2    & 7.4280 &   & 0.2    & 7.2222\\
   & 0.3    & 7.3771 &   & 0.3    & 7.3220 &   & 0.3    & 7.1219\\
   & 0.4    & 6.9207 &   & 0.4    & 6.6343 &   & 0.4    & 6.7472\\
   & 0.5    & 6.2466 &   & 0.5    & 6.1026 &   & 0.5    & 6.1049 \\
\hline
\end{tabular}
\end{center}
\caption{Tasa de compresi\'on para distintos niveles de multiresoluci\'on, hasta $t=0.5$ para la ecuaci\'on de Burgers viscosa en 1D, condici\'on inicial (\ref{discont}). }
}
\end{table}
\newpage
\begin{table}[th]
{\small
\begin{center}
\begin{tabular}{cc}
\hline
t    &  V \\
\hline
0.06 & 2.0011\\
0.12 & 1.9912\\
0.18 & 1.8123\\
0.24 & 1.7780\\
0.36 & 1.6761\\
0.42 & 1.6302\\
0.48 & 1.6079\\
\hline
\end{tabular}
\label{table:cpu1}
\caption{\small Proporci\'on $V$ entre el tiempo de CPU total de la soluci\'on num\'erica ENO2 en malla fina y el tiempo de la soluci\'on de multiresoluci\'on. $N_0=257$, $L=7$ y $\varepsilon=10^{-3}$.}
\end{center}}
\end{table}

Al aumentar el n\'umero de puntos en la malla fina, los resultados obtenidos son a\'un mejores, y en este caso la soluci\'on de multiresoluci\'on tarda menos de la mitad del tiempo total de CPU que tarda la soluci\'on num\'erica que no utiliza multiresoluci\'on.

\begin{table}[h]
{\small
\begin{center}
\begin{tabular}{cc}
\hline
t    &  V \\
\hline
0.06 & 3.0444\\
0.12 & 2.6358\\
0.18 & 2.5129\\
0.24 & 2.5089\\
0.36 & 2.4761\\
0.42 & 2.4341\\
0.48 & 2.4192\\
\hline
\end{tabular}
\label{table:cpu2}
\caption{\small Proporci\'on $V$ entre el tiempo de CPU total de la soluci\'on num\'erica ENO2 en malla fina y el tiempo de la soluci\'on de multiresoluci\'on. $N_0=513$, $L=9$ y $\varepsilon=10^{-3}$.}
\end{center}}
\end{table}

En el caso de sistemas de leyes de conservaci\'on o en el caso de problemas multidimensionales, se espera que $V$ sea a\'un m\'as significativo.
\newpage
\subsection{Ecuaci\'on de reacci\'on-difusi\'on en 1D}
Otro prototipo de una ecuaci\'on parab\'olica no lineal es la ecuaci\'on de \emph{reacci\'on-difusi\'on}. En este caso, la no linealidad no est\'a m\'as en el t\'ermino advectivo (como en la ecuaci\'on de Burgers viscosa) sino en el t\'ermino fuente. Para $(x,t)\in\ [0,20]\times[0,\infty[$, la ecuaci\'on puede ser escrita en su forma adimensional:
\begin{equation}\label{reacc-diff}
\frac{\partial u}{\partial t}=\frac{\partial^2u}{\partial x^2}+S(u),
\end{equation}
con
\begin{equation}
S(u)=\frac{\beta^2}{2}(1-u)\exp\frac{\beta(1-u)}{\alpha(1-u)-1},
\end{equation}
donde $\alpha$ es la tasa de temperatura y $\beta$ es la energ{\'i}a de activaci\'on adimensional (n\'umero de Zeldovich). Se estudia (\ref{reacc-diff}) asociada a la condici\'on inicial
\begin{equation}
u(x,0)=u_0(x)=\left\{\begin{array}{ll}
1,&\textrm{ si }x\leqslant 1,\\
\exp(1-x),&\textrm{ si }x> 1.
\end{array}\right.
\end{equation}
Esta ecuaci\'on conduce al modelo de la \emph{propagaci\'on de una llama premezclada} en 1D, donde las difusividades de masa y calor son iguales. La funci\'on $u$ representa la temperatura adimensional, que var{\'i}a entre 0 y 1. La masa parcial de gas sin quemar es $1-u$. Se elige una condici\'on de Neumann en la frontera izquierda y una condici\'on de Dirichlet en la frontera derecha.
\begin{eqnarray*}
\frac{\partial u}{\partial x}(0,t)&=&0,\\
u(20,t)&=&0.
\end{eqnarray*}
Los par\'ametros son $\alpha=0.8$ y $\beta=10$. EL tiempo final (adimensional) es $t_f=10$. En este ejemplo, la no linealidad del t\'ermino fuente implica que $\Delta t\approx O(\Delta x)$.

La \emph{velocidad de la llama}, definida por
\begin{equation}
v_f=\int_\Omega S\, dx
\end{equation}
se compara con los valores asint\'oticos dados por Peters \& Warnatz \cite{Sch}.

\begin{figure}[ht]
	\begin{center}
		\includegraphics[width=6in,height=3in]{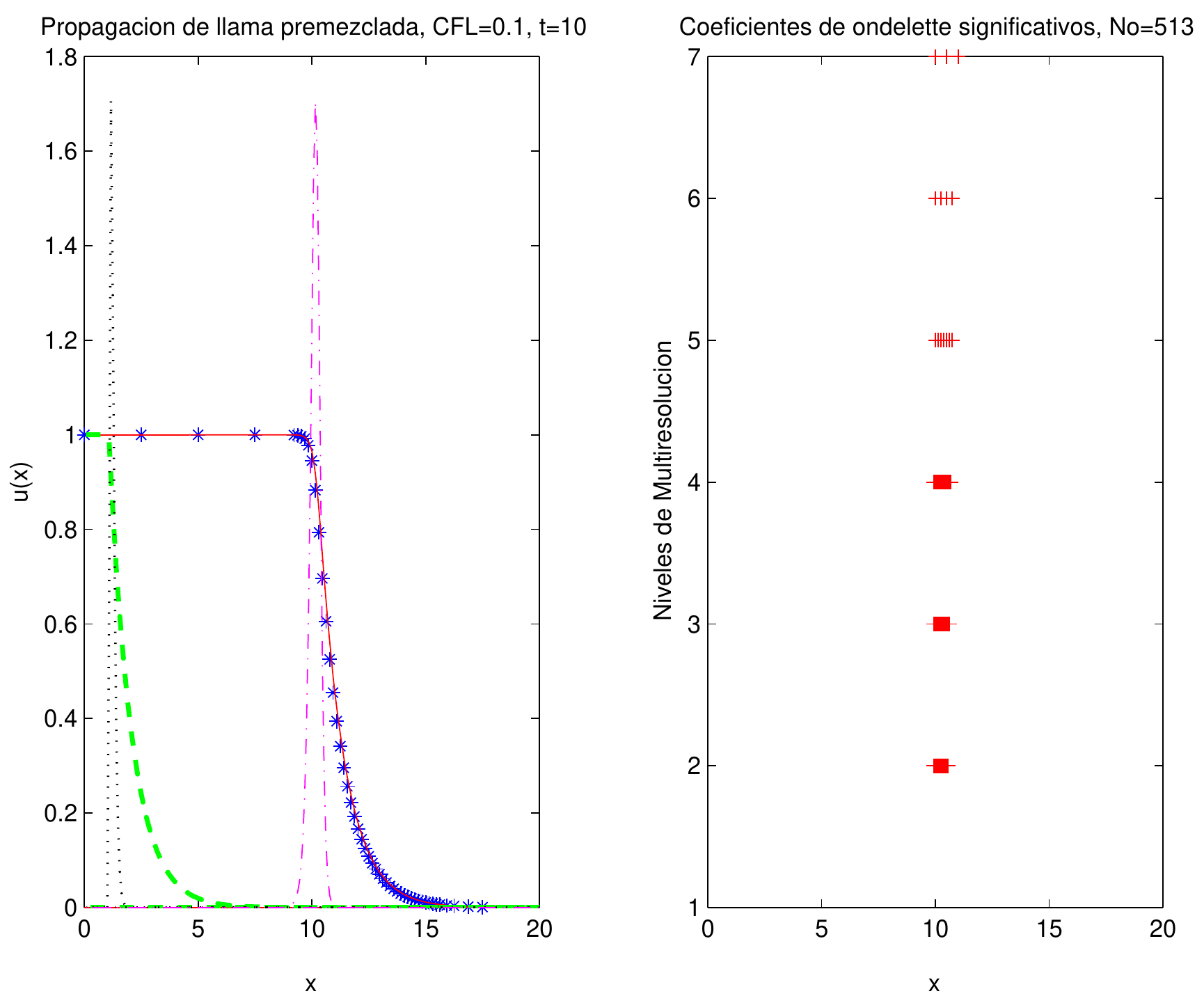}
		\caption{\small Izquierda: Condici\'on inicial \emph{(rayas)} y $S(u)$ inicial \emph{(puntos)}, soluci\'on num\'erica sin multiresoluci\'on \emph{(linea)}, soluci\'on num\'erica con multiresoluci\'on \emph{(asteriscos)} y $S(u)$ \emph{(puntos-rayas)}, en el tiempo $t=10$ para la ec. de reacci\'on-difusi\'on, con $\alpha=0.8$, $\beta=10$, $L=7$, $N_0=513$ y $\varepsilon=10^{-3}$. Derecha: Estructura de coeficientes de ondelette significativos, $t=0.5$.}
		\label{fig:flame}
	\end{center}
\end{figure}

En la figura \ref{fig:flame} se observa la propagaci\'on de la llama en la direcci\'on $x$. El mayor nivel es alcanzado en la regi\'on de la zona de reacci\'on, es decir, para $x\approx 10$.

\begin{table}[hb]
\begin{center}
\begin{tabular}{lcc}
\hline
M\'etodo                      &$v_f$    & $\mu$  \\
\hline
VF                             &0.9146  &        \\
MR $\varepsilon=5\times10^{-2}$&0.9182  & 12.5648\\
MR $\varepsilon=10^{-3}$       &0.9151  & 13.8977 (*)\\
Valor asint\'otico            &0.9080   &        \\
\hline
\end{tabular}
\end{center}
\caption{Velocidad de la llama y tasa de compresi\'on para la soluci\'on num\'erica de (\ref{reacc-diff}) sin multiresoluci\'on (VF), y a dos niveles distintos de tolerancia prescrita para el caso multiresolutivo. $N_0=513$. (*) representado en la figura \ref{fig:flame}.}
\end{table}

\clearpage{\pagestyle{empty}\cleardoublepage}
\chapter{Ecuaci\'on de convecci\'on-difusi\'on fuertemente degenerada}
En este cap{\'i}tulo se presentar\'a un m\'etodo num\'erico para obtener soluciones aproximadas de problemas provenientes de
la sedimentaci\'on de suspensiones floculadas. Estos procesos se utilizan para lograr la separaci\'on de una suspensi\'on de
peque\~nas part{\'i}culas suspendidas en un l{\'i}quido viscoso, en sus componentes s\'olido y l{\'i}quido bajo la acci\'on
de la fuerza de gravedad. Estos procesos se usan ampliamente en la industria minera, por ejemplo para recuperar el agua de
las suspensiones que salen de los procesos de flotaci\'on \cite{Thick}.

La idea principal es aplicar los m\'etodos de multiresoluci\'on a los esquemas desarrollados por B\"urger \emph{et al.}
\cite{BCS,BEK,BEKL,BK,BWC} y observar que el m\'etodo de multiresoluci\'on descrito y ejemplificado en los cap{\'i}tulos
anteriores es de gran ayuda para reducir el costo computacional en este tipo de problemas sin afectar la calidad de la soluci\'on.

Se dar\'a una breve descripci\'on del problema f{\'i}sico y su modelaci\'on mediante una ley de conservaci\'on fuertemente degenerada
con flujo no lineal \cite{BEK}. El efecto de la compresibilidad del sedimento puede ser descrito por un t\'ermino difusivo fuertemente
degenerado, mientras el flujo unidimensional contribuye una discontinuidad de flujo a la ecuaci\'on parcial diferencial. Se presentar\'a
un esquema de segundo orden desarrollado en B\"urger y Karlsen \cite{BK} para resolver este tipo de problemas y finalmente se desarrollan
ejemplos num\'ericos para comparar con los resultados publicados en \cite{BEK,BEKL,BK}.

\begin{figure}[ht]
\begin{center}
\includegraphics[width=4.9in]{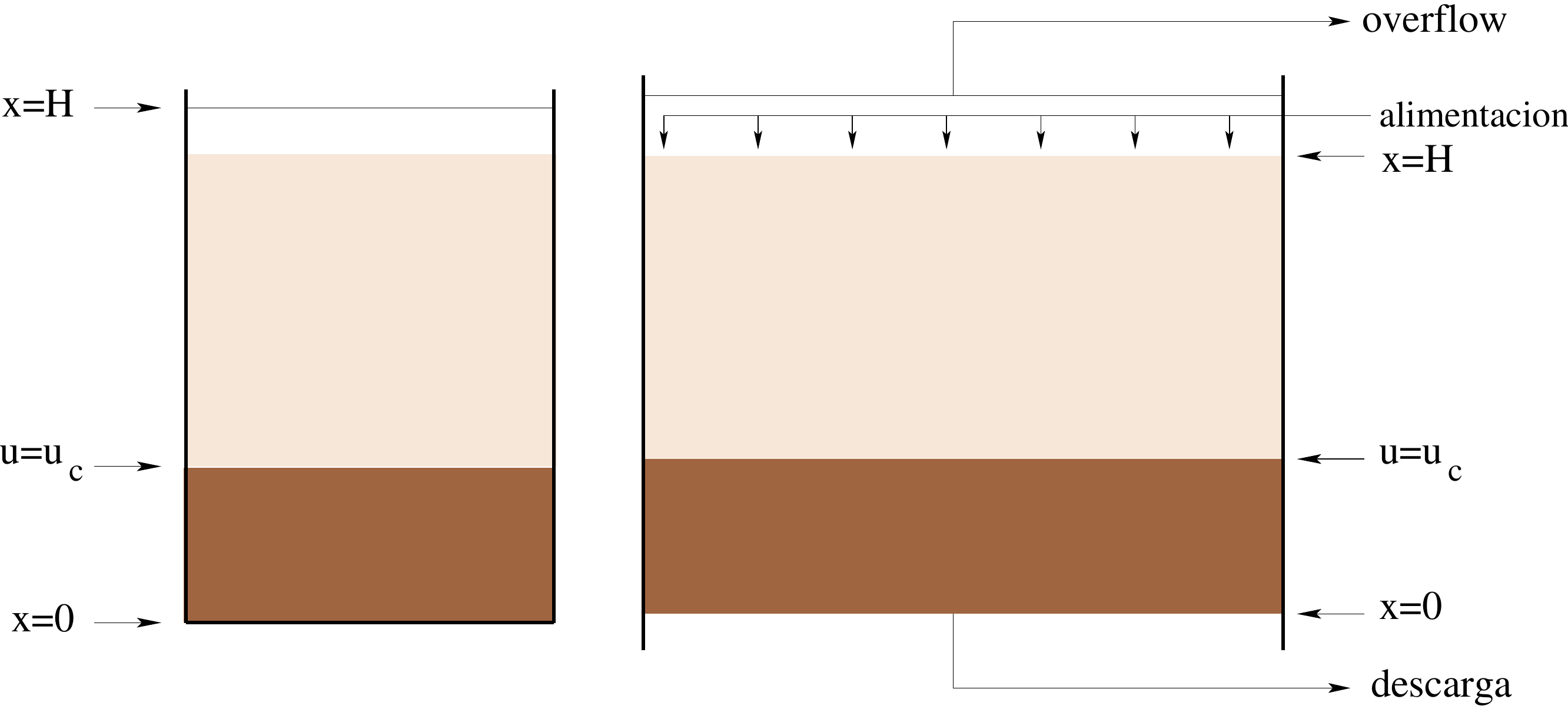}
\caption{\small Izquierda: Columna de sedimentaci\'on Batch. Derecha: ICT \emph{(Ideal Continuous Thickener)} \cite{BEKL}.}
\label{fig:ICT}
\end{center}
\end{figure}

Consid\'erese el caso de una suspensi\'on floculada en un ICT \emph{(Ideal Continuous Thickener)} como el de la figura \ref{fig:ICT}, derecha. Un ICT es un espesador cil{\'i}ndrico sin efectos de pared, en que las variables dependen s\'olo de la altura $x$ y el tiempo $t$. En $x=H$ se tiene una superficie de alimentaci\'on y en $x=0$ se tiene una superficie de descarga, lo que produce una operaci\'on continua del proceso. Esta modelaci\'on es pr\'acticamente obsoleta, pero es de gran utilidad al momento de ejemplificar el comportamiento simplificado de los procesos de sedimentaci\'on. El caso especial de sedimentaci\'on \emph{batch} se muestra en la parte izquierda de la  figura \ref{fig:ICT}. El recipiente es cerrado.

En el caso unidimensional, la teor{\'i}a de la sedimentaci\'on produce ecuaciones de equilibrio de masa y momentum lineal que pueden simplificarse \cite{Thick} hasta obtener una ecuaci\'on parab\'olica fuertemente degenerada de la forma
\begin{equation}\label{degen1}
\partial_tu+\partial_xf(u)=\partial_{xx}^2A(u),
\end{equation}
con $(x,t)\in ]0,1[\times[0,T[$ y el \emph{coeficiente de difusi\'on integrado} dado por
\begin{equation}
A(u):=\int_0^ua(s)ds,\quad a(u)\geqslant 0.
\end{equation}
En general, se permite que el \emph{coeficiente de difusi\'on} $a(u)$ sea cero sobre intervalos de $u$. En tales casos, (\ref{degen1}) es una ecuaci\'on hiperb\'olica. Por esto, (\ref{degen1}) se denomina tambi\'en ecuaci\'on hiperb\'olica-parab\'olica. A\'un cuando este tipo de ecuaciones modelan una gran variedad de fen\'omenos, se enfatizar\'a en las aplicaciones a los procesos de sedimentaci\'on-consolidaci\'on.

Las soluciones de (\ref{degen1}) desarrollan discontinuidades debido a la no linealidad de la funci\'on de densidad de flujo $f(u)$ y a la degeneraci\'on del coeficiente de difusi\'on. Esto lleva a considerar soluciones entr\'opicas para tener un problema bien puesto. A\'un m\'as, cuando (\ref{degen1}) es puramente hiperb\'olica, los valores de la soluci\'on se propagan sobre rectas caracter{\'i}sticas que podr{\'i}an intersectar las fronteras del dominio espacio-tiempo desde el interior, y esto requiere tratar a las condiciones de Dirichlet como condiciones entr\'opicas \cite{BEK}.

Una gran parte de las ecuaciones constitutivas que se proponen para estos procesos, implican que $a(u)$ tiene un \emph{comportamiento degenerado}, es decir, $a(u)=0$ para $u\leqslant u_c$ y $a(u)$ salta en $u_c$ a un valor positivo, donde $u_c$ es una constante llamada \emph{concentraci\'on cr{\'i}tica}. Se enfatiza entonces el hecho de que el coeficiente de difusi\'on $a(u)$ es degenerado, lo que hace evidente la naturaleza hiperb\'olica-parab\'olica de la ecuaci\'on diferencial (\ref{degen1}).

Consid\'erese el problema de valores iniciales y de frontera (PVIF) siguiente
\begin{eqnarray}
\partial_tu+\partial_x(q(t)u+f(u))&=&\partial_{xx}^2A(u),\quad (x,t)\in ]0,H[\times[0,T[,\label{A1}\\
u(x,0)&=&u_0(x),\quad x\in [0,H],\label{A2}\\
u(H,t)&=&0,\quad t\in]0,T]\label{A3}\\
f(u(0,t))-\partial_xA(u(0,t))&=&0,\quad t\in]0,T],\label{A4}
\end{eqnarray}
conocido como el \emph{Problema A}.
Consid\'erese adem\'as el \emph{Problema B}
\begin{eqnarray}
\partial_tu+\partial_x(q(t)u+f(u))&=&\partial_{xx}^2A(u),\quad (x,t)\in ]0,H[\times[0,T[,\label{B1}\\
u(x,0)&=&u_0(x),\quad x\in [0,H],\label{B2}\\
q(t)u(H,t)-\partial_xA(u(H,t))&=&\Psi(t),\quad t\in]0,T]\label{B3}\\
f(u(0,t))-\partial_xA(u(0,t))&=&0,\quad t\in]0,T].\label{B4}
\end{eqnarray}
Para ambos problemas, $f$ se supone continua y diferenciable a trozos, $f\leqslant 0$, $\supp(f)\subset [0,u_{\max}]$, $\|f'\|_\infty\leqslant \infty$, $a(u)\geqslant 0$, $\supp(a)\subset \supp(f)$, $a(u)=0$ para $u\leqslant u_c$, $0<u_c<u_{\max}$, $q(t)\leqslant 0,\ \forall t\in [0,T]$, $TV(q)<\infty,\ TV(q')<\infty$.

En \cite{BEK} se prueba la existencia y unicidad de soluci\'on entr\'opica para cada uno de estos problemas.

En los modelos de sedimentaci\'on-consolidaci\'on de suspensiones floculadas, la coordenada $x$ aumenta verticalmente, $u=u(x,t)$ representa la concentraci\'on volum\'etrica s\'olida local, $q(t)\leqslant0$ es la velocidad media del flujo de la mezcla (puede ser controlada externamente), $f(u)$ es una funci\'on dada que relaciona la velocidad relativa local s\'olido-fluido con la concentraci\'on de s\'olidos local, y
\begin{equation}\label{a(u)}
a(u)=-\frac{f(u)\sigma_e'(u)}{\Delta\varrho\, gu},
\end{equation}
donde $\Delta\varrho>0$ denota la diferencia de densidad de masa s\'olido-fluido, $g$ es la aceleraci\'on de gravedad, y $\sigma_e'(u)\geqslant 0$ es la derivada de la funci\'on de rigidez s\'olida efectiva.

La propiedad de mayor inter\'es, es que generalmente se supone el siguiente comportamiento para $\sigma_e(u)$:
\begin{equation}
\sigma_e(u)\left\{\begin{array}{ll}
=\textrm{cte.},& \textrm{ si } u\leqslant u_c,\\
>0,& \textrm{ si } u>u_c,\end{array}\right.\textrm{ y }\quad\sigma_e'(u):=\frac{d\sigma_e}{du}\left\{\begin{array}{ll}
=0,& \textrm{ si } u\leqslant u_c,\\
>0,& \textrm{ si } u>u_c.\end{array}\right.
\end{equation}
Notar que la naturaleza degenerada de la ecuaci\'on diferencial (\ref{degen1}) es heredada de esta propiedad.

Las propiedades materiales espec{\'i}ficas de la suspensi\'on son descritas por $f(u)$ y $\sigma_e(u)$. Ejemplos t{\'i}picos para estas funciones modelo son la funci\'on de densidad de flujo del tipo Michaels and Bolger \cite{BEKL}
\begin{equation}\label{laf}
f(u)=v_\infty u\left(1-\frac{u}{u_{\max}}\right)^C,\quad v_\infty<0,\ C>1
\end{equation}
y la funci\'on de rigidez s\'olida efectiva \emph{ley de potencia}
\begin{equation}\label{lasigma}
\sigma_e(u)=\left\{\begin{array}{ll}
0,& \textrm{ si } u\leqslant u_c,\\
\sigma_0\left(\left(\frac{u}{u_c}\right)^n-1\right),& \textrm{ si } u>u_c,\end{array}\right.\quad \sigma_0>0,\ n>1.
\end{equation}

 Las condiciones (\ref{A2}) corresponden a una distribuci\'on inicial de concentraci\'on dada, la condici\'on (\ref{A3}) corresponde a prescribir el valor de la concentraci\'on en $x=L$, las condiciones (\ref{A4}) y (\ref{B4}) equivalen a reducir la densidad de flujo del volumen s\'olido en el fondo del recipiente a su parte convectiva $q(t)u(0,t)$ y la condici\'on (\ref{B3}) corresponde a una condici\'on de flujo en $x=L$.

\section{Esquemas de segundo orden}
Para el esquema expl{\'i}cito a desarrollar, se utilizar\'a una discretizaci\'on similar a la utilizada en la secci\'on \ref{section:fn} (ver detalles en \cite{BEK}). Los t\'erminos advectivo y difusivo son aproximados de diferente forma, con el fin de obtener una discretizaci\'on que mantenga la conservatividad en ambos t\'erminos. Para la parte advectiva puede utilizarse el esquema de Roe cl\'asico con una interpolaci\'on ENO de segundo orden, ya utilizado en los cap{\'i}tulos anteriores, o bien puede utilizarse un esquema de Engquist-Osher \cite{EO} modificado para ser de segundo orden \cite{BEKL,BK,EK}. Para la parte difusiva, se necesita un esquema centrado de segundo orden que mantenga la conservatividad \cite{BEKL}.

Dado que el principal inter\'es se encuentra en la discretizaci\'on del t\'ermino difusivo, consid\'erese la siguiente ecuaci\'on puramente difusiva:
\begin{eqnarray}
\partial_tu&=&\partial_{xx}^2A(u),\\
A(u)&=&\int_0^ua(s)ds.
\end{eqnarray}
Una formulaci\'on conservativa de diferencias finitas para esta ecuaci\'on es
\begin{equation}
\frac{u_j^{n+1}-u_j^n}{\Delta t}=\frac{A(u_{j-1}^n)-2A(u_j^n)+A(u_{j+1}^n)}{(\Delta x)^2}.
\end{equation}
Este esquema es estable y convergente bajo la condici\'on $CFL$ (ver \cite{BK})
\begin{equation}
2\max_u|a(u)|\frac{\Delta t}{(\Delta x)^2}\leqslant 1.
\end{equation}
Adem\'as, debe recordarse que el esquema expl{\'i}cito utilizado para la ecuaci\'on puramente hiperb\'olica es estable bajo la condici\'on $CFL$ (ver \cite{GR})
\begin{equation}
\max_u|f'(u)|\frac{\Delta t}{\Delta x}\leqslant 1.
\end{equation}
El esquema interior resultante para la ecuaci\'on (\ref{degen1}) (si se utiliza  el esquema de Roe cl\'asico (\ref{roe}) con una interpolaci\'on ENO de segundo orden para la parte advectiva) es:
\begin{equation}\label{esqbek}
\frac{u_j^{n+1}-u_j^n}{\Delta t}+q(n\Delta t)\frac{u^-_{j+1}-u^+_{j-1}}{\Delta x}+\frac{F_{j+\frac{1}{2}}-F_{j-\frac{1}{2}}}{\Delta x}=\frac{A(u_{j-1}^n)-2A(u_j^n)+A(u_{j+1}^n)}{(\Delta x)^2}
\end{equation}
con
\begin{equation}
F_{j+\frac{1}{2}}=f^R\left(u^-_{j+1},u^+_{j+1}\right).
\end{equation}

La evoluci\'on temporal se har\'a mediante el m\'etodo Runge-Kutta de segundo orden utilizado en los cap{\'i}tulos anteriores.

Las condiciones de borde (\ref{A4}) y (\ref{B4}) prescritas en $x=0$ se discretizan utilizando (\ref{esqbek}) haciendo:
\begin{equation}
f(u(0,t^n))-\partial_xA(u(0,t^n))\approx F^n_{-\frac{1}{2}}-\frac{A(u_0^n)-A(u_{-1}^n)}{\Delta x}=0,
\end{equation}
de donde se obtiene la expresi\'on para la actualizaci\'on del flujo en $u_0^n$
\begin{equation}
\frac{u_0^{n+1}-u_0^n}{\Delta t}+q(n\Delta t)\frac{u_1^n-u_0^n}{\Delta x}+ \frac{F^n_{\frac{1}{2}}}{\Delta x}=\frac{A(u_1^n)-A(u_0^n)}{(\Delta x)^2}.
\end{equation}
Esta formulaci\'on evita utilizar un valor artificial $u_{-1}^n$.

Para el problema A, la condici\'on de borde en $x=H$ se aproxima simplemente poniendo $u_{N_0}^n=0$, en cambio para el problema B, (\ref{B3}) se aproxima haciendo
\begin{equation}
q(n\Delta t)u^n_{N_0}+F^n_{N_0+\frac{1}{2}}-\frac{A(u^n_{N_0+1})-A(u^n_{N_0})}{\Delta x}=\Psi(n\Delta t).
\end{equation}
Con esto, se obtiene la expresi\'on para la actualizaci\'on del flujo en $u_{N_0}^n$
\begin{equation}
\frac{u_{N_0}^{n+1}-u_{N_0}^n}{\Delta t}+\frac{\Psi(n\Delta t)-q(n\Delta t)u_{N_0}^n}{\Delta x}-\frac{F^n_{N_0-\frac{1}{2}}}{\Delta x}=\frac{A(u_{N_0-1}^n)-A(u_{N_0}^n)}{(\Delta x)^2}.
\end{equation}

Como alternativa a la discretizaci\'on de la parte advectiva, puede utilizarse un esquema de Engquist-Osher modificado mediante extrapolaci\'on de variables MUSCL (\emph{Monotonic Upwind Scheme for Conservation Laws}) para lograr un esquema de segundo orden \cite{BEKL,BK,EK,GR,FV}. Para ello se introduce una funci\'on $u^n(x)$ lineal a trozos definida por
$$u^n(x)=u_j^n+s_j^n(x-x_j),\quad x\in ]x_{j-1/2},x_{j+1/2}[,$$
donde $s_j^n$ es una pendiente adecuada, construida a partir de $u^n$. En las regiones donde $s_j^n=1$, la reconstrucci\'on es lineal y el error de truncamiento es $O((\Delta x)^2)$. En las regiones donde $s_j^n=0$, la reconstrucci\'on es constante a trozos y el error de truncamiento es $O(\Delta x)$. Es necesario utilizar limitadores de pendiente para forzar la monoton{\'i}a de la reconstrucci\'on. En este caso, se utilizar\'a el \emph{$\theta-$limitador} (ver \cite{GR,GR2})
$$ s_j^n=MM\left(\theta\frac{u_j^n-u_{j-1}^n}{\Delta x},\, \frac{u_{j+1}^n-u_{j-1}^n}{2\Delta x},\, \theta\frac{u_{j+1}^n-u_j^n}{\Delta x}\right),\quad \theta\in[0,2],$$
donde $MM$ es otra funci\'on tipo \emph{Min-Mod} definida por
\begin{equation}
MM(a,b,c):=\left\{\begin{array}{ll}
\min(a,b,c),&\textrm{ si } a,b,c>0,\\
\max(a,b,c),&\textrm{ si } a,b,c<0,\\
0,&\textrm{e.o.c.}\end{array}\right.
\end{equation}
Luego se extrapola la informaci\'on hacia la frontera de cada volumen de control, con lo que
\begin{equation}
u_j^L:=u_j^n-\frac{\Delta x}{2}s_j^n,\qquad u_j^R:=u_j^n+\frac{\Delta x}{2}s_j^n.
\end{equation}
As{\'i}, el esquema upwind interior de segundo orden correspondiente se escribe
\scriptsize
\begin{equation}\label{esqbek2}
\frac{u_j^{n+1}-u_j^n}{\Delta t}+q(n\Delta t)\frac{u_{j+1}^L-u_j^R}{\Delta x}+\frac{f^{EO}(u_j^R,u_{j+1}^L)-f^{EO}(u_{j-1}^R,u_j^L)}{\Delta x}=\frac{A(u_{j-1}^n)-2A(u_j^n)+A(u_{j+1}^n)}{(\Delta x)^2},
\end{equation}
\normalsize
donde $f^{EO}(u_j^n,u_{j+1}^n):=f^+(u_j^n)+f^-(u^n_{j+1})$ es el flujo num\'erico de Engquist-Osher \cite{EO},
\begin{equation}
f^+(u)=f(0)+\int_0^u\max(f'(s),0)\, ds,\qquad f^-(u)=\int_0^u\min(f'(s),0)\, ds.
\end{equation}
Este esquema es estable bajo la condici\'on $CFL$ (ver \cite{EK})
\begin{equation}
\max_u|f'(u)|\frac{\Delta t}{\Delta x}+2\max_u|a(u)|\frac{\Delta t}{(\Delta x)^2}\leqslant 1.
\end{equation}
Las condiciones de borde (\ref{A4}) y (\ref{B4}) prescritas en $x=0$ quedan entonces
\begin{equation}
\frac{u_0^{n+1}-u_0^n}{\Delta t}+ q(n\Delta t)\frac{u_1^n-u_0^n}{\Delta x}+ \frac{f^{EO}(u_0^n,u_1^n)}{\Delta x}=\frac{A(u_1^n)-A(u_0^n)}{(\Delta x)^2}
\end{equation}
y la condici\'on de borde (\ref{B3}) queda
\begin{equation}
\frac{u_{N_0}^{n+1}-u_{N_0}^n}{\Delta t}+\frac{\Psi(n\Delta t)-q(n\Delta t)u_{N_0}^n}{\Delta x}-\frac{f^{EO}(u_{N_0-1}^n,u_{N_0}^n)}{\Delta x}=\frac{A(u_{N_0-1}^n)-A(u_{N_0}^n)}{(\Delta x)^2}.
\end{equation}
\newpage
\section{Un algoritmo de multiresoluci\'on}
Se presenta a continuaci\'on una breve descripci\'on de un
algoritmo de multiresoluci\'on para resolver num\'ericamente una
ecuaci\'on parab\'olica fuertemente degenerada.
\begin{enumerate}
\item \emph{Inicializaci\'on de par\'ametros y variables}:
\begin{itemize}
\item Longitud del dominio $H$,
\item concentraci\'on cr{\'i}tica $u_c$,
\item orden de la interpolaci\'on de multiresoluci\'on $r$,
\item niveles de multiresoluci\'on $L$,
\item n\'umero de puntos y paso en la malla fina $N_0$ y $h_0$, y en cada nivel $N_k$ y $h_k$,
\item tolerancia prescrita $\varepsilon$ y estrategia
de truncamiento $\varepsilon_k$,
\item tiempo de simulaci\'on $t_f$,
\item constantes de Lipschitz para $a(u)$ y $f'(u)$,
\item condici\'on $CFL$:
$$\max_u|f'(u)|\frac{\Delta t}{h_0}+2\max_u|a(u)|\frac{\Delta t}{h_0^2}\leqslant 1.$$
\item paso temporal $\Delta t$,
$$\Delta t=\frac{CFL\cdot h_0}{\max_u|f'(u)|+2\max_u|a(u)|/h_0}.$$
\item condiciones iniciales $u_0$ y
\item otros par\'ametros del modelo (\ref{parameters}): $v_\infty$, $C$, $n$,
$u_{\max}$, $\Delta \varrho$, etc.
\item Inicializaci\'on de la estructura de datos. (En este caso, estructura esparsa).
\end{itemize}
\item \emph{Aplicaci\'on de la codificaci\'on a la condici\'on inicial}:
Este proceso entrega los coeficientes de ondelette significativos
y los valores de la soluci\'on en las posiciones correspondientes
a coeficientes de ondelette significativos. Se incluyen los \emph{safety points}.
\item \emph{Evoluci\'on temporal}: Se utiliza un m\'etodo
Runge-Kutta de segundo orden.
\begin{itemize}
\item Primer paso intermedio Runge-Kutta,
\item Segundo paso Runge-Kutta,
\item Actualizaci\'on de los flujos y actualizaci\'on de la soluci\'on,
\item Imposici\'on de condiciones de contorno fijas y condiciones de flujo,
\item Se aplica el paso 2. a la soluci\'on actual y se itera hasta alcanzar el tiempo final.
\end{itemize}
\item \emph{Salidas}: Se realizan gr\'aficos de la soluci\'on
num\'erica y coeficientes de ondelette significativos
correspondientes. Se calculan adem\'as tasas de compresi\'on y
tiempos de CPU para comparar con la resoluci\'on obtenida sin utilizar
multiresoluci\'on.
\end{enumerate}

\newpage
\section{Ejemplos num\'ericos}
Se calculan soluciones de los problemas A y B utilizando los
esquemas num\'ericos descritos en la secci\'on anterior, con
una discretizaci\'on para el flujo de tipo Enqguist-Osher, dada por (\ref{esqbek2}). Se
reproducen algunos resultados num\'ericos obtenidos por B\"urger
\emph{et al.} \cite{BEK,BEKL,BK} y Bustos \emph{et al.}
\cite{Thick}.
\subsection{Sedimentaci\'on batch de suspensi\'on ideal}
Considerar en primer lugar, el proceso de sedimentaci\'on batch de suspensi\'on ideal en una columna de asentamiento \cite{CB}. El caso ideal permite formular el proceso como
\begin{eqnarray*}
\frac{\partial u}{\partial t}+\frac{\partial f(u)}{\partial x}&=&0, \quad x\in \RR,\ t>0,\\
u(x,0)&=&u_0(x),\quad x\in[0,H[,\\
u(0,t)&=&u_{\infty},\quad t>0,\\
u(L,t)&=&u_L,\quad t>0.
\end{eqnarray*}
En el ejemplo se considera una columna de asentamiento de longitud $H=1$, una concentraci\'on inicial $u_0(x)=0.25$, condiciones de borde $u_\infty=0.642$ y $u_0=0$. Se elige una ecuaci\'on constitutiva para la funci\'on de densidad de flujo s\'olido. Se utiliza la funci\'on descrita por Shannon (1963, consultar \cite{Thick})
\small\begin{equation}
f(u)=\left(-0.33843u+1.37672u^2-1.62275u^3-0.11264u^4+0.902253u^5\right)\times10^{-2}\, [m/s].
\end{equation}
\normalsize

\begin{figure}[h]
\begin{center}
\includegraphics[width=3.5in,height=2in]{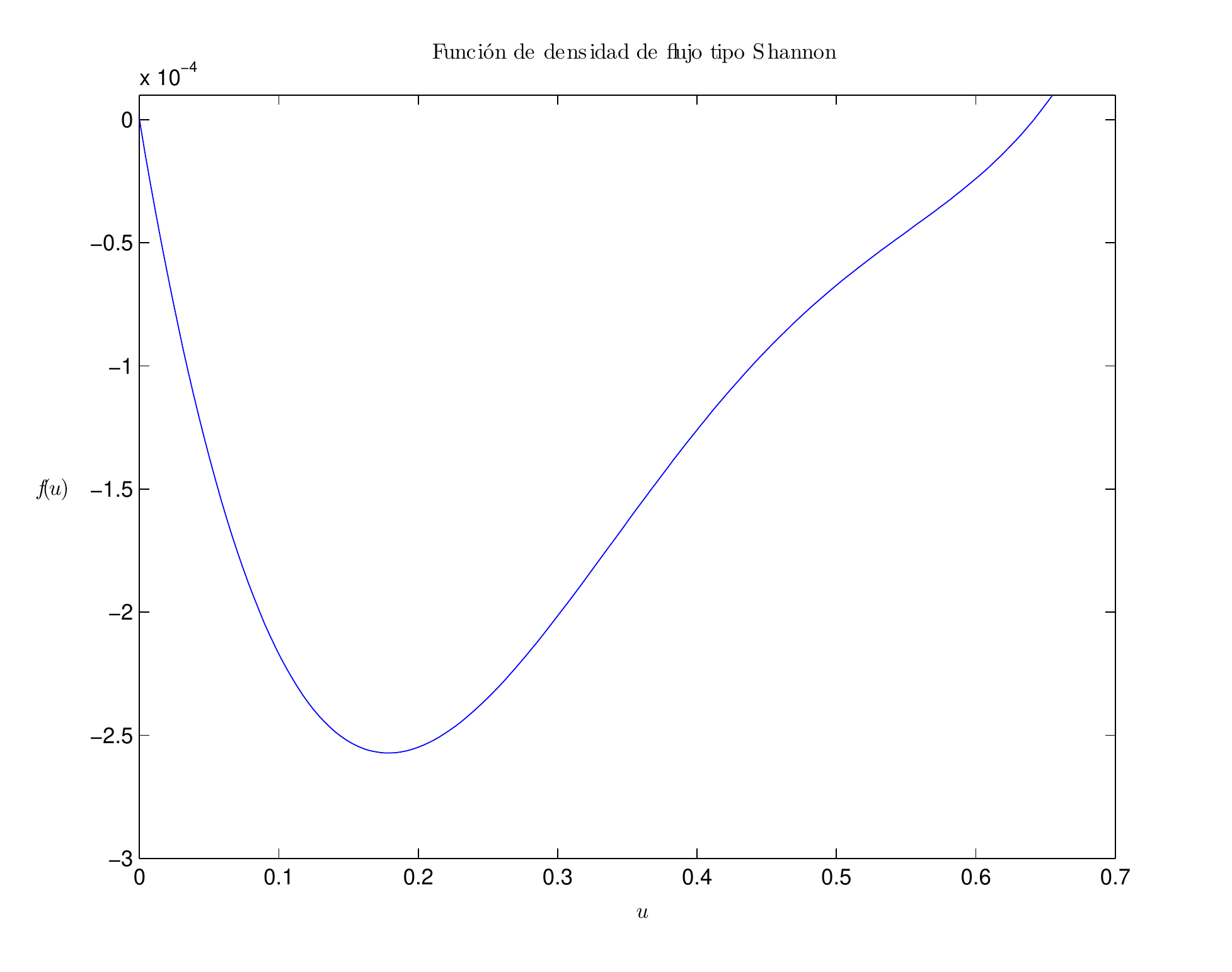}
\caption[Funci\'on modelo para el \emph{problema de sedimentaci\'on batch de suspensi\'on ideal}.]{\small Funci\'on de densidad de flujo $f(u)$ para el problema de sedimentaci\'on batch de suspensi\'on ideal. Unidad: $[m/s]$.}
\label{fig:density_flux}
\end{center}
\end{figure}

\begin{figure}[tp]
\begin{center}
\includegraphics[width=5in,height=3.3in]{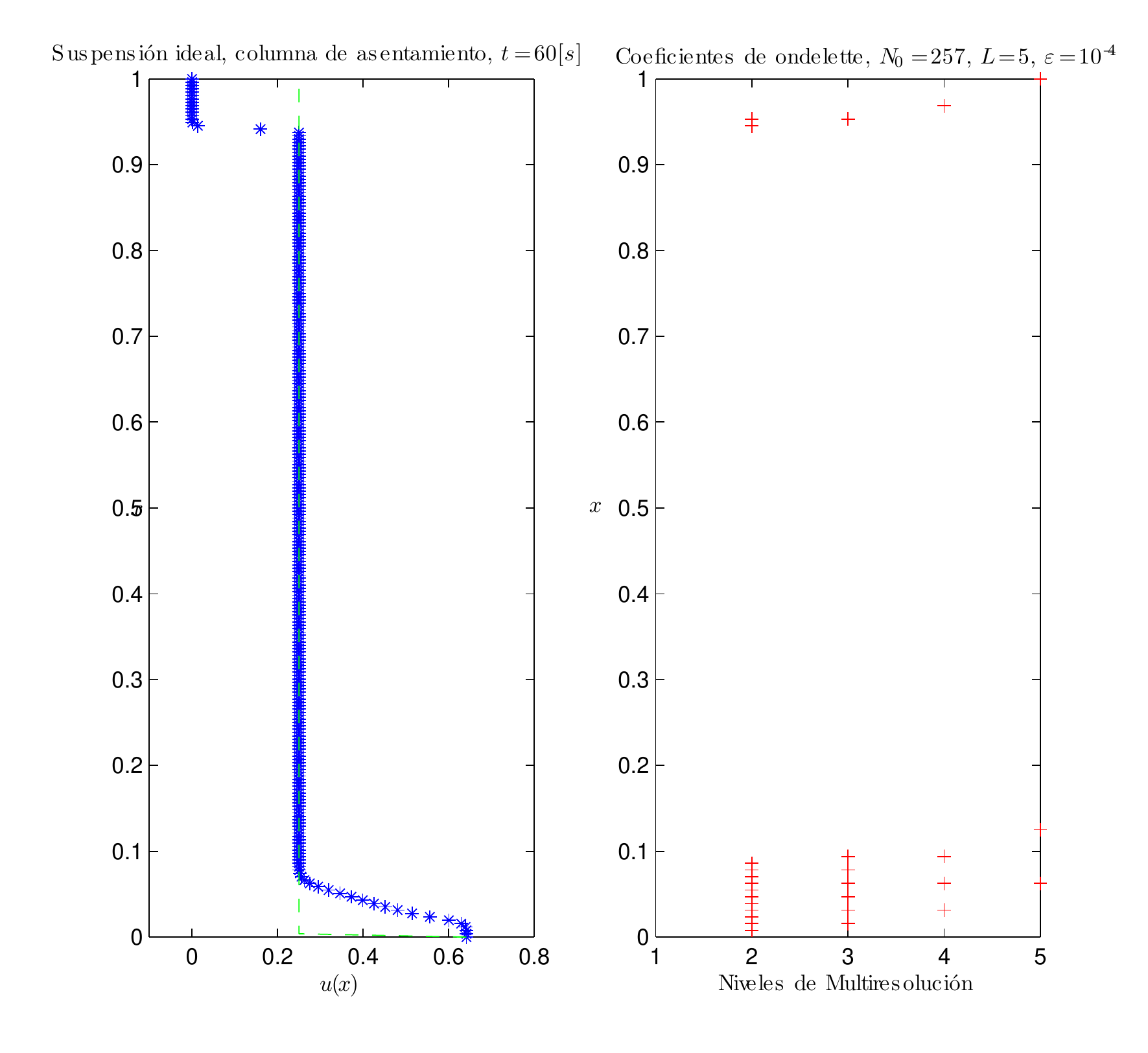}
\caption{\scriptsize Izquierda: Condici\'on inicial \emph{(rayas)} y perfil de concentraci\'on a $t=60[s]$ para el problema de sedimentaci\'on  batch de suspensi\'on ideal \emph{(Asteriscos)}. Derecha: Coeficientes de ondelette significativos correspondientes.}
\label{fig:caso0_60s}
\end{center}
\end{figure}

\begin{figure}[bp]
\begin{center}
\includegraphics[width=5in,height=3.3in]{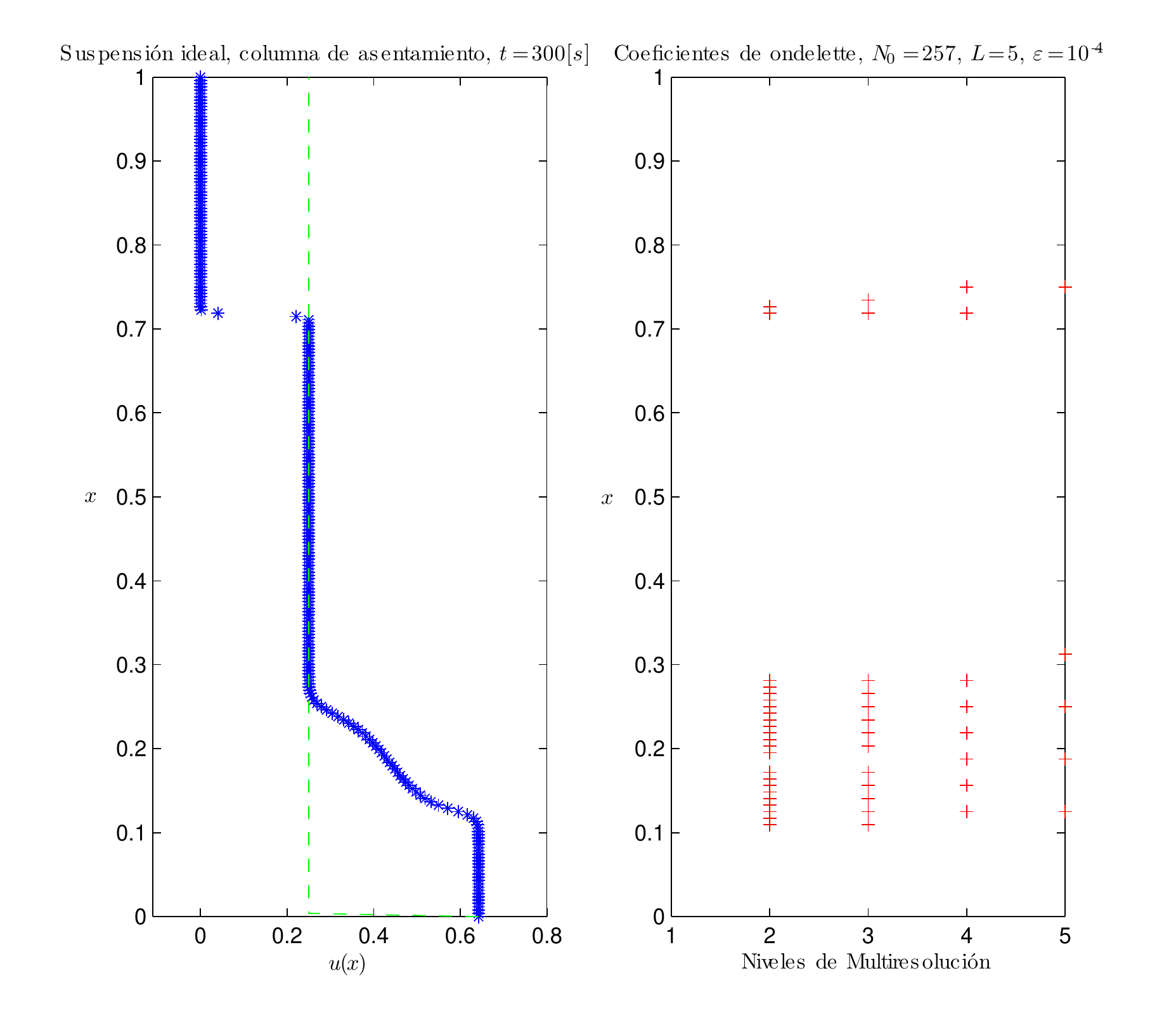}
\caption{\scriptsize Izquierda: Condici\'on inicial \emph{(rayas)} y perfil de concentraci\'on a $t=300[s]$ para el problema de sedimentaci\'on  batch de suspensi\'on ideal \emph{(Asteriscos)}. Derecha: Coeficientes de ondelette significativos correspondientes.}
\label{fig:caso0_300s}
\end{center}
\end{figure}

\begin{figure}[ht]
\begin{center}
\includegraphics[width=5in,height=3.3in]{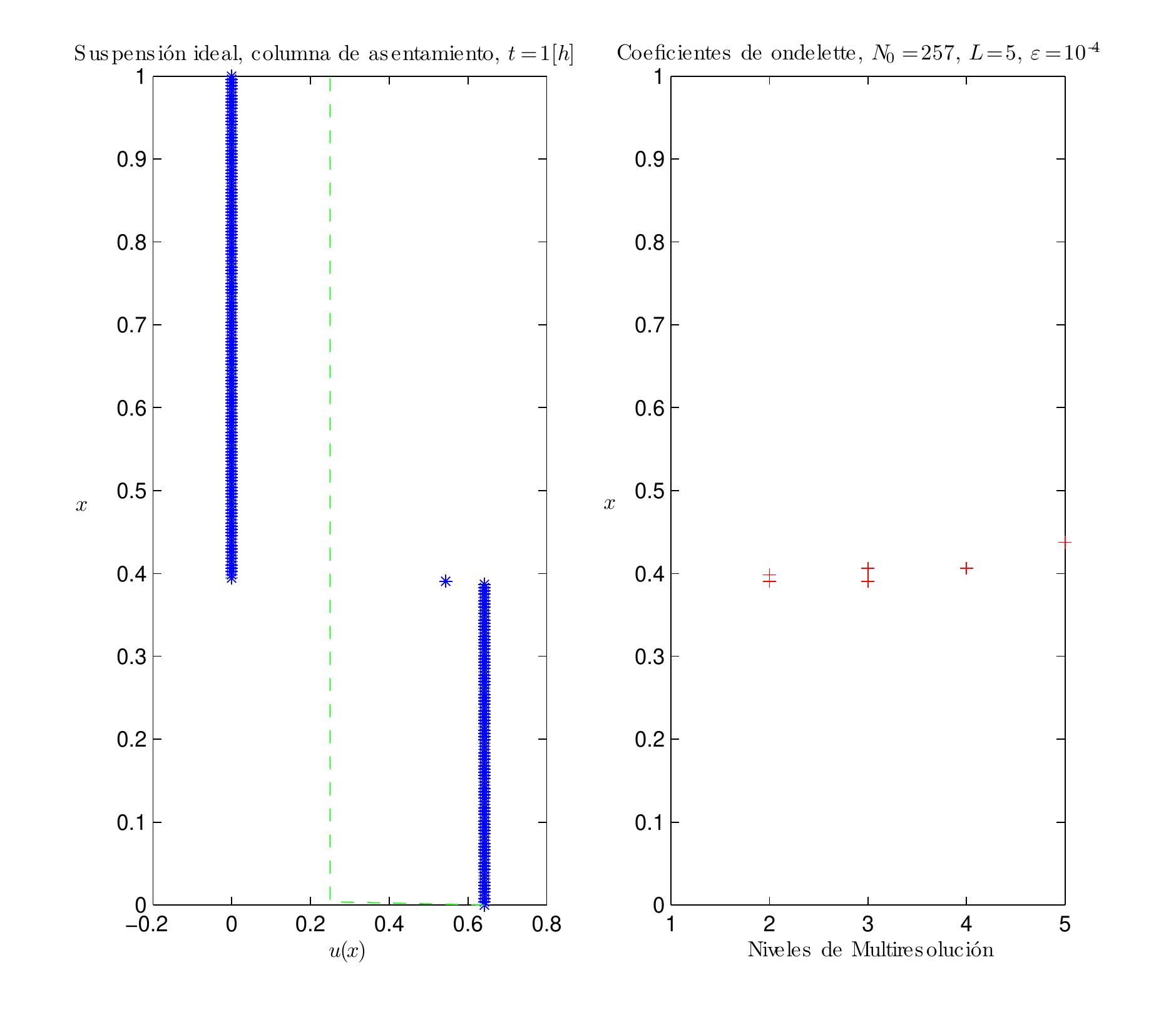}
\caption{\small Izquierda: Condici\'on inicial \emph{(rayas)} y perfil de concentraci\'on a $t=3600[s]$ para el problema de sedimentaci\'on  batch de suspensi\'on ideal \emph{(Asteriscos)}. Derecha: Coeficientes de ondelette significativos correspondientes.}
\label{fig:caso0_1hora}
\end{center}
\end{figure}

En las figuras \ref{fig:caso0_60s}-\ref{fig:caso0_1hora} se muestran soluciones num\'ericas para $t=60\, [s]$, $t=300\, [s]$ y $t=3600\, [s]$ obtenidas mediante el esquema de segundo orden descrito en la secci\'on anterior, aplicando multiresoluci\'on. En $t=3600\, [s]$ la soluci\'on ya alcanz\'o un estado estacionario.

En la tabla \ref{tabla:caso0} se muestran la proporci\'on $V$, tasa de compresi\'on y errores entre la soluci\'on calculada utilizando multiresoluci\'on y la soluci\'on calculada sin multiresoluci\'on (ver secci\'on \ref{sec:hip-num}).

\begin{table}[!h]
\begin{center}
\begin{tabular}{lccccc}
\hline
$t\, [s]$ & $V$    & $\mu$  &            $e_1$   &        $e_2$     &     $e_\infty$     \\
\hline
60  &4.3457 & 7.8456 & 2.64$\times10^{-5}$&6.54$\times10^{-6}$&9.03$\times10^{-6}$\\
300 &5.6212 & 5.8456 & 1.70$\times10^{-5}$&6.39$\times10^{-6}$&1.12$\times10^{-5}$\\
1800&5.9443 & 14.9168 & 7.28$\times10^{-5}$&2.98$\times10^{-5}$&4.35$\times10^{-5}$\\
3600&6.1385 & 29.8479& 8.89$\times10^{-5}$&4.04$\times10^{-5}$&6.50$\times10^{-5}$\\
\hline
\end{tabular}
\end{center}
\caption{\small Sedimentaci\'on de suspension ideal. $\varepsilon=1.0\times10^{-4}$, $N_0=257$ y $L=5$.}
\label{tabla:caso0}
\end{table}

Notar que los errores permanecen siempre bajo la tolerancia prescrita $\varepsilon=1.0\times10^{-4}$.
\newpage
\subsection{Caso batch de suspensiones floculadas: primer ejemplo}
En este ejemplo se considera el caso batch de suspensi\'on homog\'enea de concentraci\'on inicial $u_0(x)=0.15$ en un a columna de asentamiento cerrada, es decir, se considera el caso de $q\equiv 0$, con una concentraci\'on prescrita en $x=1$ dada por (\ref{A3}). El dominio espacial es $[0,1]$ y la concentraci\'on cr{\'i}tica es $u_c=0.23$.  Notar que la discontinuidad entre $u=0$ y $u=u_0$ es un choque. A\'un m\'as, el problema (\ref{A1})-(\ref{A4}) es un problema de Riemann, en el sentido de que el dato inicial consiste en dos estados constantes y la soluci\'on, en general, consistir\'a en ondas elementales: choques, ondas de rarefacci\'on y discontinuidades de contacto \cite{Thick}.

Como funci\'on de densidad de flujo, se utiliza una funci\'on Kynch batch Richardson-Zaki con par\'ametros correspondientes a suspensi\'on de cobre \cite{BEKL}.
\begin{equation}\label{zaki}
f(u)=-6.05\times10^{-4}u(1-u)^{12.59}\, [m/s].
\end{equation}
Se utilizar\'a la funci\'on $\sigma_e'(u)$ dada por (\cite{BK,Thick})
\begin{equation}
\sigma_e'(u)=\frac{d}{d\,u}\left(100(u/u_c)^8-1\right)[Pa], \textrm{ si } u> u_c.
\end{equation}
Luego
\begin{equation}
\sigma_e'(u)=\left\{\begin{array}{ll}
0,&\textrm{ si } u\leqslant u_c=0.23,\\
\frac{800}{u_c}\left(\frac{u}{u_c}\right)^7\, [Pa],& \textrm{ si } u> u_c.
\end{array}\right.
\end{equation}

La funci\'on $a(u)$ (\ref{a(u)}) est\'a dada entonces por
\begin{equation}
a(u)=\left\{\begin{array}{ll}
0,&\textrm{ si } u\leqslant u_c=0.23,\\
\frac{4.84\times10^{-1}u^7(1-u)^{12.59}}{u_c^8\Delta\varrho\, g},& \textrm{ si } u> u_c,
\end{array}\right.
\end{equation}
con $\Delta\varrho=1500\, [Kg/m^3]$ y $g=9.81\, [Kg\, m/s^2]$.

\begin{figure}[!h]
\begin{center}
\includegraphics[width=6in,height=2.5in]{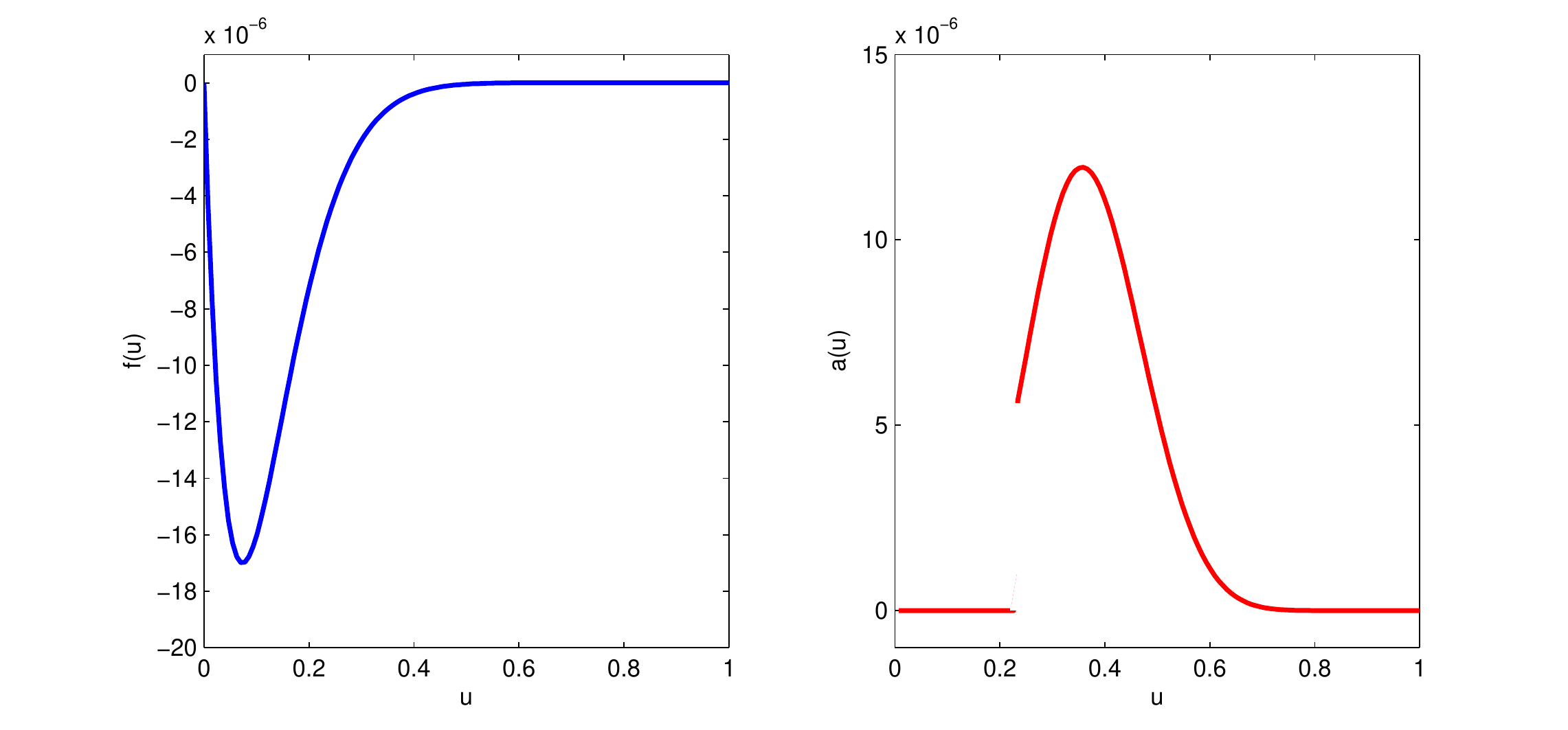}
\caption[Funciones modelo para el \emph{problema de sedimentaci\'on-consolidaci\'on}.]{\small Funciones modelo $f(u)$ (izquierda) y $a(u)$ (derecha) para el problema de sedimentaci\'on-consolidaci\'on. Las unidades son $[m/s]$ para $f(u)$ y $[m^2/s]$ para $a(u)$.}
\label{fig:bek1}
\end{center}
\end{figure}

La figura \ref{fig:bek1} muestra las funciones modelo $f(u)$ y $a(u)$. La funci\'on $A(u)$ correspondiente al t\'ermino difusivo integrado, se calcula mediante las f\'ormulas (\ref{dif1})-(\ref{dif2}).

En la tabla \ref{tabla:caso1} se muestran la proporci\'on $V$, tasa de compresi\'on y errores entre la soluci\'on obtenida utilizando multiresoluci\'on y la soluci\'on obtenida sin multiresoluci\'on.

Notar de la tabla \ref{tabla:caso1}, que los errores se encuentran
por debajo de la tolerancia prescrita. Notar adem\'as los
excelentes resultados en cuanto a proporci\'on $V$
(correspondiente al tiempo total de CPU en ambos casos). Los
resultados en cuanto a tasa de compresi\'on no son excelentes,
pero hay que tomar en cuenta que se est\'a considerando una malla
de 129 puntos.

\begin{table}[!h]
\begin{center}
{\small
\begin{tabular}{lccccc}
\hline
$t\, [s]$ & $V$    & $\mu$  &            $e_1$   &        $e_2$     &     $e_\infty$     \\
\hline
60       &6.5737 & 17.8796& 1.29$\times10^{-4}$&8.72$\times10^{-5}$&5.33$\times10^{-5}$\\
1800 (*) &5.7349 & 9.4132 & 1.99$\times10^{-4}$&9.06$\times10^{-5}$&7.42$\times10^{-5}$\\
3600 (*) &6.1982 & 9.1246 & 2.77$\times10^{-4}$&2.67$\times10^{-4}$&9.61$\times10^{-5}$\\
7200 (*) &6.2110 & 9.1246 & 3.21$\times10^{-4}$&4.67$\times10^{-4}$&2.41$\times10^{-4}$\\
14400(*) &7.9244 & 9.4132 & 8.92$\times10^{-4}$&7.81$\times10^{-4}$&6.18$\times10^{-4}$\\
\hline
\end{tabular}}
\end{center}
\caption{Suspensiones floculadas, primer ejemplo.
Multiresoluci\'on utilizando $\varepsilon=10^{-3}$,
$N_0=129$ y $L=5$. (*): figuras \ref{fig:caso1_mediahora} -
\ref{fig:caso1_4horas}.} \label{tabla:caso1}
\end{table}

\begin{figure}[!h]
\begin{center}
\includegraphics[width=6in,height=2.7in]{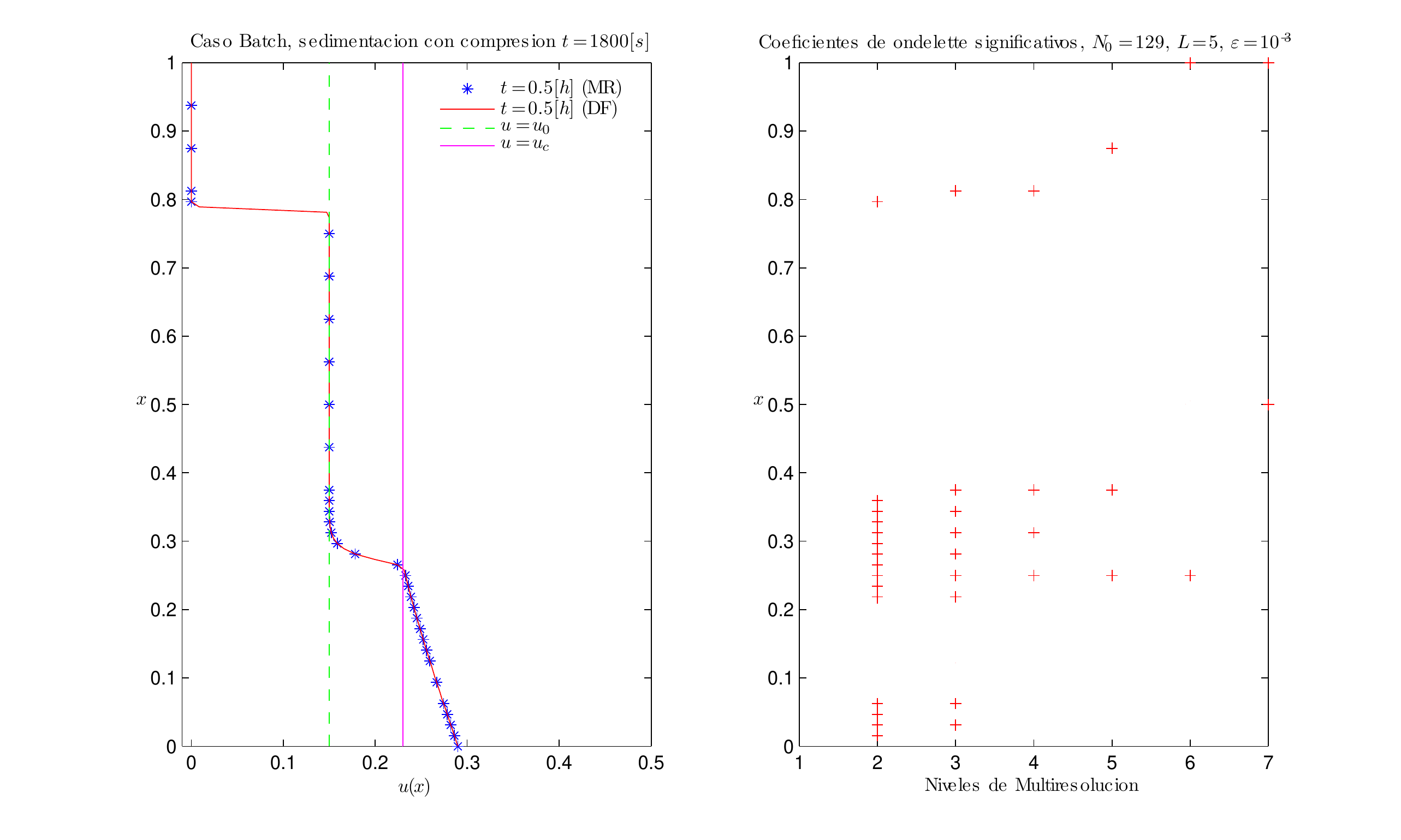}
\caption{\small Izquierda: Condici\'on inicial \emph{(rayas)} y perfil de
concentraci\'on a $t=1800[s]$ para el problema de
sedimentaci\'on-consolidaci\'on  \emph{(asteriscos)}. Derecha:
Coeficientes de ondelette significativos correspondientes.
$\varepsilon=10^{-3}$, $N_0=129$ y $L=5$.}
\label{fig:caso1_mediahora}
\end{center}
\end{figure}

En la figura \ref{fig:caso1_mediahora} se presenta un perfil de
concentraci\'on en un tiempo $t=1800[s]$, utilizando
multiresoluci\'on. La soluci\'on se calcula utilizando $129$
puntos en la malla fina, con una estrategia de truncamiento
$\varepsilon_k=\frac{\varepsilon}{2^{L-k}}$. Se presenta adem\'as
la configuraci\'on de los coeficientes de ondelette significativos. Notar que cuanto m\'as perfilada es la discontinuidad, menor es el n\'umero de coeficientes de ondelette
 significativos asociados a tal discontinuidad.

\begin{figure}[!h]
\begin{center}
\includegraphics[width=6in,height=2.7in]{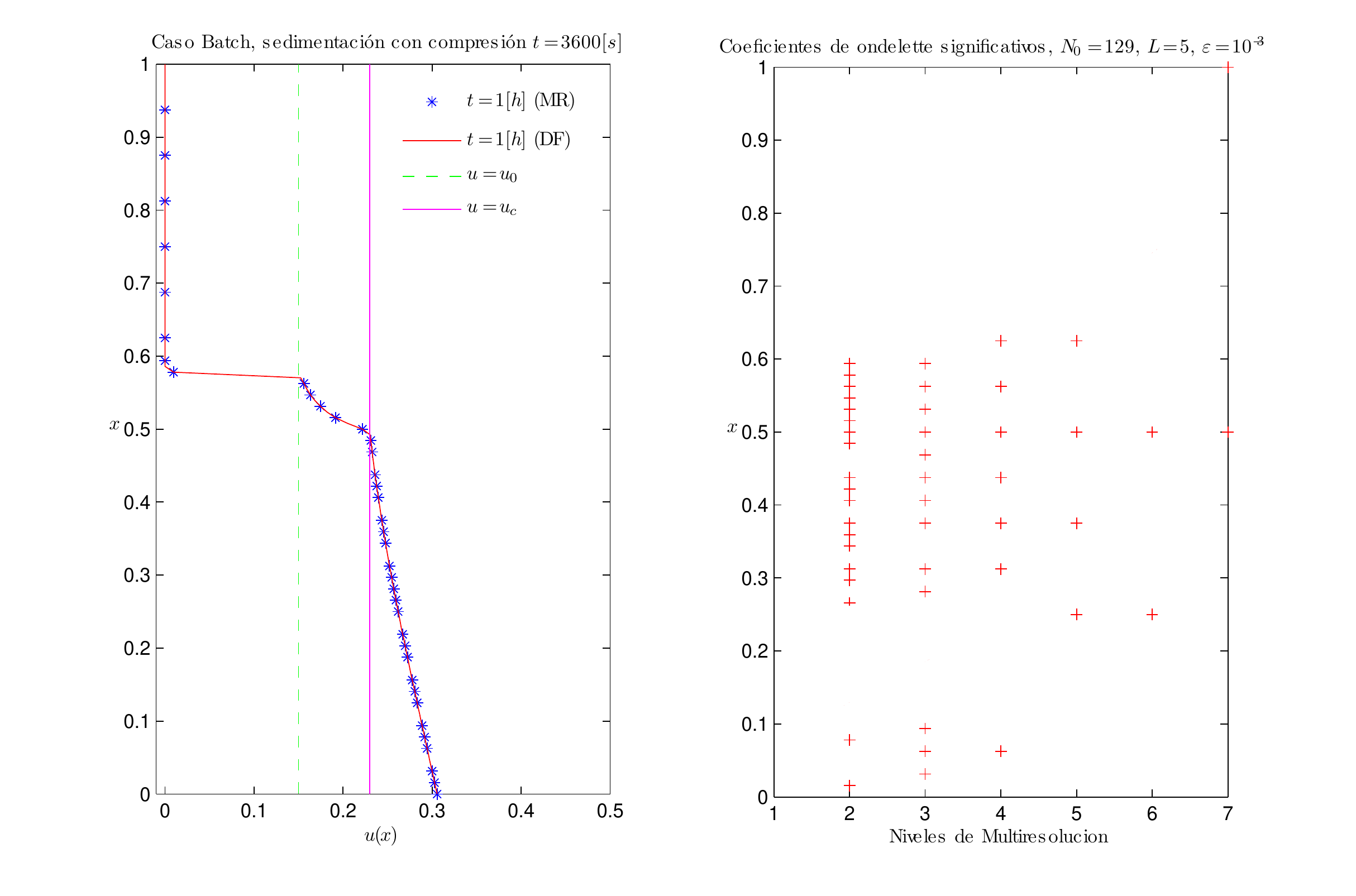}
\caption{\small Izquierda: Condici\'on inicial \emph{(rayas)} y perfil de
concentraci\'on a $t=3600[s]$ para el problema de
sedimentaci\'on-consolidaci\'on  \emph{(asteriscos)}. Derecha:
Coeficientes de ondelette significativos correspondientes.
$\varepsilon=10^{-3}$, $N_0=129$ y $L=5$.}
\label{fig:caso1_1hora}
\end{center}
\end{figure}

En la figura \ref{fig:caso1_1hora} se presenta un perfil de
concentraci\'on en un tiempo $t=1[h]$,  utilizando
multiresoluci\'on. La soluci\'on se calcula utilizando $129$
puntos en la malla fina. Se presenta adem\'as la configuraci\'on
de los coeficientes de ondelette significativos correspondientes.

\begin{figure}[!ht]
\begin{center}
\includegraphics[width=6in,height=2.5in]{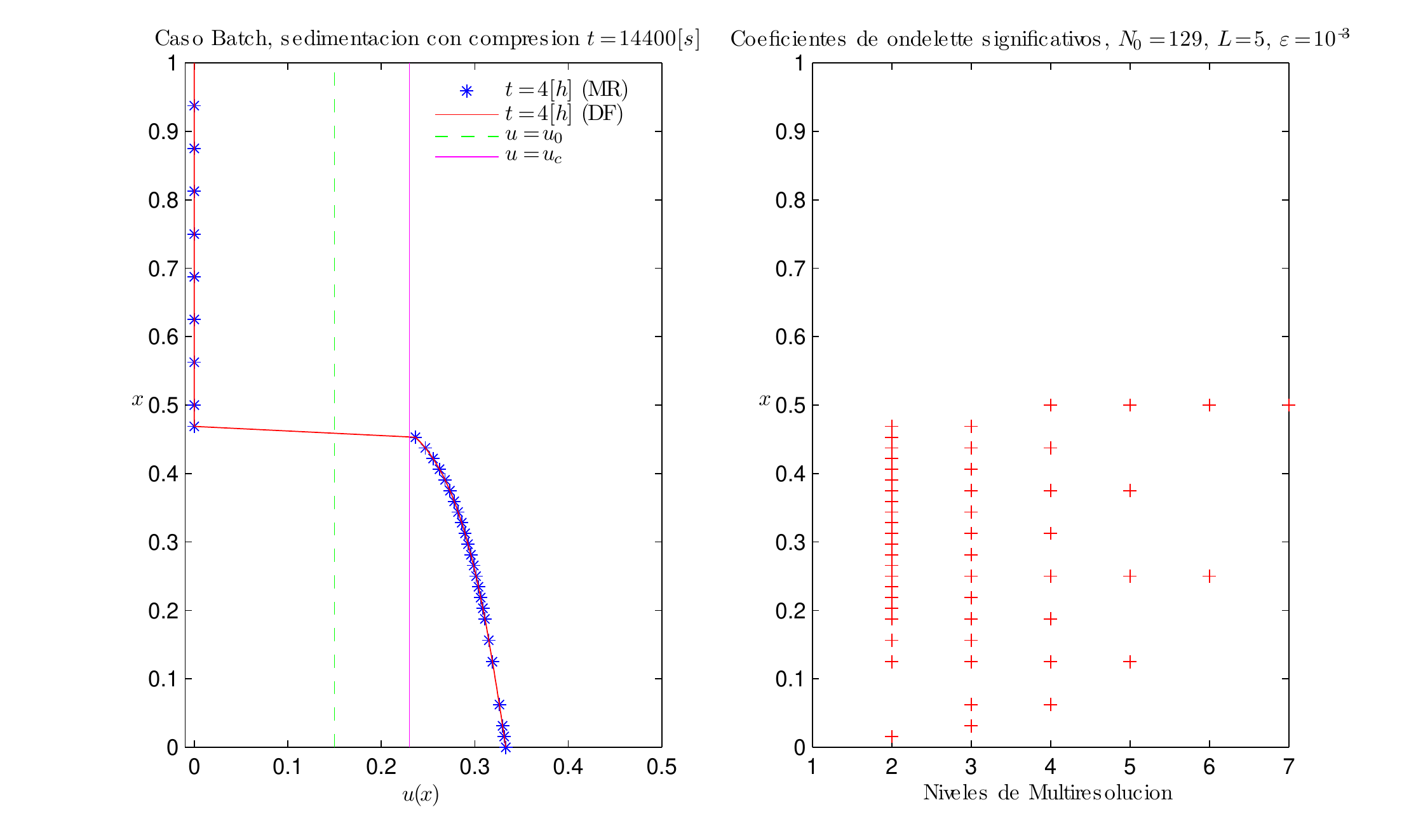}
\caption[Perfil de concentraci\'on.
Sedimentaci\'on-consolidaci\'on, primer ejemplo
$t=14400$.]{\small Izquierda: Condici\'on inicial \emph{(rayas)} y perfil
de concentraci\'on a $t=4[h]$ para el problema de
sedimentaci\'on-consolidaci\'on, caso \emph{Batch. (asteriscos)}.
Derecha: Coeficientes de ondelette significativos
correspondientes. $\varepsilon=10^{-3}$, $N_0=129$ y
$L=5$.} \label{fig:caso1_4horas}
\end{center}
\end{figure}

\begin{figure}[!ht]
\begin{center}
\includegraphics[width=5.2in,height=2.5in]{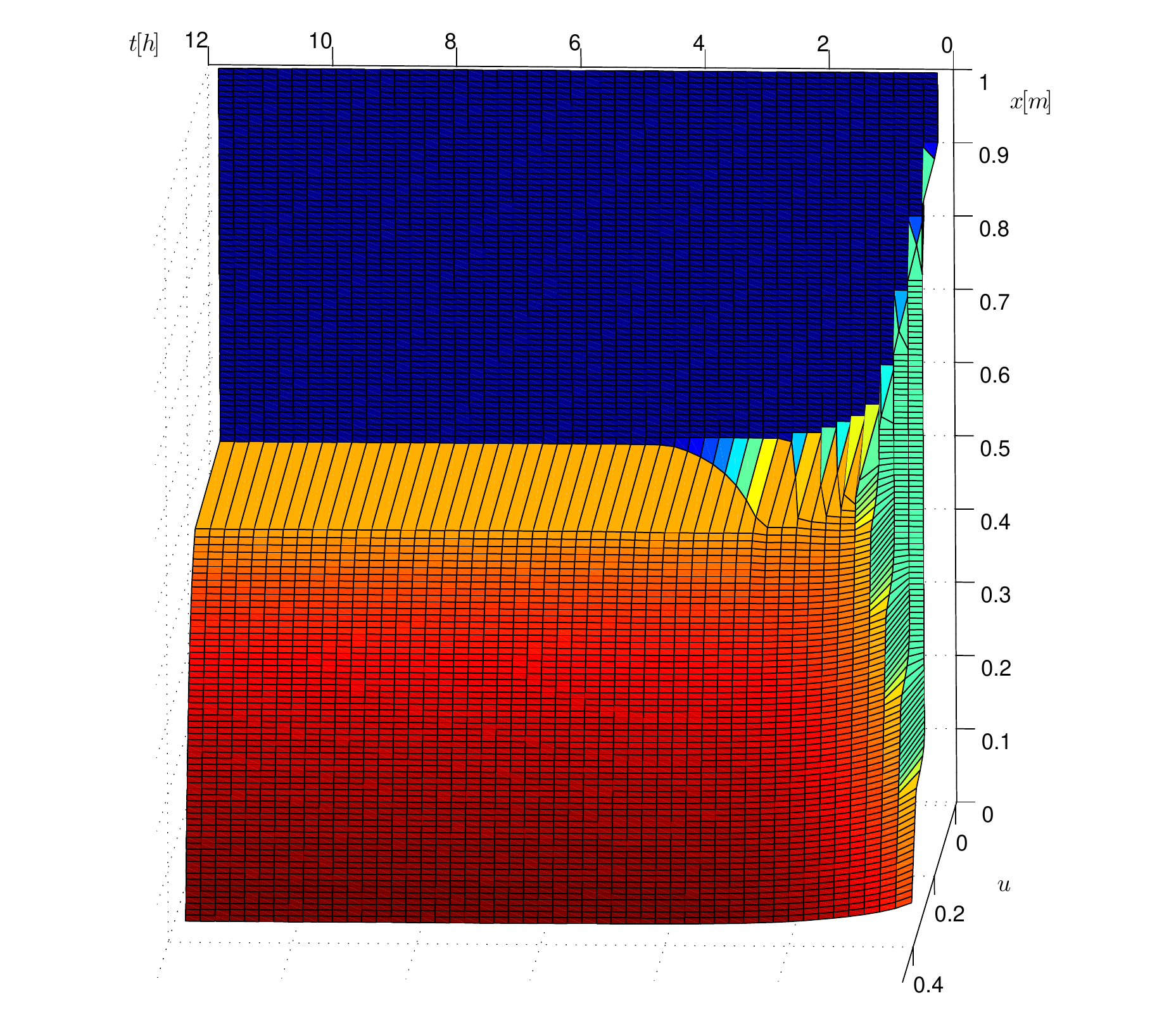}
\caption[Perfiles de concentraci\'on. Sedimentaci\'on-consolidaci\'on, primer ejemplo hasta $t=43200$.]{\small Perfiles de concentraci\'on hasta $t=12[h]$ para el problema de sedimentaci\'on-consolidaci\'on, caso \emph{Batch.}  $\varepsilon=10^{-3}$, $N_0=129$ y $L=5$.}
\label{fig:caso1_plot}
\end{center}
\end{figure}

En la figura \ref{fig:caso1_4horas} se presenta un perfil de
concentraci\'on en un tiempo $t=4[h]$,  utilizando
multiresoluci\'on. La soluci\'on se calcula utilizando $129$
puntos en la malla fina. Se presenta adem\'as la configuraci\'on
de los coeficientes de ondelette significativos correspondientes.
En este tiempo, la soluci\'on ya se encuentra en un estado
estacionario (ver adem\'as \ref{fig:caso1_plot} y \cite{BEK}).

Finalmente se presenta en la figura \ref{fig:caso1_plot} la
soluci\'on num\'erica del problema de sedimentaci\'on
consolidaci\'on en asentamiento tipo batch, hasta el tiempo $t=12[h]$.

Los resultados num\'ericos concuerdan con los resultados obtenidos por  B\"urger \emph{et al.} \cite{BEK}.
\newpage
\subsection{Caso batch de suspensiones floculadas: segundo ejemplo}
En este ejemplo se considera el caso batch de suspensi\'on
homog\'enea de concentraci\'on inicial $u_0(x)=0.05$ en una
columna de asentamiento cerrada ($q\equiv 0$) de menor longitud:
$H=0.16 [m]$ (ver \cite{BK}). La concentraci\'on cr{\'i}tica es
$u_c=0.07$.

Como funci\'on de densidad de flujo, se utiliza la funci\'on
$f(u)$ dada por (\ref{laf}) y como funci\'on de rigidez s\'olida
efectiva, se utiliza la funci\'on $\sigma_e'(u)$ dada por
(\ref{lasigma}), donde los par\'ametros necesarios
\begin{equation}\label{parameters}
v_\infty=-2.7\times 10^{-4}[ms^{-1}],\ C=21.5,\ u_{max}=0.5,\ \sigma_0=5.7[Pa] \textrm{ y }n=5,
\end{equation}
corresponden al modelo de suspensi\'on con compresi\'on tipo Kaolin (ver \cite{BK}). Adem\'as $\Delta\varrho=1690\, [Kg/m^3]$ y $g=9.81\, [Kg\, m/s^2]$.

\begin{figure}[!ht]
\begin{center}
\includegraphics[width=6in,height=2.5in]{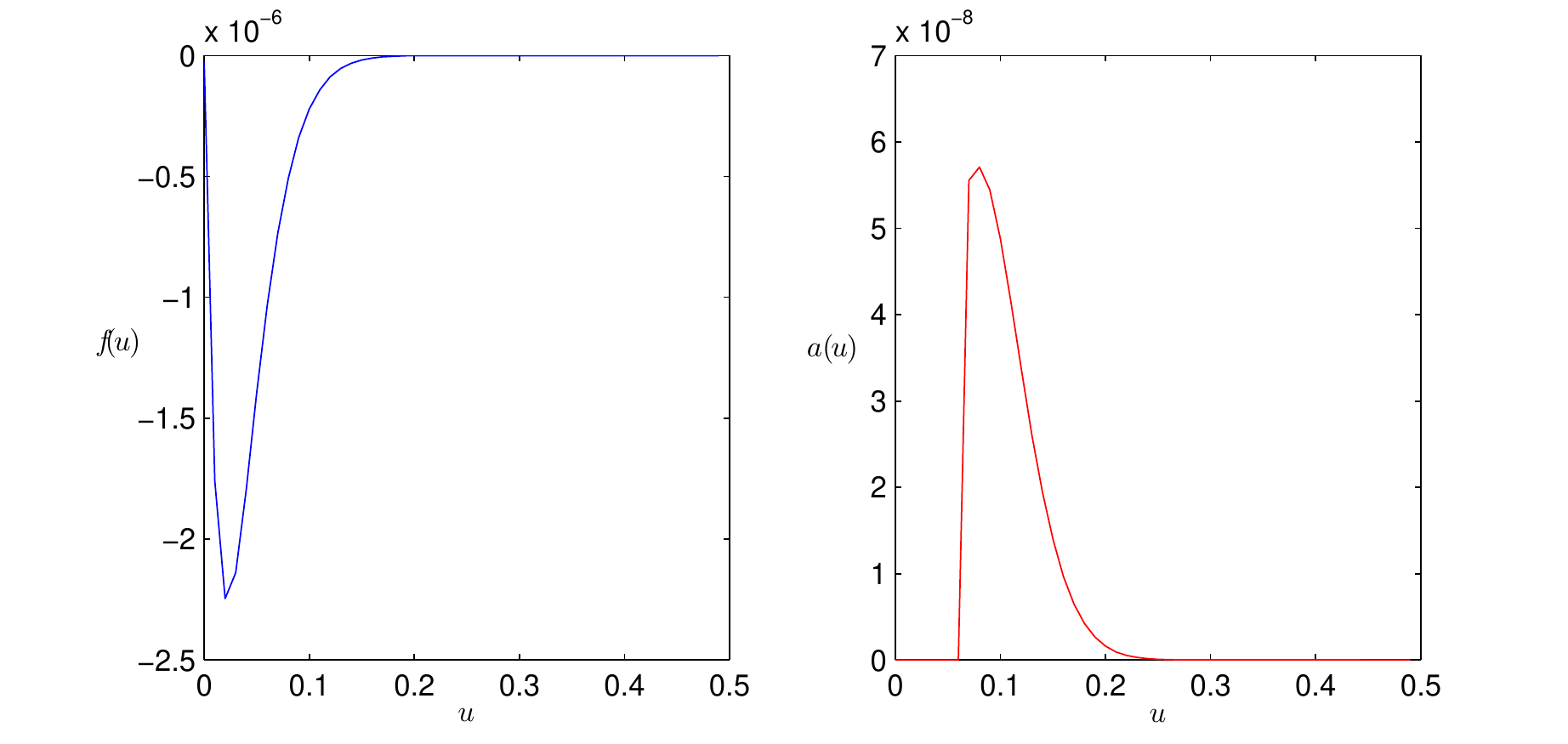}
\caption[Funciones modelo para el \emph{problema de sedimentaci\'on-consolidaci\'on}.]{\small Funciones modelo $f(u)$ (izquierda) y $a(u)$ (derecha) para el problema de sedimentaci\'on-consolidaci\'on, segundo ejemplo. Las unidades son $[m/s]$ para $f(u)$ y $[m^2/s]$ para $a(u)$.}
\label{fig:bek2}
\end{center}
\end{figure}

La figura \ref{fig:bek2} muestra las funciones modelo $f(u)$ y $a(u)$ para este caso.

En \cite{BK} se da la siguiente expresi\'on para el t\'ermino difusivo integrado:
\begin{equation}\label{dif1}
A(u)=\left\{\begin{array}{ll}
0,&\textrm{ si } u\leqslant u_c,\\
\mathcal{A}(u)-\mathcal{A}(u_c),&\textrm{ si } u>u_c,
\end{array}\right.
\end{equation}
donde
\begin{equation}\label{dif2}
\mathcal{A}(u)=\frac{v_\infty\sigma_0}{\Delta\varrho g u_c^n}\left(1-\frac{u}{u_{\max}}\right)^Cu^n\sum_{j=1}^n\left(\prod_{l=i}^j\frac{n+1-l}{C+l}\right)\left(\frac{u_{\max}}{u}-1\right)^j,
\end{equation}
cuya gr\'afica se muestra a continuaci\'on.

\begin{figure}[!h]
\begin{center}
\includegraphics[width=4.5in,height=2.5in]{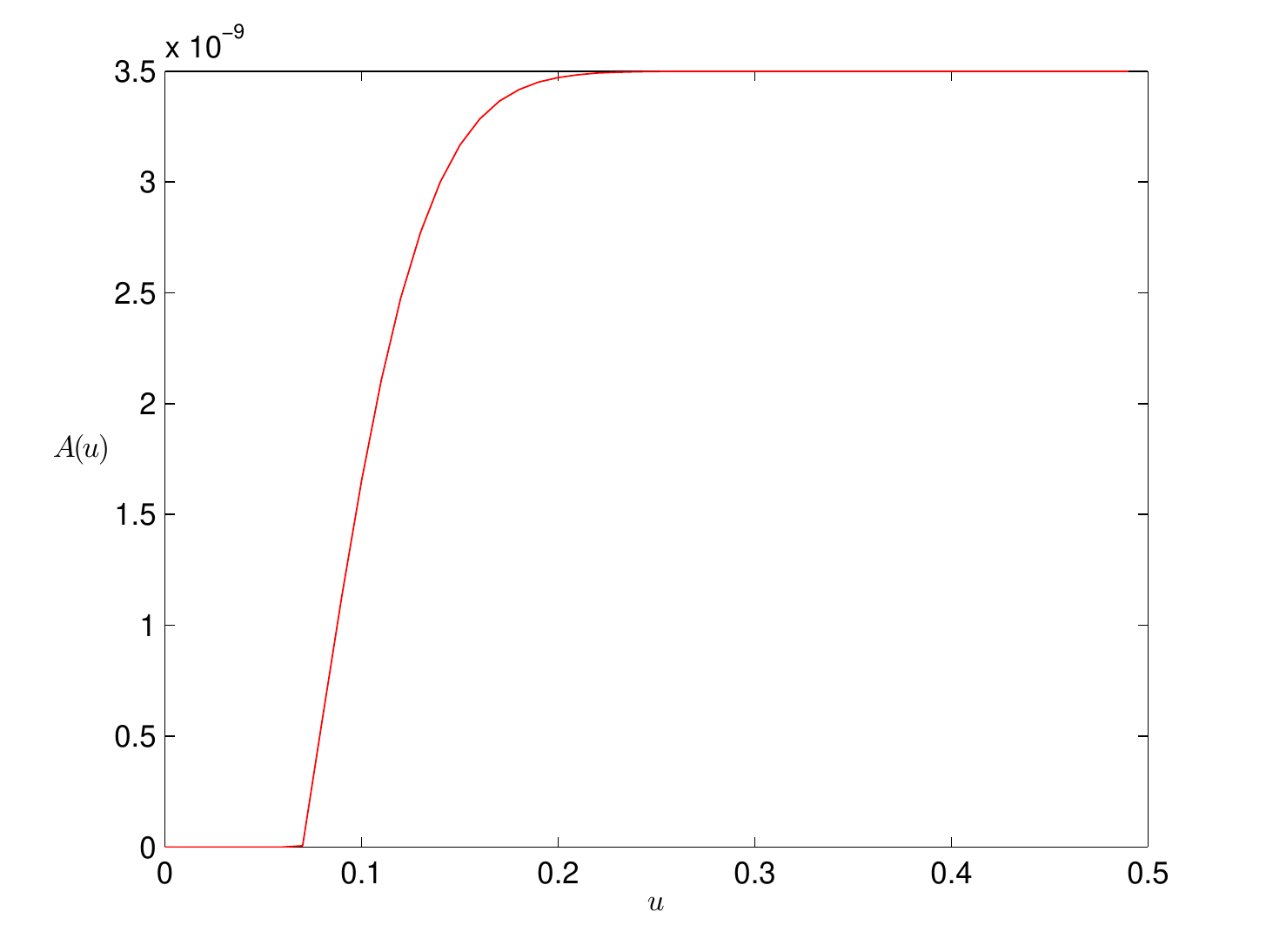}
\caption[T\'ermino  difusivo integrado, segundo \emph{problema de sedimentaci\'on-consolidaci\'on}.]{\small T\'ermino  difusivo integrado $A(u)$ para el problema de sedimentaci\'on-consolidaci\'on, segundo ejemplo.}
\label{fig:difu2}
\end{center}
\end{figure}

En la tabla \ref{tabla:caso2} se muestran la proporci\'on $V$, tasa de compresi\'on y errores entre la soluci\'on obtenida utilizando multiresoluci\'on y la soluci\'on obtenida sin multiresoluci\'on.

\begin{table}[!h]
\begin{center}
{\small
\begin{tabular}{lccccc}
\hline
$t\, [s]$ &  $V$    & $\mu$  &            $e_1$   &        $e_2$     &     $e_\infty$     \\
\hline
60        &1.4109 & 5.6100 & 4.31$\times10^{-5}$&2.34$\times10^{-4}$&1.46$\times10^{-4}$\\
2000  (*) &4.4782 & 7.1542 & 6.87$\times10^{-5}$&5.78$\times10^{-4}$&7.88$\times10^{-4}$\\
6000  (*) &7.2384 & 10.7245 & 1.36$\times10^{-4}$&9.45$\times10^{-4}$&9.65$\times10^{-4}$\\
10000 (*) &10.4568& 10.9781 & 6.74$\times10^{-4}$&1.32$\times10^{-3}$&1.03$\times10^{-3}$\\
\hline
\end{tabular}}
\end{center}
\caption{\small Caso batch de suspensiones floculadas, segundo ejemplo. Tolerancia prescrita $\varepsilon=10^{-3}$, $N_0=129$ puntos en la malla fina y $L=5$ niveles de multiresoluci\'on. (*): figuras \ref{fig:caso2_2000s} - \ref{fig:caso2_10000s}.}
\label{tabla:caso2}
\end{table}

An\'alogamente al primer ejemplo, en la tabla \ref{tabla:caso2} puede verse que los errores entre la soluci\'on obtenida utilizando multiresoluci\'on y la soluci\'on obtenida sin multiresoluci\'on, est\'an por debajo de la tolerancia prescrita. De igual modo, se ve una gran rebaja en costo computacional, dada por la alta tasa de compresi\'on y proporci\'on $V$.

\begin{figure}[!ht]
\begin{center}
\includegraphics[width=6in,height=2.7in]{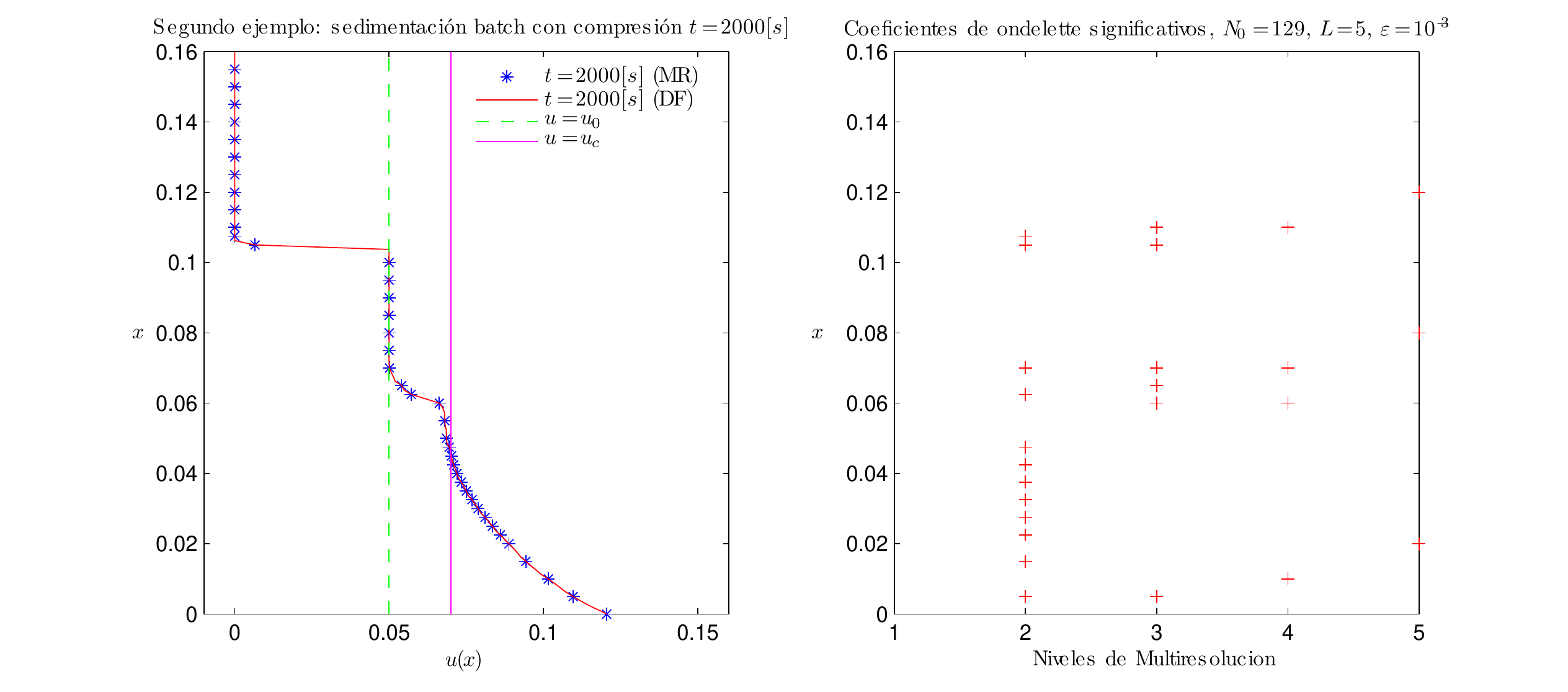}
\caption[Perfil de concentraci\'on. Sedimentaci\'on-consolidaci\'on, $t=2000$.]{\small Izquierda: Condici\'on inicial \emph{(rayas)} y perfil de concentraci\'on a $t=2000[s]$ para el problema de sedimentaci\'on-consolidaci\'on, segundo caso \emph{(asteriscos)}. Derecha: Coeficientes de ondelette significativos correspondientes. $\varepsilon=10^{-3}$, $N_0=129$ y $L=5$.}
\label{fig:caso2_2000s}
\end{center}
\end{figure}

En la figura \ref{fig:caso2_2000s} se presenta un perfil de concentraci\'on en $t=2000[s]$, para la soluci\'on utilzando multiresoluci\'on, y la soluci\'on sin multiresoluci\'on. La soluci\'on se calcula utilizando $129$ puntos en la malla fina. Se presenta adem\'as la configuraci\'on correspondiente de los coeficientes de ondelette significativos.

\newpage
\begin{figure}[!h]
\begin{center}
\includegraphics[width=6in,height=2.7in]{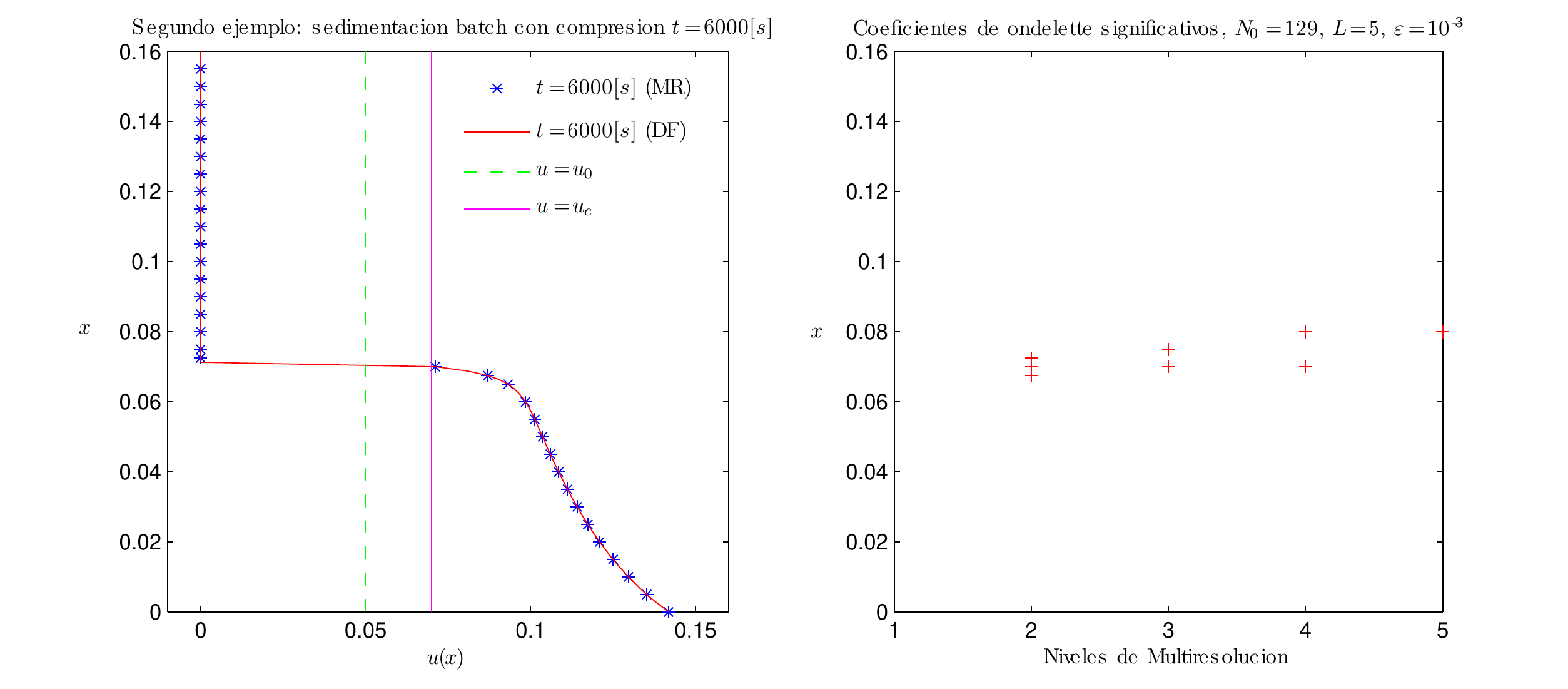}
\caption[Perfil de concentraci\'on. Sedimentaci\'on-consolidaci\'on, $t=6000$.]{\small Izquierda: Condici\'on inicial \emph{(rayas)} y perfil de concentraci\'on a $t=6000[s]$ para el problema de sedimentaci\'on-consolidaci\'on, segundo caso \emph{(asteriscos)}. Derecha: Coeficientes de ondelette significativos correspondientes. $\varepsilon=10^{-3}$, $N_0=129$ y
$L=5$.}
\label{fig:caso2_6000s}
\end{center}
\end{figure}

\begin{figure}[!ht]
\begin{center}
\includegraphics[width=6in,height=2.5in]{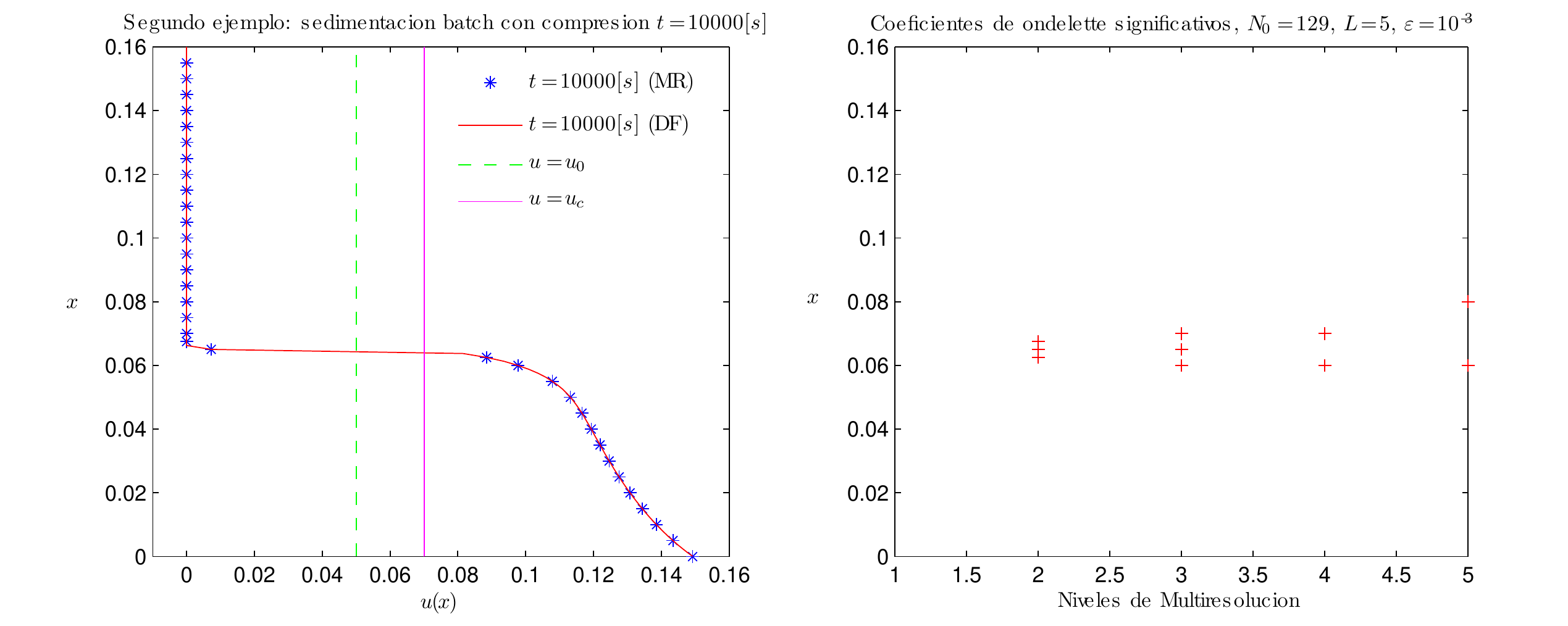}
\caption[Perfil de concentraci\'on. Sedimentaci\'on-consolidaci\'on, $t=10000$.]{\small Izquierda: Condici\'on inicial \emph{(rayas)} y perfil de concentraci\'on a $t=10000[s]$ para el problema de sedimentaci\'on-consolidaci\'on, segundo caso \emph{(asteriscos)}. Derecha: Coeficientes de ondelette significativos correspondientes. $\varepsilon=10^{-3}$, $N_0=129$ y
$L=5$.}
\label{fig:caso2_10000s}
\end{center}
\end{figure}

\small En las figuras \ref{fig:caso2_6000s} y \ref{fig:caso2_10000s} se presentan perfiles de concentraci\'on en tiempos $t=6000[s]$ y $t=10000[s]$, utilizando multiresoluci\'on, y la configuraci\'on de coeficientes de ondelette significativos.

\newpage

\begin{figure}[!ht]
\begin{center}
\includegraphics[width=4in,height=3in]{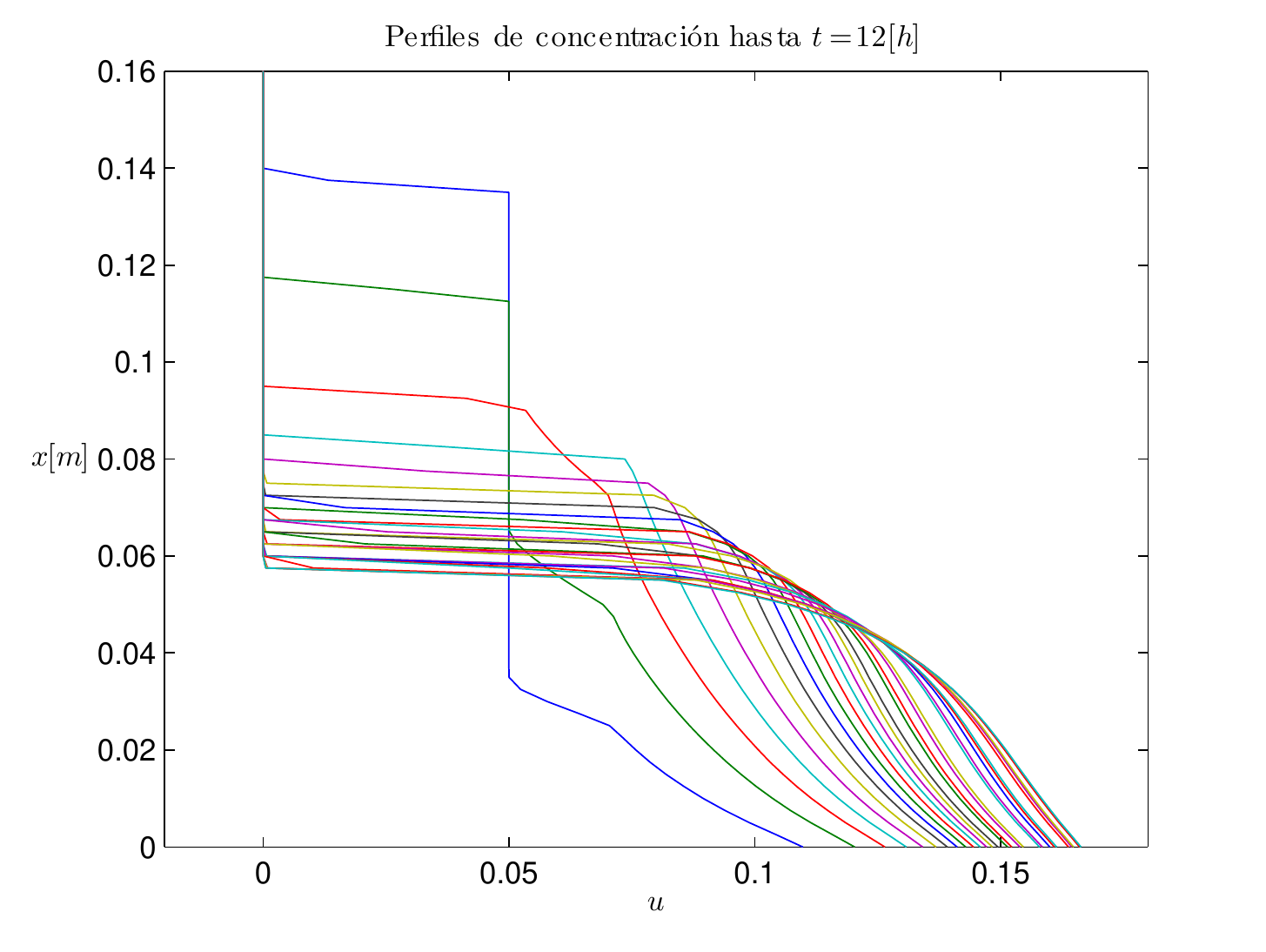}
\caption[Perfiles de concentraci\'on. Sedimentaci\'on-consolidaci\'on, segundo ejemplo hasta $t=43200$.]{\small Perfiles de concentraci\'on hasta $t=12[h]$ para el segundo problema de sedimentaci\'on-consolidaci\'on, asentamiento Batch. $\varepsilon=10^{-3}$, $N_0=129$ y
$L=5$.}
\label{fig:caso2_plot}
\end{center}
\end{figure}

Finalmente se presenta en la figura \ref{fig:caso2_plot} la soluci\'on num\'erica del problema de sedimentaci\'on consolidaci\'on obtenida utilizando el m\'etodo de multiresoluci\'on, hasta el tiempo $t=12[h]$.

Los resultados num\'ericos concuerdan con los resultados obtenidos por  B\"urger y Karlsen \cite{BK}.
\normalsize
\newpage
\subsection{Simulaci\'on de sedimentaci\'on continua}
Se modela un ICT de longitud 2, con una concentraci\'on inicial, $u_0=0.052$. En $x=1$ se prescribe una alimentaci\'on dada por $\Psi(t)=-8.55\times10^{-7}$. Se supone el ICT cerrado, es decir, $q\equiv0$ y se simula el proceso de llenado hasta antes que el nivel de concentraci\'on en $x=0$ alcance el valor $u(0,t)=0.171$. En ese momento, el recipiente se abre, y se hace $q(t)=-5\times10^{-6}[m/s]$. Notar que desde ese momento, $\Psi(t)=0.171\cdot q(t)$, es decir, el flujo en la alimentaci\'on es igual al flujo de descarga y el perfil de concentraci\'on entra en estado constante \cite{BK}.

Notar que en este caso se utiliza como modelo el problema B (\ref{B1})-(\ref{B4}). Se utiliza una funci\'on de densidad de flujo dada por
\begin{equation}
f(u)=-1.98\times10^{-4}u\left(1-\frac{u}{0.3}\right)^{5.647},
\end{equation}
y una funci\'on de rigidez s\'olida efectiva dada por
\begin{equation}
\sigma_e(u)=\left\{\begin{array}{ll}
0,&\textrm{ si }u\leqslant u_c:=0.1,\\
5.7\left[\left(\frac{u}{u_c}\right)^9-1\right],&\textrm{ si }u> u_c:=0.1.\end{array}\right.
\end{equation}
Estas aproximan a las funciones modelo determinadas para suspensi\'on de carbonato de calcio \cite{BK}.

\begin{figure}[h]
\begin{center}
\includegraphics[width=6in,height=2.5in]{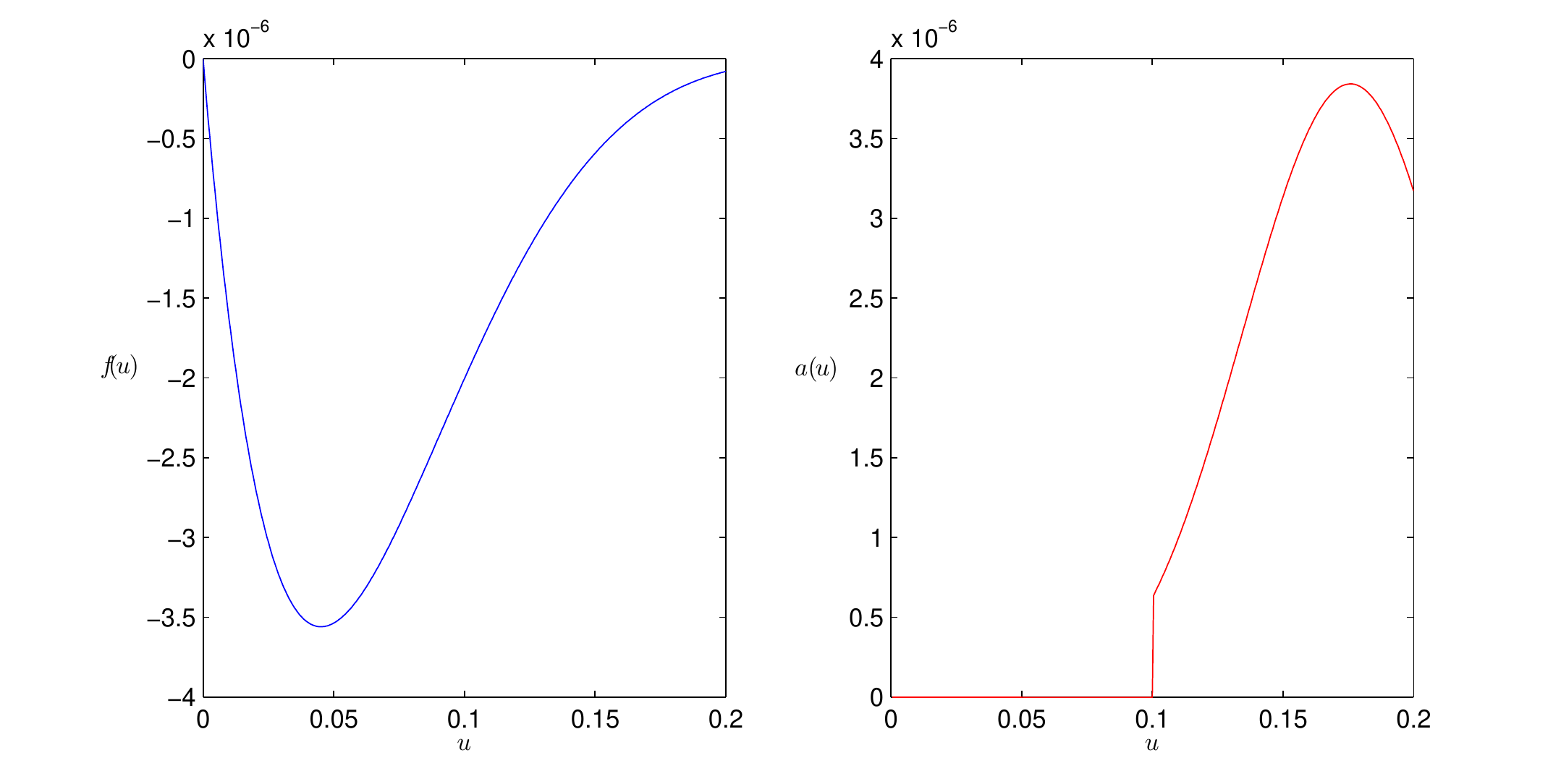}
\caption[Funciones modelo para la simulaci\'on de \emph{sedimentaci\'on continua}.]{\small Funciones modelo $f(u)$ (izquierda) y $a(u)$ (derecha) para para la simulaci\'on de sedimentaci\'on continua. Las unidades son $[m/s]$ para $f(u)$ y $[m^2/s]$ para $a(u)$.}
\label{fig:bek3}
\end{center}
\end{figure}

\begin{figure}[!h]
\begin{center}
\includegraphics[width=4.5in,height=2.5in]{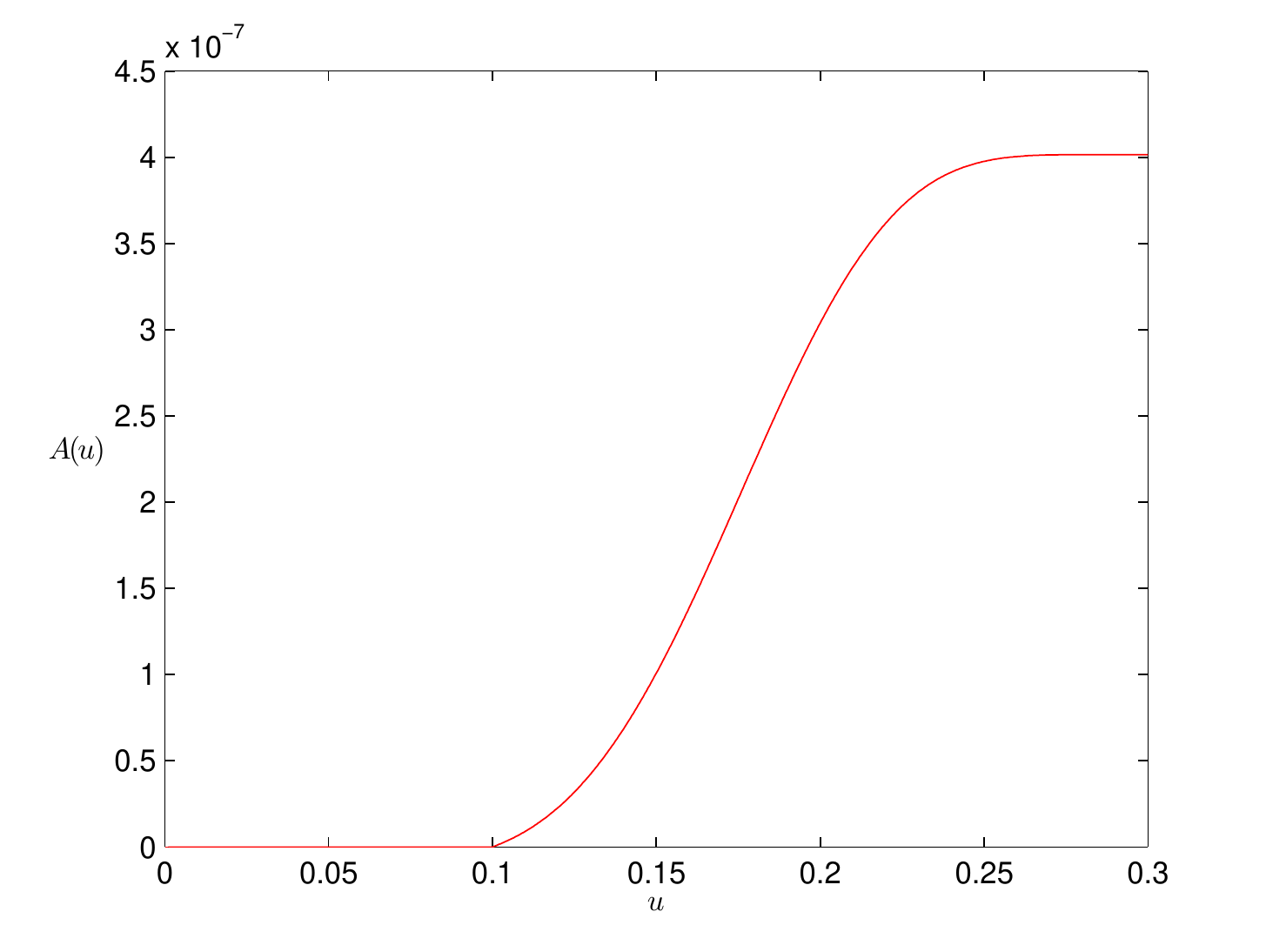}
\caption[T\'ermino  difusivo integrado, segundo \emph{problema de sedimentaci\'on-consolidaci\'on}.]{\small T\'ermino difusivo integrado $A(u)$ para el problema de sedimentaci\'on continua.}
\label{fig:difu3}
\end{center}
\end{figure}

En este caso, $\Delta\varrho=1690\ [Kg\, m^{-3}]$.
Adem\'as
\begin{equation}
a(u)=\left\{\begin{array}{ll}
0,&\textrm{ si }u\leqslant u_c:=0.1,\\
6.1267\times10^{2}\cdot u^8\left(1-\frac{u}{0.3}\right)^{5.647},&\textrm{ si }u> u_c:=0.1,\end{array}\right.
\end{equation}
y para el t\'ermino difusivo integrado $A(u)$ se utiliza (\ref{dif1}), (\ref{dif2}). Su gr\'afica se muestra en la figura \ref{fig:difu3}.

En la tabla \ref{tabla:caso3} se muestran la proporci\'on $V$, tasa de compresi\'on y errores entre la soluci\'on obtenida utilizando multiresoluci\'on y la soluci\'on obtenida sin multiresoluci\'on.

\begin{table}[h]
\begin{center}
{\small
\begin{tabular}{lccccc}
\hline
$t\, [s]$&  $V$   & $\mu$  &            $e_1$   &        $e_2$     &     $e_\infty$     \\
\hline
1800     & 6.6818  & 16.0156 & 7.81$\times10^{-5}$&5.83$\times10^{-5}$&1.80$\times10^{-5}$\\
3600  (*)& 7.0845  & 16.0156 & 1.61$\times10^{-4}$&6.77$\times10^{-5}$&4.01$\times10^{-5}$\\
7200  (*)& 7.6731  & 15.3010 & 2.44$\times10^{-4}$&9.05$\times10^{-5}$&6.46$\times10^{-5}$\\
14400 (*)& 9.5790  & 14.6441 & 4.92$\times10^{-4}$&1.64$\times10^{-4}$&1.84$\times10^{-4}$\\
43200 (*)& 14.0489 & 19.6441 & 5.10$\times10^{-4}$&4.26$\times10^{-4}$&4.76$\times10^{-4}$\\
\hline
\end{tabular}
}
\end{center}
\caption{\small Simulaci\'on de sedimentaci\'on continua. Tolerancia
prescrita $\varepsilon=5\times10^{-4}$, $N_0=513$ puntos en la
malla fina y $L=5$ niveles de multiresoluci\'on. (*): figuras \ref{fig:caso3_1hora} -
\ref{fig:caso3_12horas}.}
\label{tabla:caso3}
\end{table}

Al mirar la tabla \ref{tabla:caso3}, de nuevo los errores entre la soluci\'on obtenida utilizando multiresoluci\'on y la soluci\'on obtenida sin multiresoluci\'on, se encuentran por debajo de la tolerancia prescrita. Una alta tasa de compresi\'on y proporci\'on $V$ de tiempo total de CPU delatan la importancia del m\'etodo de multiresoluci\'on en la aplicaci\'on de este tipo de problemas.


\begin{figure}[!ht]
\begin{center}
\includegraphics[width=6.2in,height=2.5in]{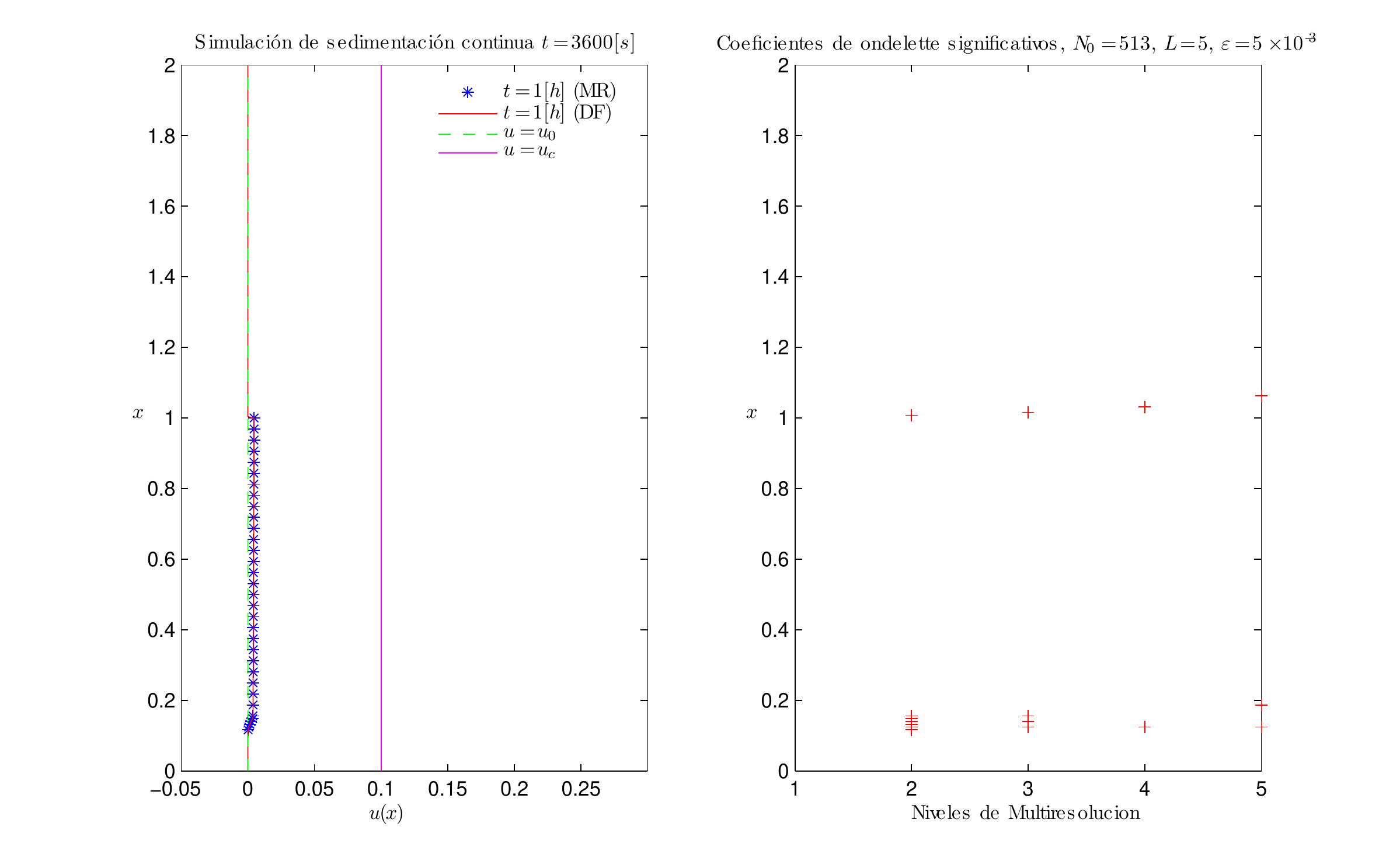}
\caption[Perfil de concentraci\'on. Sedimentaci\'on continua,
$t=3600$.]{\small Izquierda: Condici\'on inicial \emph{(rayas)} y
perfil de concentraci\'on a $t=1[h]$ para el problema de
sedimentaci\'on continua \emph{(asteriscos)}. Derecha:
Coeficientes de ondelette significativos correspondientes. $\varepsilon=5\times10^{-4}$, $N_0=513$ y $L=5$.}
\label{fig:caso3_1hora}
\end{center}
\end{figure}

\begin{figure}[!h]
\begin{center}
\includegraphics[width=6.2in,height=2.5in]{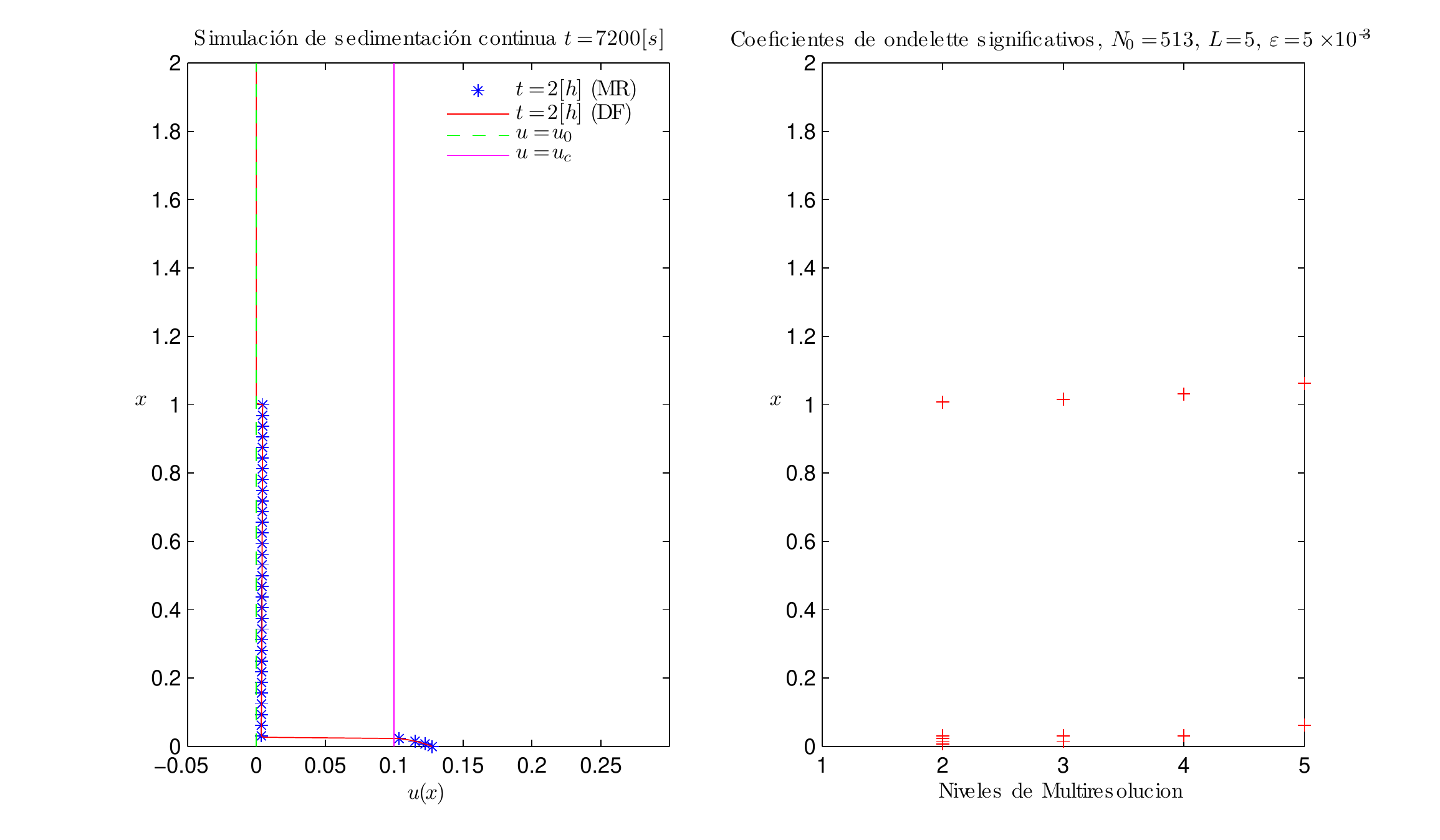}
\caption[Perfil de concentraci\'on. Sedimentaci\'on continua, $t=7200$.]{\small Izquierda: Condici\'on inicial \emph{(rayas)} y perfil de concentraci\'on a $t=2[h]$ para el problema de sedimentaci\'on continua \emph{(asteriscos)}. Derecha: Coeficientes de ondelette significativos correspondientes. $\varepsilon=5\times10^{-4}$, $N_0=513$ y $L=5$.}
\label{fig:caso3_2horas}
\end{center}
\end{figure}

\begin{figure}[!ht]
\begin{center}
\includegraphics[width=6.2in,height=2.7in]{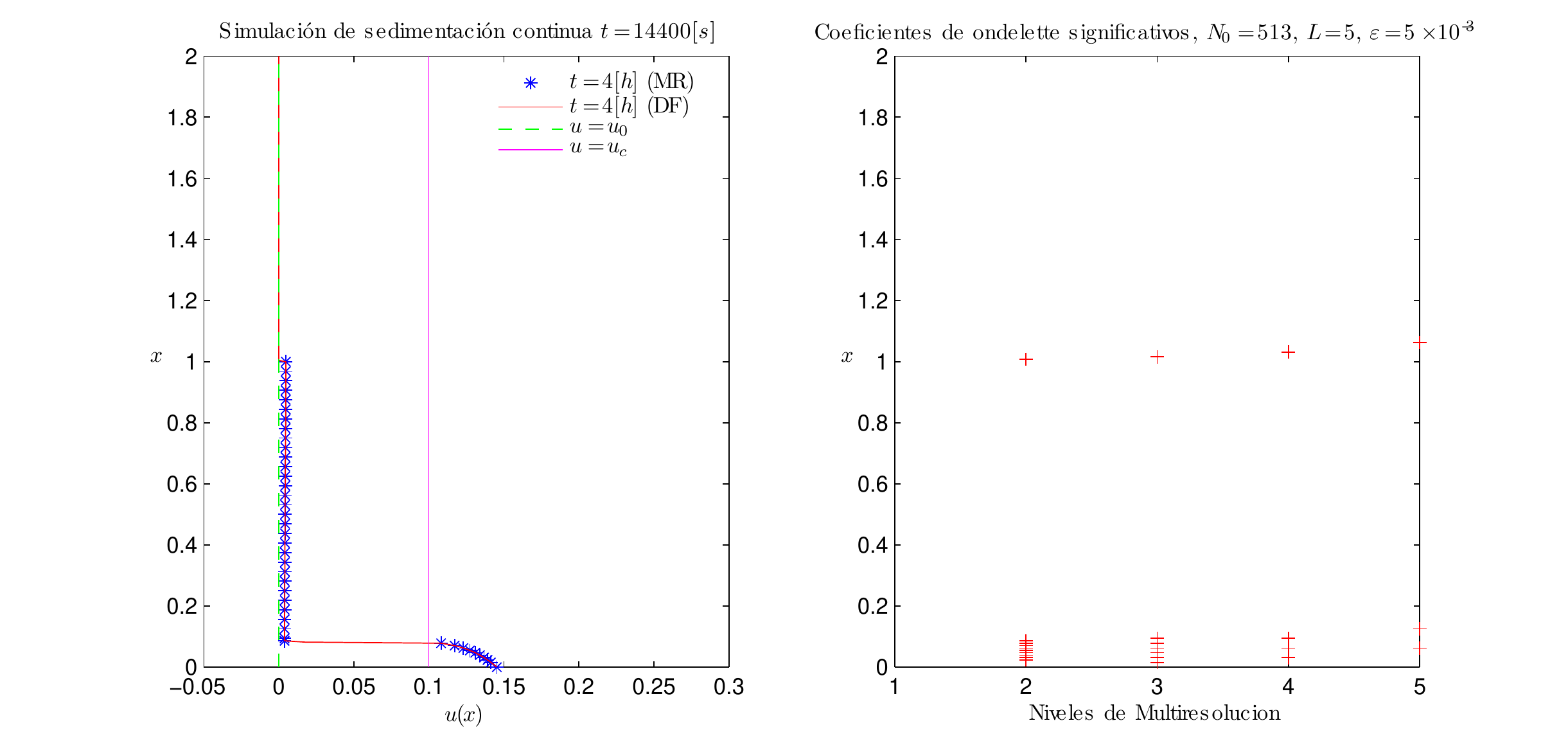}
\caption[Perfil de concentraci\'on. Sedimentaci\'on continua, $t=14400$.]{\small Izquierda: Condici\'on inicial \emph{(rayas)} y perfil de concentraci\'on a $t=4[h]$ para el problema de sedimentaci\'on continua \emph{(asteriscos)}. Derecha: Coeficientes de ondelette significativos correspondientes. $\varepsilon=5\times10^{-4}$, $N_0=513$ y $L=5$.}
\label{fig:caso3_4horas}
\end{center}
\end{figure}

\begin{figure}[!h]
\begin{center}
\includegraphics[width=6.2in,height=2.7in]{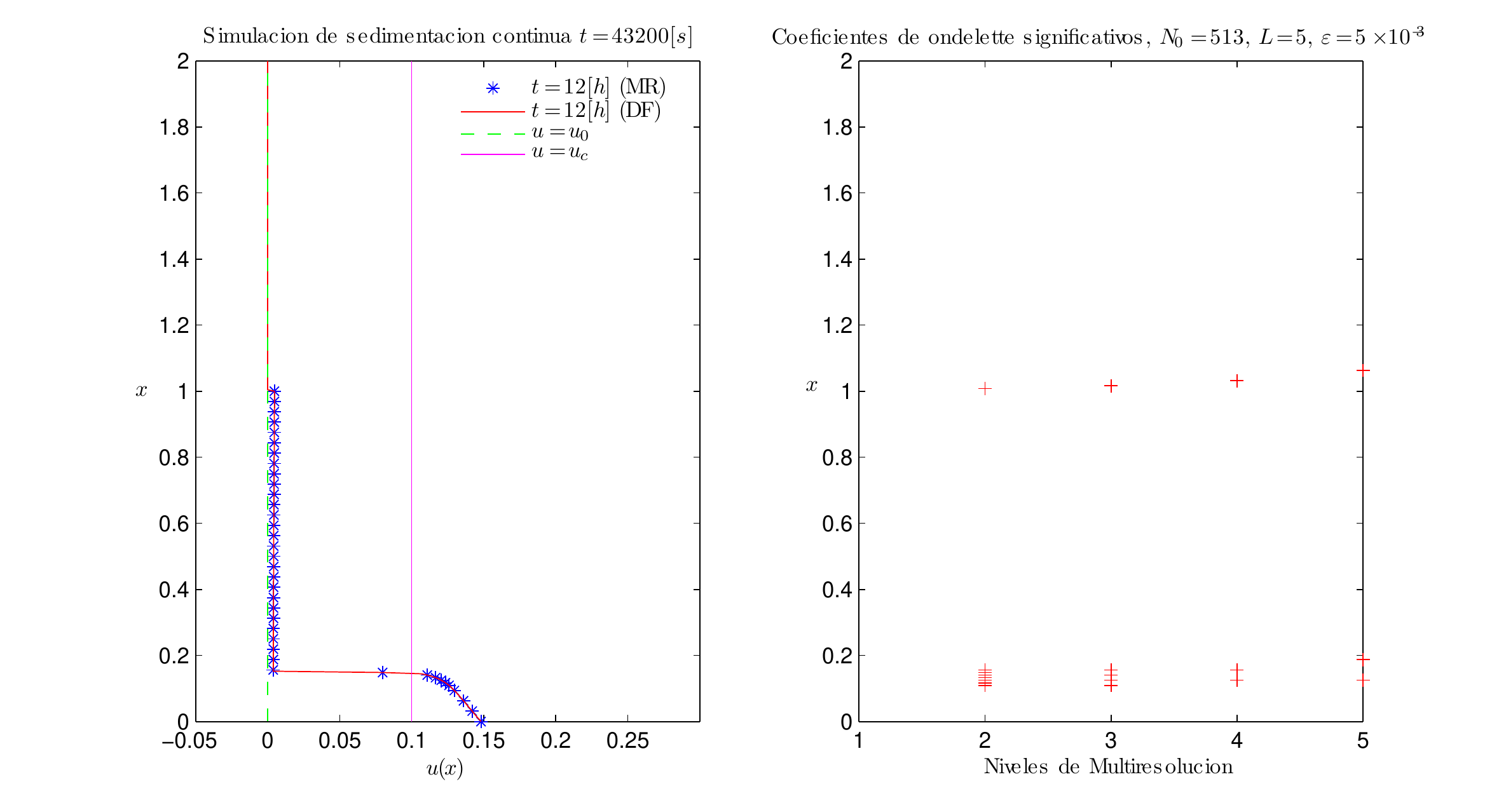}
\caption[Perfil de concentraci\'on. Sedimentaci\'on continua, $t=43200$.]{\small Izquierda: Condici\'on inicial \emph{(rayas)} y perfil de concentraci\'on a $t=12[h]$ para el problema de sedimentaci\'on continua \emph{(asteriscos)}. Derecha: Coeficientes de ondelette significativos correspondientes. $\varepsilon=5\times10^{-4}$, $N_0=513$ y $L=5$.}
\label{fig:caso3_12horas}
\end{center}
\end{figure}

En la figura \ref{fig:caso3_12horas} se presenta un perfil de concentraci\'on para el modelo de sedimentaci\'on continua, a $t=43200 [s]$. Notar que en este tiempo la soluci\'on ya entra en un estado estacionario, pues el flujo de alimentaci\'on es igual al flujo de descarga.

\begin{figure}[!ht]
\begin{center}
\includegraphics[width=5in,height=3.5in]{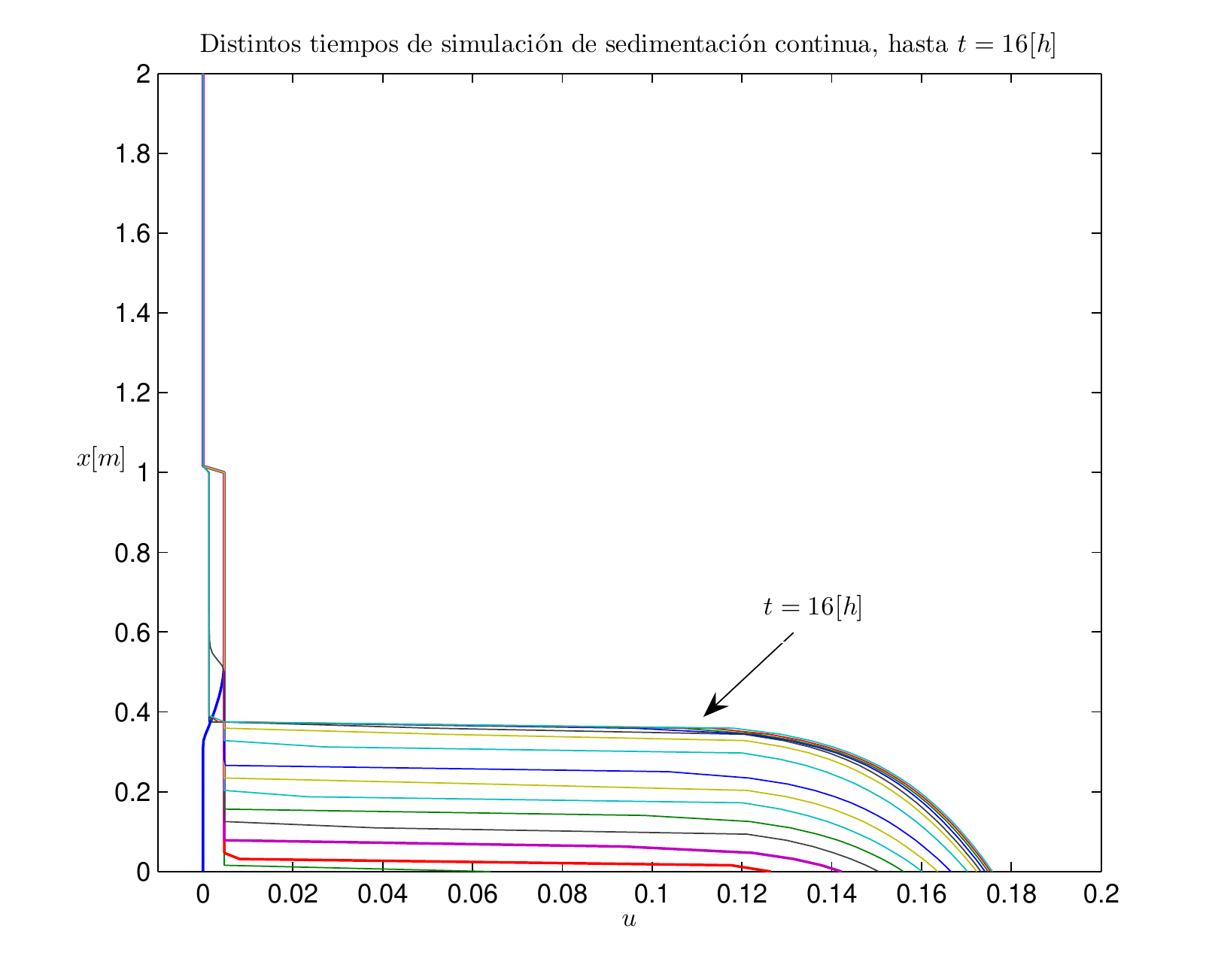}
\caption[Perfiles de concentraci\'on. Sedimentaci\'on continua
hasta $t=57600$.]{\small Perfiles de concentraci\'on hasta
$t=16[h]$ para el problema de sedimentaci\'on continua. $\varepsilon=5\times10^{-4}$, $N_0=513$ y $L=5$.}
\label{fig:caso3_plot}
\end{center}
\end{figure}

Finalmente se presenta en la figura \ref{fig:caso3_plot} la soluci\'on num\'erica del problema de sedimentaci\'on continua, obtenida utilizando el m\'etodo de multiresoluci\'on, hasta el tiempo $t=16[h]$.

Los resultados num\'ericos concuerdan con los resultados obtenidos por  B\"urger y Karlsen \cite{BK}.

\clearpage{\pagestyle{empty}\cleardoublepage}
\chapter{Conclusiones y perspectivas}
\section{Conclusiones}
En el presente trabajo se desarroll\'o un esquema num\'erico completamente adaptativo para acelerar los c\'alculos de vol\'umenes finitos de ecuaciones diferenciales parab\'olicas (originalmente desarrollado para leyes de conservaci\'on hiperb\'olicas) y ecuaciones parab\'olicas fuertemente degeneradas en una dimensi\'on espacial. Se estudiaron varios casos test de ecuaciones hiperb\'olicas,  parab\'olicas linales y no lineales, y ecuaciones parab\'olicas fuertemente degeneradas provenientes de la teor{\'i}a de procesos de sedimentaci\'on-consolidaci\'on. 

Generalmente, al a\~nadir un t\'ermino viscoso a un esquema, la soluci\'on tiende a suavizar y en algunos casos puede estabilizar un esquema num\'erico originalmente inestable. Se pudo ver que excepto por una limitaci\'on de paso temporal (que en el caso inv{\'i}scido es diferente) el problema viscoso no implica mayores complicaciones desde el punto de vista num\'erico.

El an\'alisis de multiresoluci\'on se mantiene inalterado, pues s\'olo tiene que ver con la regularidad de los valores puntuales o medias en celda de la soluci\'on.

Es importante destacar que en el cap{\'i}tulo 5 se utilizaron esquemas de diferencias finitas, por lo que en los algoritmos de multiresoluci\'on empleados se considera un an\'alisis de multiresoluci\'on para valores puntuales. 
 
Se comienza con una discretizaci\'on de vol\'umenes finitos (o diferencias finitas) en una malla uniforme, y una integraci\'on expl{\'i}cita en tiempo, ambas de segundo orden. Mediante t\'ecnicas de an\'alisis de multiresoluci\'on, se reduce el tama\~no de la malla, eliminando los puntos con detalles no significativos, pero manteniendo siempre un esquema de segundo orden.

La actualizaci\'on temporal de la malla se realiza mediante una estrategia de adaptaci\'on din\'amica que aprovecha la representaci\'on puntual esparsa, agregando coeficientes vecinos en escala y espacio para mejorar la captura de la informaci\'on. 

Para la evaluaci\'on de los flujos num\'ericos, en la malla localmente refinada, se utilizaron esquemas ENO de segundo orden y esquemas de Engquist-Osher modificados de segundo orden. 

Los algoritmos empleados son generalizables al caso de otras condiciones de borde (simplemente modificando el interpolador intermallas y el c\'alculo de los flujos en los puntos de frontera), otra elecci\'on para la condici\'on inicial, otro tipo de evoluci\'on temporal, otra elecci\'on para los predictores intermallas, otra elecci\'on para el orden de las interpolaciones ENO, otra elecci\'on para el c\'alculo del flujo num\'erico, otro tipo de estructura de datos, etc. 


La eficiencia del algoritmo fue medida mediante la tasa de compresi\'on y el tiempo de CPU. La diferencia de tiempo total de CPU entre la soluci\'on num\'erica que no utiliza multiresoluci\'on y la que utiliza multiresoluci\'on est\'a directamente relacionada con el hecho de que en una, la soluci\'on num\'erica sin multiresoluci\'on se evaluan todos los flujos num\'ericos mientras que en la otra soluci\'on num\'erica con multiresoluci\'on, s\'olo se calculan los flujos num\'ericos donde existen coeficientes de ondelette significativos. L\'ogicamente esta diferencia se ve incrementada cuando el flujo num\'erico es m\'as costoso. 

La aplicaci\'on del m\'etodo de multiresoluci\'on resulta a\'un m\'as provechosa en la simulaci\'on de procesos de sedimentaci\'on de suspensiones floculadas. El que las ecuaciones sean de naturaleza m\'as compleja, se suma el hecho de que los resultados experimentales publicados requieren un tiempo de simulaci\'on de varias \emph{horas}, en contraste con las fracciones de segundo suficientes para estudiar la soluci\'on num\'erica de los problemas hiperb\'olicos y parab\'olicos incluidos en este trabajo. Adem\'as, la condici\'on $CFL$ en este caso, hace que $\Delta t$ sea muy peque\~no. Esto hace pensar en la utilizaci\'on de un esquema impl{\'i}cito o semi- impl{\'i}cito \cite{BCS2}.   

La gran desventaja de utilizar algoritmos de multiresoluci\'on, es quiz\'as el hecho de que los resultados en cuanto a convergencia a\'un no tienen un gran auge. Una gran parte de los argumentos del an\'alisis de multiresoluci\'on desarrollado por Harten es de naturaleza heur{\'i}stica.   

En la parte final se present\'o un m\'etodo num\'erico para obtener soluciones aproximadas de problemas provenientes de fen\'omenos de sedimentaci\'on. La idea desarrollada fue aplicar los m\'etodos de multiresoluci\'on a los esquemas dise\~nados por B\"urger \emph{et al.} \cite{BCS,BEK,BEKL,BK,BWC} y se observ\'o que el m\'etodo de multiresoluci\'on es de gran ayuda para reducir el costo computacional en este tipo de problemas sin afectar la calidad de la soluci\'on.  

Todos los experimentos se realizaron en equipos con procesadores Pentium 4 de 1.6 Mhz, con 1GB de memoria RAM, tanto en plataforma Linux como Windows.

\section{Perspectivas}
\begin{itemize}
\item Para el caso de ecuaciones parab\'olicas fuertemente degeneradas, la perspectiva a m\'as corto plazo es modificar el algoritmo para poder aplicarlo a las ecuaciones que modelan otros tipos de fen\'omenos de sedimentaci\'on. 
\item Aplicar m\'etodos de multiresoluci\'on a la resoluci\'on de problemas inversos. 
\item Utilizar esquemas ENO de orden superior a dos. Combinar esto con la utilizaci\'on de esquemas con varios \emph{switches} \cite{Harten4}.     
\item Aplicar m\'etodos de multiresoluci\'on a problemas que modelan la separaci\'on de suspensiones polidispersas \cite{BBK}.
\item Realizar los experimentos del cap{\'i}tulo 5 utilizando esquemas semi-impl{\'i}citos e impl{\'i}citos. Esto se traduce en pasar de un $\Delta t$ de orden de $(\Delta x)^2$ a un orden de $\Delta x$. Sin embargo las complicaciones est\'an en tener que resolver un sistema de ecuaciones no-lineales en cada iteraci\'on. Adem\'as el proceso de multiresoluci\'on para esquemas semi-impl{\'i}citos se complica bastante.
\item Extender los resultados de los puntos anteriores al caso de sistemas y ecuaciones multidimensionales.
\item Los c\'odigos pueden ser f\'acilmente traducidos a un lenguaje m\'as robusto como FORTRAN, C, o C++, dado que las funciones y subrutinas en la implementaci\'on no abusan de las funciones impl{\'i}citas de MATLAB (excepto en la estructura SPARSE de los datos).
\end{itemize}

\clearpage{\pagestyle{empty}\cleardoublepage}
\begin{appendix}
\chapter{C\'alculo de los coeficientes de interpolaci\'on en la multiresoluci\'on}
\section{Multiresoluci\'on de valores puntuales}\label{ap1}
En la secci\'on correspondiente se ha mencionado que 
\begin{equation}\label{xxa}
\mathcal{I}(x_{2j-1}^{k-1},u^k)=\sum_{l=1}^s\beta_l(u_{j+l-1}^k+u_{j-l}^k)
\end{equation}
es el polinomio de grado $r-1$ que interpola los puntos $(u_{j-s}^k,\ldots,u_{j+s-1}^k)$. Para ver esto, y encontrar los valores de los coeficientes $\beta_l$, se utiliza el polinomio interpolador de Lagrange
\begin{equation*}
P(x)=\sum_{l=j-s}^{j+s-1}u(x_l^k)\prod_{l=j-s}^{j+s-1}\frac{x-x_i^k}{x_l^k-x_i^k},\quad i\neq l
\end{equation*}
en el punto $x_{2j-1}^{k-1}$
\begin{equation*}
\mathcal{I}(x_{2j-1}^{k-1},u^k)=\sum_{l=j-s}^{j+s-1}u_l^k\prod_{l=j-s}^{j+s-1}\frac{x_{2j-1}^{k-1}-x_i^k}{x_l^k-x_i^k},\quad i\neq l,
\end{equation*}
donde $x_{2j-1}^{k-1}=(2j-1)\cdot h_{k-1}=(j-1/2)\cdot h_k$. Luego
\begin{equation*}
\mathcal{I}(x_{2j-1}^{k-1},u^k)=\sum_{l=j-s}^{j+s-1}u_l^k\prod_{l=j-s}^{j+s-1}\frac{j-\frac{1}{2}-i}{l-i},\quad i\neq l.
\end{equation*}
Si se toma en cuenta que los pares de valores puntuales $(u_{j-1},u_j),\ (u_{j-2}u_{j+1}),\ldots$ est\'an multiplicados por el mismo factor, se tiene (\ref{xxa}), con
\begin{equation*}
\beta_l=\frac{1^2\cdot3^2\cdots(2l-3)^2\cdot(2l-1)\cdot(2l+1)^2\cdots(2s-1)^2}{2^{2s-1}\cdot(s+l-1)!\cdot(s-l)!}\cdot(-1)^{l+1}.
\end{equation*}
Por lo tanto se tienen los siguientes coeficientes para cada $r=2s$ mencionado:
\begin{itemize}
\item $r=2,\ s=1$, $\beta_1=\frac{1}{2\cdot1!\cdot0!}\cdot(-1)^2=\frac{1}{2}$.
\item $r=4,\ s=2$, $\beta_1=\frac{1\cdot3^2}{2^3\cdot2!\cdot1!}\cdot(-1)^2=\frac{9}{16}$, $\beta_2=\frac{1^2\cdot3}{2^3\cdot3!\cdot0!}\cdot(-1)^3=\frac{-1}{16}$.
\end{itemize}
\section{Multiresoluci\'on de medias en celda}\label{ap2}
De manera an\'aloga al caso anterior, en la secci\'on correspondiente se ha mencionado que 
\begin{equation}\label{xxb}
d_j^k=\bar{u}^{k-1}_{2j-1}-\tilde{u}^{k-1}_{2j-1}=\bar{u}^{k-1}_{2j-1}-\frac{\mathcal{I}(x^{k-1}_{2j-1},U^k)-U^k_{j-1}}{h_{k-1}}
\end{equation}
Se aplica el caso anterior (para valores puntuales) al est\'encil $(U^k_{j-s},\ldots,U^k_{j+s-1})$, por tanto
\begin{equation*}
\mathcal{I}(x_{2j-1}^{k-1};U^k)=\sum_{l=1}^s\beta_l(U^k_{j+l-1}+U^k_{j-l}),
\end{equation*}
con los mismos $\beta_l$ calculados en el ap\'endice \ref{ap1}, por tanto
\begin{equation*}
\tilde{u}^{k-1}_{2j-1}=\frac{\sum_{l=1}^s\beta_l(U^k_{j+l-1}+U^k_{j-l})-U^k_{j-1}}{h_{k-1}},
\end{equation*}
y utilizando $h_k=2h_{k-1}$, se tiene
\begin{eqnarray*}
\tilde{u}^{k-1}_{2j-1}&=&\frac{\sum_{l=1}^s\beta_l(U^k_{j+l-1}+U^k_{j-l})-U^k_{j-1}}{2^{-1}h_k}\\
&=&\frac{\sum_{l=1}^s2\beta_l(U^k_{j+l-1}+U^k_{j-l})-2U^k_{j-1}}{h_k}\\
&=&\frac{2\beta_1(U^k_j+U^k_{j-1})+2\beta_2(U^k_{j+1}+U^k_{j-2})+\cdots+2\beta_s(U^k_{j+s-1}+U^k_{j-s})-2U^k_{j-1}}{h_k},
\end{eqnarray*}
o equivalentemente,
\begin{eqnarray*}
h_k\cdot\tilde{u}^{k-1}_{2j-1}&=&\cdots+(2\beta_1-1+2\beta_2)\cdot(U^k_{j-2}-U^k_{j-3})+(2\beta_1-1)\cdot(U^k_{j-1}-U^k_{j-2})+\\
&&1\cdot(U^k_j-U^k_{j-1})+(1-2\beta_1)\cdot(U^k_{j+1}-U^k_j)+\\
&&(1-2\beta_1-2\beta_2)\cdot(U^k_{j+2}-U^k_{j+1})+\cdots
\end{eqnarray*}
Adem\'as, si se toma en cuenta la relaci\'on 
$$\bar{u}_j^k=\frac{U^k_j-U^k_{j+1}}{h_k},$$
y que las medias en celda $\bar{u}_{j+l}^k$ y $\bar{u}_{j+l}^k$, $l=1,\ldots,s-1$, est\'an multiplicados por el mismo factor (s\'olo cambia de signo), se llega a la expresi\'on
\begin{equation}
\tilde{u}_j^k=\sum_{l=1}^{s-1}\gamma_l(\bar{u}_{j+l}^k-\bar{u}_{j-l}^k)+\bar{u}_j^k
\end{equation}
con
$$\gamma_l=-(2\cdot\beta_l-\gamma_{l-1}),\quad \gamma_0=1.$$
Por lo tanto se obtiene para cada $\bar{r}=2s-1$ mencionado:
\begin{itemize}
\item $\bar{r}=3,\ s=2$, 

$\gamma_1=-(2\cdot\beta_1-\gamma_0)=-(2\cdot\frac{9}{16}-1)=-\frac{1}{8}$.
\item $\bar{r}=5,\ s=3$, 

$\gamma_1=-(2\cdot\beta_1-\gamma_0)=-(2\cdot\frac{150}{256}-1)=-\frac{22}{128}$, 

$\gamma_2=-(2\cdot\beta_2-\gamma_1)=-(2\cdot\frac{-25}{256}+\frac{22}{128})=\frac{3}{128}$.
\end{itemize}

\clearpage{\pagestyle{empty}\cleardoublepage}
\chapter{An\'alisis de estabilidad para el caso parab\'olico no lineal}\label{anexob}
A continuaci\'on se analizar\'a la estabilidad en el sentido de la variaci\'on total del esquema num\'erico presentado en la secci\'on \ref{estab1}. Este an\'alisis puede aplicarse al caso de flujos lineales o no lineales. Se quiere encontrar una condici\'on CFL que pueda utilizarse para el esquema ENO-TVD de segundo orden.

Un esquema se dice \emph{TV-estable} si la variaci\'on total 
$$TV(v_h(x,t))=TV(v^n):=\sum_{j=0}^{N-1}|v^n_{j+1}-v^n_j|$$
de una sucesi\'on de aproximaciones num\'ericas $v_h(x,t)$ est\'a acotada uniformemente en $h=\Delta x$ y $t=n\Delta t$, con $h\to 0$ y $0\leqslant t\leqslant T$. A\'un m\'as, el esquema es TVD si
$$TV(v^{n+1})\leqslant TV(v^n).$$
Claramente, un esquema TVD es TV-estable.

Con estas definiciones b\'asicas, Harten \cite{Harten3} prob\'o el siguiente
\begin{Le} Si un esquema escrito en la forma 
\begin{equation}\label{leb}
v_j^{n+1}=v_j^n+C_j^+\Delta_+v_j^n-C_{j-1}^-\Delta_-v_j^n,
\end{equation}
satisface, para todo $j$,
\begin{eqnarray}
C_j^+&\geqslant&0,\label{lemaa}\\
C_j^-&\geqslant&0,\label{lemab}\\
C_j^++C_j^-&\leqslant& 1\label{lemac},
\end{eqnarray}
entonces el esquema es TVD.
\end{Le}
\'El introduce el siguiente esquema expl{\'i}cito, de primer orden,
\begin{eqnarray}
v_j^{n+1}&=&=v_j^n\hat{C}_j^+\Delta_+v_j^n-\hat{C}_{j-1}^-\Delta_-v_j^n,\\
\hat{C}_j^\pm &=&\frac{1}{2}[|\omega+\xi|\mp(\omega+\xi)]_j\label{hatC}\\
\omega_j&=&\lambda\frac{\Delta_+f_j}{\Delta_-v_j}\\
\xi_j&=&\frac{\Delta_+g_j}{\Delta_-v_j},
\end{eqnarray}
donde $f_j=f(u_j)$, y $g_j=g(u_j)$ es elegida tal que
\begin{equation}\label{propiedadx}
|\xi_j|\leqslant \rho(\omega_j)
\end{equation}
donde $\rho$ es el cl\'asico limitador de flujo de Harten
\begin{equation*}
\rho(a)=\left\{\begin{array}{ll}
0,&\textrm{ para primer orden en espacio}\\
\frac{1}{2}(|a|-a^2),&\textrm{ para segundo orden en espacio.}
\end{array}\right.
\end{equation*}
Con estas definiciones, Harten prueba que para esquemas de primer y segundo orden, una condici\'on suficiente para que el esquema sea TVD es la condici\'on tipo $CFL$
\begin{equation}
\max_j|\omega_j|\leqslant 1,
\end{equation}
pues $\omega_j$ es el coeficiente CFL medio local. 

Se quiere modificar la demostraci\'on hecha por Harten \cite{Harten3} para el caso de esquemas de segundo orden, con el fin de aplicarla al caso viscoso, para ello, Bihari \cite{Bihari} prob\'o el siguiente
\begin{Th}
Un esquema escrito en la forma (\ref{leb}), con $C_j^\pm$ definido por
\begin{equation}\label{leb2}
C_j^\pm=\hat{C}_j^\pm+\lambda\frac{\nu}{\Delta x},
\end{equation}
es TVD si 
\begin{equation}
\sigma\leqslant \frac{Re}{Re+4}
\end{equation}
con 
\begin{equation}\label{defis}
\sigma=\max_j|\omega_j|,\quad Re=\max_j|\omega_j|\frac{\Delta x}{\lambda\nu}.
\end{equation}
\end{Th}
Notar que el esquema (\ref{leb}), (\ref{leb2}) es de segundo orden en espacio al aproximar la soluci\'on del problema, puesto que se ha incluido un t\'ermino viscoso con una discretizaci\'on central al esquema original TVD de segundo orden.

Notar adem\'as que la definici\'on dada para $\sigma$ y $Re$ difieren de las definiciones dadas en (\ref{defisig}) y (\ref{defire}) para el caso lineal. Sin embargo el significado cualitativo de estas cantidades es el mismo. Es decir, (\ref{defis}) es la definici\'on equivalente para $\sigma$ y $Re$ en el caso no lineal.

\dem (Del teorema) Se mostrar\'a que se satisfacen las condiciones del Lema 1. Con la definici\'on dada de $\hat{C}_j^\pm$ (\ref{hatC}), es claro que se satisfacen las condiciones (\ref{lemaa}) y (\ref{lemab}). Falta entonces mostrar que
\begin{equation}\label{propiedady}
\hat{C}_j^++\hat{C}_j^-+2\lambda\frac{\nu}{\Delta x}\leqslant 1.
\end{equation}
De las definiciones dadas y de la propiedad (\ref{propiedadx}) se sigue que (\ref{propiedady}) se satisfar\'a si
\begin{equation}\label{propiedadz}
\frac{3}{2}\sigma-\frac{1}{2}\sigma^2+2\sigma\frac{1}{Re}\leqslant 1.
\end{equation}
Ahora, dado que $\sigma\leqslant 1$ (necesario para que se satisfaga (\ref{propiedadx})), se tiene que $\sigma\frac{4}{Re}(\sigma-1)\leqslant 0$. Luego, es posible obtener una versi\'on levemente m\'as restrictiva que (\ref{propiedadz}):
\begin{equation*}
(\sigma-2)\left[\sigma\left(1+\frac{4}{Re}\right)-1\right]\geqslant 0,
\end{equation*}
la cual se satisface si se satisface 
\begin{equation*}
\sigma\leqslant \frac{Re}{Re+4}.
\end{equation*}

\qedo
\clearpage{\pagestyle{empty}\cleardoublepage}
\chapter{C\'odigo y documentaci\'on}

Tanto los c\'odigos en MATLAB para cada experimento, la documentaci\'on respectiva, como una versi\'on electr\'onica de este informe pueden ser obtenidos en forma gratuita, desde el sitio \href{http://www.udec.cl/~riruiz/tesis.html}{http://www.udec.cl/$\sim$riruiz/tesis.html}.
\end{appendix}
\clearpage{\pagestyle{empty}\cleardoublepage}
\listoffigures
\listoftables

\end{document}